\documentclass[11pt]{amsart}
\usepackage[T1]{fontenc}
\usepackage[latin9]{inputenc}
\usepackage{amsmath}
\usepackage{tikz}
\usetikzlibrary{arrows.meta,decorations.pathreplacing,calc}
\usepackage{amsthm}
\usepackage{amssymb}
\usepackage{setspace}

\makeatletter
\theoremstyle{plain}
\newtheorem{thm}{\protect\theoremname}[section]
\theoremstyle{plain}
\newtheorem{lem}[thm]{\protect\lemmaname}
\theoremstyle{plain}
\newtheorem{prop}[thm]{\protect\propositionname}
\theoremstyle{definition}
\newtheorem{remark}[thm]{\protect\remarkname}
\theoremstyle{plain}

\theoremstyle{plain}
\newtheorem{corollary}[thm]{\protect\corollaryname}

\newtheorem{mainthm}{Theorem}

\providecommand{\corollaryname}{Corollary}
\newtheorem{maincorollary}{Corollary}

\usepackage{hyperref}
\usepackage{a4wide}
\@ifundefined{definecolor}
 {\usepackage{color}}{}

\newcommand{\Z}{\mathbb{Z}}

\makeatother

\usepackage{babel}
\providecommand{\lemmaname}{Lemma}
\providecommand{\theoremname}{Theorem}
\providecommand{\propositionname}{Proposition}
\providecommand{\remarkname}{Remark}
\providecommand{\examplename}{Example}

\begin{document}
\raggedbottom

\title[Spectral zeta values of metric graphs]{An evaluation of spectral zeta functions of metric graphs at any integer argument}
\author{Jay Jorgenson}\thanks{The first-named author acknowledges grant support from PSC-CUNY Awards 67415-00 55 and 68462-00 56, which are jointly funded by the Professional Staff Congress and The City University of New York.}
\author{Anders Karlsson}\thanks{The second-named author was supported by the Swiss NSF Grants 200020-200400 and 200021-212864, and by the Swedish Research Council Grant 104651320.}
\author{Lejla Smajlovi\'{c}}
\date{\today}
\begin{abstract}\noindent
Let $\Gamma$ be any finite metric graph, possibly with loops and multiple edges, with underlying
combinatorial graph $G_{0}$, and let $\zeta_{\Gamma}$ be the spectral zeta function of the Laplacian on
$\Gamma$ with standard vertex conditions.  Then for every positive integer $n$ we prove an exact, closed
form evaluation of $\zeta_{\Gamma}(n)$ as a ratio of weighted sums over the spanning forests of $G_{0}$ with at
most $n+1$ components, where the weights are products of edge lengths with coefficients built from
Bernoulli numbers.  The evaluation comes from a limit theorem, in which the suitably normalized
characteristic polynomials of the graphs obtained by subdividing the edges of $\Gamma$ converge to the
zeta-regularized determinant $\text{\rm det}(\Delta_{\Gamma}+x)$.  In the case $\Gamma$ is a dipole
graph with two edges of equal length, the evaluation reduces to Euler's computation of the special values
of the Riemann zeta function at the even integers.  In recent work, Kurasov and Sarnak showed that, for
rationally independent edge lengths, the spectrum of the Laplacian on $\Gamma$ spans an
infinite-dimensional vector space over $\mathbb{Q}$, and they expect, conditionally on Schanuel's
conjecture, that the eigenvalues are algebraically independent, hence that the field they generate has
infinite transcendence degree.  Nonetheless, we derive a precise formula for
$\zeta_{\Gamma}(n)$ as the value of a homogeneous rational function of degree $2n$ in the edge lengths, with coefficients in $\mathbb{Q}$ determined by $G_{0}$.
\end{abstract}

\maketitle

\section{Introduction}\label{sec:intro}

For almost a century, metric graphs have appeared in chemistry, physics, evolutionary biology, and mathematics within topics such as spectral theory, arithmetic geometry, and geometric group theory \cite{Ku24}. From a geometric viewpoint, one can consider metric graphs as an interesting class of $1$-dimensional spaces, with combinatorial and topological complexity 
that closed $1$-dimensional manifolds, the circle, do not have. 
Like for Riemannian manifolds, to any metric graph there is a naturally associated Laplacian. For a finite, compact, connected metric graph $\Gamma$, the spectrum of this operator contains only eigenvalues that we denote $0 = \mu_{0} < \mu_{1}\leq\mu_{2}\leq\cdots\rightarrow\infty$, and the \emph{spectral zeta function}
$$
\zeta_{\Gamma}(s) = \sum_{j\geq 1}\frac{1}{\mu_{j}^{s}}.
$$
The special values at integers are of particular interest, for example in the study of polymers \cite{Y71}. For the use of spectral zeta functions in physics in the specific context of metric graphs, we refer to \cite{KoSm99, ACDMT00, Ho08}.


Our first main result determines these special values in terms of the underlying combinatorial graph $G_0$ and with the lengths $\ell_j$ of $\Gamma$.  Define the \emph{metric forest function}
\begin{equation}\label{eq:mforestfunction}
\Phi(x) = \sum_{F}\ \prod_{\mathcal{E}_{j}\notin F}\frac{\sinh(\sqrt{x}\,\ell_{j})}{\sqrt{x}}\ \prod_{C\in F}\Big(\sqrt{x}\sum_{h\text{ \rm at }C}\tanh\frac{\sqrt{x}\,\ell_{e(h)}}{2}\Big),
\end{equation}
where $F$ runs over all spanning forests of $G_0$, $C$ are the connected components of $F$, and the inner
sum is over the half-edges $h$ of $G_{0}$ at the vertices of $C$, with $e(h)$ the edge carrying $h$. Let the total length, or volume, $L$ be the finite sum of the edge lengths $\ell_{1}+\ell_{2}+\cdots+\ell_{E_{0}}$.

\begin{mainthm} \label{mthm:special}
Let $\Gamma$ be a finite, compact, connected metric graph with underlying combinatorial graph $G_0$ and edge
lengths $\ell_j$.  Write $\Phi(x) = \sum_{m\geq 1}P_{m}x^{m}$, where each $P_{m}$ is an explicit finite sum over
the spanning forests of $G_{0}$ with at most $m$ components, and put $e_{m} = P_{m+1}/P_{1}$.  Then $\Phi$ is
entire and the following hold.
\begin{enumerate}
\item $\zeta_{\Gamma}$ has a meromorphic continuation to $\mathbb{C}$, whose only pole is a simple pole at
$s = 1/2$ with residue $L/2\pi$.
\item At the non-positive integers,
$$
\zeta_{\Gamma}(0) = \frac{\chi(G_{0})}{2}-1,\qquad \zeta_{\Gamma}(-n) = 0\quad(n\geq 1).
$$
\end{enumerate}
At the positive integers, the following hold for every integer $n\geq 1$.
\begin{enumerate}
\item[(3)] (\emph{logarithmic derivative})
$$
\zeta_{\Gamma}(n) = \frac{(-1)^{n-1}}{(n-1)!}\,\frac{d^{\,n-1}}{dx^{n-1}}\Big(\frac{\Phi'(x)}{\Phi(x)}-\frac{1}{x}\Big)\Big|_{x=0};
$$
\item[(4)] (\emph{Newton's recursion})
$$
\zeta_{\Gamma}(n) = \sum_{i=1}^{n-1}(-1)^{i-1}e_{i}\,\zeta_{\Gamma}(n-i)+(-1)^{n-1}\,n\,e_{n};
$$
\item[(5)] (\emph{closed form})
$$
\zeta_{\Gamma}(n) = \sum_{(m_{1},\dots,m_{n})}(-1)^{n+m_{1}+\cdots+m_{n}}\,
\frac{n\,(m_{1}+\cdots+m_{n}-1)!}{m_{1}!\cdots m_{n}!}\prod_{i=1}^{n}e_{i}^{m_{i}},
$$
the sum being over all $(m_{1},\dots,m_{n})\in\mathbb{Z}_{\geq 0}^{n}$ with $m_{1}+2m_{2}+\cdots+nm_{n} = n$.
\end{enumerate}
\end{mainthm}

Note that in particular, $\zeta_{\Gamma}(n)$ is a rational function of the edge lengths with rational coefficients, determined by the spanning forests of $G_{0}$ with at most $n+1$ components. For example, for $n = 1$ the recursion gives $\zeta_{\Gamma}(1) = e_{1} = P_{2}/P_{1}$, and the formulas of
Section~\ref{sec:explicit} for $P_{1}$ and $P_{2}$ give an expression in terms of spanning trees and two-component forests.
Similarly, $\zeta_{\Gamma}(2) = e_{1}^{2}-2e_{2}$ involves forests with at most three components.

Theorem \ref{mthm:special} has the following arithmetic consequence.
\begin{maincorollary} \label{cor:arithmetic}
For every integer $n$,
$$
\zeta_{\Gamma}(n)\in \mathbb{Q}(\ell_{1},\dots,\ell_{E_{0}}),
$$
the field generated by the edge lengths.
Consequently, the following hold.
\begin{enumerate}
\item if all edge lengths are algebraic numbers, then every $\zeta_{\Gamma}(n)$ is an algebraic number;
\item if all edge lengths are rational multiples of a single length $a$, then $\zeta_{\Gamma}(n)\in\mathbb{Q}\,a^{2n}$ for every
integer $n$.
\end{enumerate}

\end{maincorollary}

Apart from spectra that are arithmetic progressions, where
the evaluation is Euler's, and the inverse power sums of the zeros of a Bessel function, evaluated by Rayleigh
and recalled below, we know of no setting with an infinite spectrum in which a spectral zeta function is
evaluated in closed form at every positive integer.  For the circle, Theorem~\ref{mthm:special} produces $\zeta_{\Gamma}(n)$ from the Taylor
coefficients of $\Phi$ alone, and comparing the result with
$\zeta_{\Gamma}(s) = 2(L/2\pi)^{2s}\zeta_{\mathbb{Q}}(2s)$ returns Euler's evaluation of
$\zeta_{\mathbb{Q}}(2n)$, where $\zeta_{\mathbb{Q}}$ denotes the Riemann zeta function; see
Section~\ref{subsec:circle}.

A classical
setting is provided by the $J$-Bessel functions.  The squares of the zeros $j_{n, k}$ of $J_{n}$ are the eigenvalues for the Dirichlet problem on the unit disk. Rayleigh \cite{Ra1874} evaluated the sums of the reciprocals; for instance, it is proved that
\begin{equation}\label{eq:Rayleigh}
\sum_{k=1}^{\infty} \frac{1}{j_{n ,k}^{10}}=\frac{2(7n+19)}{2^{10}(n+1)^5(n+2)^2(n+3)(n+4)(n+5)},
\end{equation}
which were then used to obtain precise estimates for the smallest non-zero eigenvalue; see \cite[15.51]{Wa66}.
 On the other hand, Siegel proved in a celebrated memoir \cite{Si29} that the Bessel zeros are transcendental numbers.  Therefore, \eqref{eq:Rayleigh} shows that the sum of reciprocal powers of the transcendental numbers $\{j_{n,k}\}$ is rational.
Furthermore, recent papers by Kurasov--Sarnak \cite{KS20, KS23} investigate arithmetic aspects of the spectrum of finite metric graphs, see also
\cite[Theorem~10.1 and (10.11)]{Ku24}. If the edge lengths are pairwise rationally dependent, then after
inserting vertices of degree two the graph becomes equilateral, and the set of frequencies $\sqrt{\mu_j}$ is a
finite union of arithmetic progressions. On the other hand, excluding the interval, circle and figure-eight graphs, the spectrum is highly algebraically independent
over $\mathbb{Q}$. Their proof uses the irreducibility of the reduced secular polynomial and diophantine analysis like Lang's conjecture on
torsion points of subvarieties of tori, which is a theorem; see \cite[Chapter~10]{Ku24} and
\cite{KS23}.  Most strikingly, assuming Schanuel's conjecture, it is shown in \cite{KS20} and \cite{KS23} that with the exception of the three graphs mentioned, whenever the lengths of the edges are algebraic numbers that are linearly independent over $\mathbb Q$, the spectrum has infinite transcendence degree. That is, there is an infinite set such that there are no polynomial relations between any finite subset.

For essentially all metric graphs, then, no closed formula for the individual eigenvalues is known or
expected, since they are the roots of transcendental secular equations.  Yet
Theorem~\ref{mthm:special} evaluates every power sum of these roots in closed form, and
Corollary~\ref{cor:arithmetic} places each one in the field generated by the edge lengths.  The
comparison is sharpest once the arithmetically transparent part of the spectrum is removed.  A loop of
length $\ell$ always carries the frequencies $2\pi n/\ell$, an arithmetic progression whose
contribution is Euler's.  Removing one such progression for each loop, and writing
$\text{\rm Spec}^{*}(\Gamma)$ for what is left and $J$ for the set of loops,
Theorem~\ref{thm:reduced} gives
\begin{equation}\label{eq:reduced-intro}
\sum_{\kappa\in\text{\rm Spec}^{*}(\Gamma)}\kappa^{-2r}
= \zeta_{\Gamma}(r)-\sum_{j\in J}\frac{(-1)^{r+1}B_{2r}}{2\,(2r)!}\,\ell_{j}^{2r}
\;\in\;\mathbb{Q}(\ell_{1},\dots,\ell_{E_{0}})
\end{equation}
for every integer $r\geq 1$.  When the edge lengths are rationally independent, and $\Gamma$ is not a
segment, a cycle or the figure-eight graph, $\text{\rm Spec}^{*}(\Gamma)$ is the reduced spectrum of
\cite{KS20} and \cite[Theorem~10.1]{Ku24}, which contains no arithmetic progression and spans an
infinite-dimensional vector space over $\mathbb{Q}$.  So a set of numbers carrying no arithmetic
structure, and expected on Schanuel's conjecture to be algebraically independent, has every one of its
inverse even power sums in the field generated by the edge lengths.

If $\Gamma$ has no loops, then
$J$ is empty and $\text{\rm Spec}^{*}(\Gamma)$ is the whole spectrum, so that
\eqref{eq:reduced-intro} is Corollary~\ref{cor:arithmetic} itself.  The subtraction is needed exactly when $\Gamma$ has loops,
since the progressions they carry are arithmetically transparent and would otherwise be suspected of
accounting for the rationality on their own.  In this we see
Theorem~\ref{mthm:special} and
Corollary~\ref{cor:arithmetic} as a vast generalization of Rayleigh's theorem to any finite metric
graph, and all of it is unconditional.

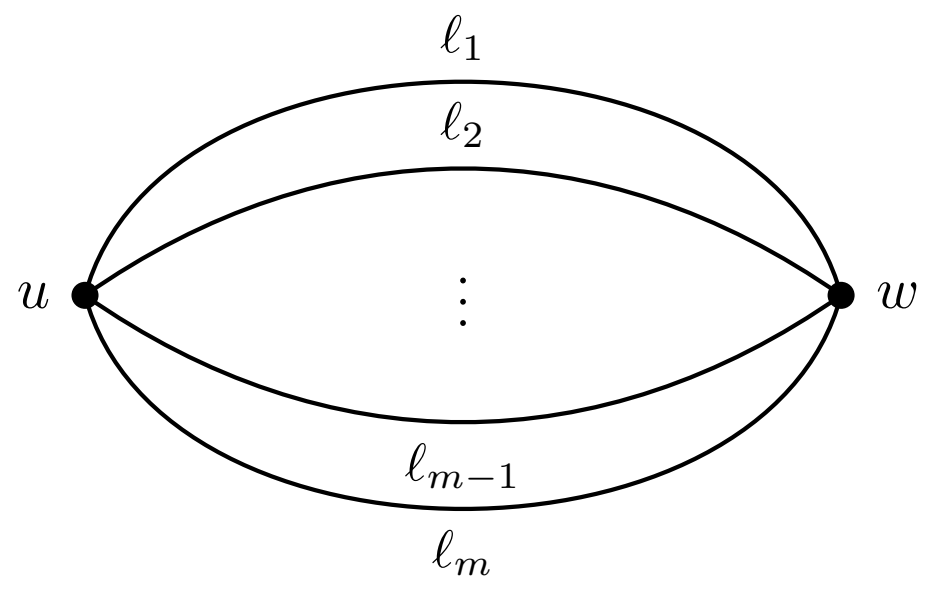
\begin{figure}[ht]
\centering
\begin{tikzpicture}[scale=1.25, v/.style={circle,fill,inner sep=1.8pt}, e/.style={thick}, lab/.style={font=\small}]
  \coordinate (U) at (-2,0); \coordinate (W) at (2,0);
  \draw[e] (U) to[bend left=75] node[lab,above] {$\ell_{1}$} (W);
  \draw[e] (U) to[bend left=35] node[lab,above] {$\ell_{2}$} (W);
  \node at (0,0.05) {$\vdots$};
  \draw[e] (U) to[bend right=35] node[lab,below] {$\ell_{m-1}$} (W);
  \draw[e] (U) to[bend right=75] node[lab,below] {$\ell_{m}$} (W);
  \node[v,label=left:{$u$}] at (U) {}; \node[v,label=right:{$w$}] at (W) {};
\end{tikzpicture}
\caption{The watermelon graph $D_{m}$, with two vertices $u$ and $w$ joined by $m$ edges of lengths $\ell_{1},\dots,\ell_{m}$.
}
\label{fig:watermelon-intro}
\end{figure}

As an example, consider the watermelon, or dipole graphs, $D_m$, see Fig. \ref{fig:watermelon-intro}.
For the circle of length $L$, or $D_2$ with total length $L=\ell_1+ \ell_2$, the eigenvalues are $(2\pi j/L)^{2}$, $j\geq 0$. For $L$ algebraic, Lindemann's theorem shows that these numbers are transcendental, and the field they generate has transcendence degree one over $\mathbb{Q}$, since each of them is a rational multiple of $\pi^{2}$. Moreover, the special values $\zeta_{D_2}(m)$ in Theorem \ref{mthm:special} are up to a normalization Euler's famous values for the Riemann zeta function $\zeta_{\mathbb{Q}}(s)$ at $s=2m>0$.
For $m>2$, with algebraic and rationally independent lengths, the eigenvalues should have infinite
transcendence degree.  Yet $\zeta_{D_{m}}(1) = \mu_{1}^{-1}+\mu_{2}^{-1}+\cdots$ is a rational expression in
the lengths, hence an algebraic number; the formula is in Section~\ref{subsec:dipole}.

From theoretical chemistry, the value $\zeta_{\Gamma}(1)$ is the Kirchhoff index of the metric graph, a much used molecular
descriptor \cite{KR93, GM96}.  As an example, we use the formula for $\zeta_{\Gamma}(1)$ to give an evaluation of the Kirchhoff index for the molecular
graph of naphthalene; see Section~\ref{subsec:naphthalene}.

\subsection{Further main results}

The study of a continuous object by finite approximations of increasing size is a major topic in
mathematical physics, of importance for statistical physics and quantum field theory, as well as a
basic approach in numerical analysis.  For the determinant of the Laplacian $\text{\rm det}^{*}\Delta$,
with the zero eigenvalue removed, this began with Kasteleyn's work \cite{Ka61} in the 1960s.  In our own work we
have carried it out for tori in all dimensions, and for the full spectral zeta function rather than
only the value $-\zeta'(0)$ which corresponds to $\log\text{\rm det}^{*}\Delta$, see
\cite{CJK10, CJK12, FK17}.

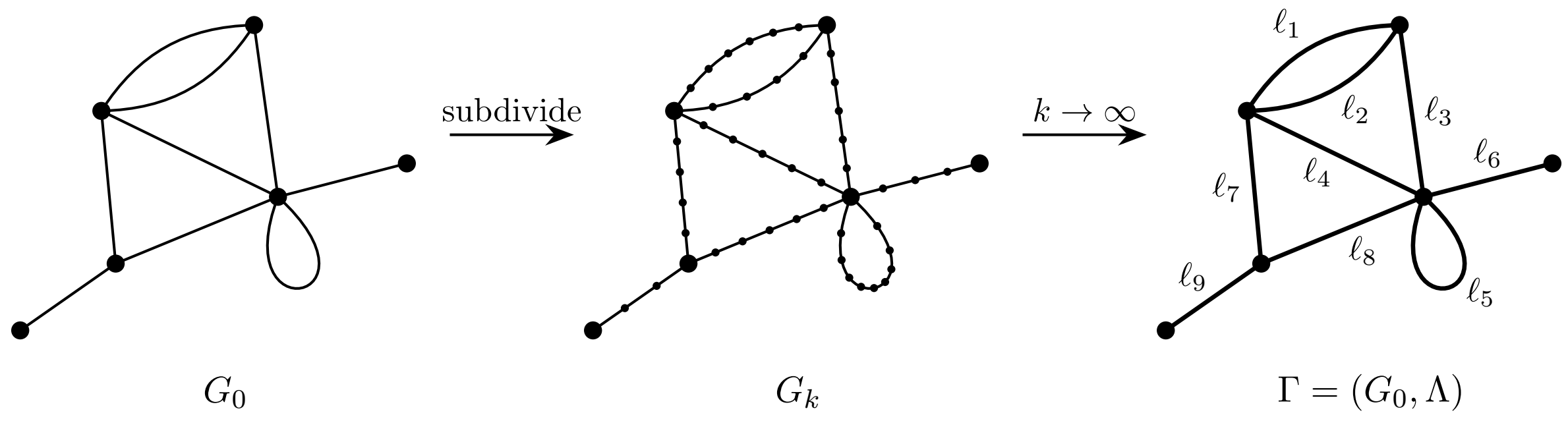
\begin{figure}[!ht]
\centering
\resizebox{\textwidth}{!}{%
\begin{tikzpicture}[v/.style={circle,fill,inner sep=1.9pt}, sv/.style={circle,fill,inner sep=0.85pt},
  e/.style={thick}, lab/.style={font=\small}]
\newcommand{\jksCoords}{%
  \coordinate (a) at (-1.5,0.35); \coordinate (b) at (0.1,1.25);
  \coordinate (c) at (0.35,-0.55); \coordinate (d) at (1.7,-0.2);
  \coordinate (e) at (-1.35,-1.25); \coordinate (f) at (-2.35,-1.95);}
\newcommand{\jksEdges}[1][]{%
  \draw[e,#1] (a) to[bend left=28] (b);
  \draw[e,#1] (a) to[bend right=28] (b);
  \draw[e,#1] (b) -- (c);
  \draw[e,#1] (a) -- (c);
  \draw[e,#1] (c) .. controls +(-0.55,-1.45) and +(1.2,-1.1) .. (c);
  \draw[e,#1] (c) -- (d);
  \draw[e,#1] (a) -- (e);
  \draw[e,#1] (e) -- (c);
  \draw[e,#1] (e) -- (f);}
\begin{scope}[shift={(0,0)}]
  \jksCoords \jksEdges
  \foreach \p in {a,b,c,d,e,f} \node[v] at (\p) {};
  \node at (-0.2,-2.6) {$G_0$};
\end{scope}
\draw[-{Stealth[length=3mm]},thick] (2.15,0.1) -- node[above,lab] {subdivide} (3.45,0.1);
\begin{scope}[shift={(6.0,0)}]
  \jksCoords \jksEdges
  \foreach \i in {1,...,6} \path (a) to[bend left=28] node[sv,pos=\i/7] {} (b);
  \foreach \i in {1,...,4} \path (a) to[bend right=28] node[sv,pos=\i/5] {} (b);
  \foreach \i in {1,...,5} \path (b) -- node[sv,pos=\i/6] {} (c);
  \foreach \i in {1,...,5} \path (a) -- node[sv,pos=\i/6] {} (c);
  \foreach \i in {1,...,9} \path (c) .. controls +(-0.55,-1.45) and +(1.2,-1.1) .. node[sv,pos=\i/10] {} (c);
  \foreach \i in {1,...,3} \path (c) -- node[sv,pos=\i/4] {} (d);
  \foreach \i in {1,...,4} \path (a) -- node[sv,pos=\i/5] {} (e);
  \foreach \i in {1,...,5} \path (e) -- node[sv,pos=\i/6] {} (c);
  \foreach \i in {1,...,2} \path (e) -- node[sv,pos=\i/3] {} (f);
  \foreach \p in {a,b,c,d,e,f} \node[v] at (\p) {};
  \node at (-0.2,-2.6) {$G_{k}$};
\end{scope}
\draw[-{Stealth[length=3mm]},thick] (8.15,0.1) -- node[above,lab] {$k\rightarrow\infty$} (9.45,0.1);
\begin{scope}[shift={(12.0,0)}]
  \jksCoords \jksEdges[line width=1.3pt]
  \path (a) to[bend left=28] node[lab,above left=-1pt] {$\ell_{1}$} (b);
  \path (a) to[bend right=28] node[lab,below right=-2pt] {$\ell_{2}$} (b);
  \path (b) -- node[lab,right] {$\ell_{3}$} (c);
  \path (a) -- node[lab,below,pos=0.4] {$\ell_{4}$} (c);
  \node[lab] at (0.95,-1.55) {$\ell_{5}$};
  \path (c) -- node[lab,above] {$\ell_{6}$} (d);
  \path (a) -- node[lab,left] {$\ell_{7}$} (e);
  \path (e) -- node[lab,below right=-2pt] {$\ell_{8}$} (c);
  \path (e) -- node[lab,above left=-2pt] {$\ell_{9}$} (f);
  \foreach \p in {a,b,c,d,e,f} \node[v] at (\p) {};
  \node at (-0.2,-2.6) {$\Gamma = (G_0,\Lambda)$};
\end{scope}
\end{tikzpicture}}
\caption{A combinatorial graph $G_0$ with a loop and a double edge; a discrete approximation $G_{k}$, in which each edge of length
$\ell_{j}$ is subdivided into about $k\ell_{j}$ segments; and the metric graph $\Gamma$ obtained as $k\rightarrow\infty$.}
\label{fig:three-step}
\end{figure}

In the present work, we approach the metric graph through finite graph approximations in a natural manner, see Figure~\ref{fig:three-step}. 
Let $G_0$ be a finite connected (combinatorial) graph, possibly with loops and multiple
edges, and let $\Gamma$ be the metric graph obtained by giving the edges arbitrary positive lengths
$\ell_{1},\dots,\ell_{E_{0}}$.

For $k\geq 1$, let $G_{k}$ be the discrete graph obtained by subdividing each
edge into approximately $k\ell_{j}$ segments, so that $G_{k}$ is a discretization of $\Gamma$ with mesh
$1/k$.  The Laplacian $\Delta_{\Gamma}$ has standard (Kirchhoff) vertex conditions, which are exactly the
conditions induced in the limit by the combinatorial Laplacians of $G_{k}$.  The Laplacian $\Delta_{\Gamma}$ is self-adjoint
and non-negative, and its spectrum is discrete, consisting of infinitely many eigenvalues of finite multiplicity
accumulating only at infinity \cite{BK13, Ku24}.

Let $0 = \lambda_{0,k}<\lambda_{1,k}\leq\cdots$ and $\zeta_{G_{k}}(s) = \sum_{j\geq 1}\lambda_{j,k}^{-s}$
be the sequence of eigenvalues and spectral zeta function
associated to the combinatorial Laplacian of $G_{k}$.  Let $d_{v}$ be the degree of the vertex $v$ of $G_{0}$, and $L$ the total length of $\Gamma$,
where $G_{0}$ has $E_{0}$ edges and $V_{0}$ vertices.


\begin{mainthm} \label{mthm:convergence}
Let $\Gamma = (G_{0},\Lambda)$ be a compact metric graph as above, with standard vertex
conditions, and let $G_{k}$ be its discrete approximations as described above.  Then the following hold.
\begin{enumerate}
\item For every $x\in\mathbb{C}$,
$$
 \lim_{k\rightarrow\infty}k^{V_{0}-E_{0}}\det\big(\Delta_{G_{k}}+x/k^{2}\big)=\Phi(x),
$$
the metric forest function, uniformly for $x$ in compact subsets of $\mathbb{C}$.  In particular $\Phi$ is entire.
\item $\Phi$ has order at most $1/2$, and its zeros are exactly the numbers $-\mu_{j}$, $j\geq 0$, counted with
multiplicity.
\item We have
$$
\Phi(x) = P_{1}\,x\prod_{j\geq 1}\Big(1+\frac{x}{\mu_{j}}\Big),\qquad P_{1} = \Phi'(0) = L\sum_{T}\prod_{\mathcal{E}_{j}\notin T}\ell_{j},
$$
where $T$ runs over the spanning trees of $G_{0}$.
\item For every $j$ we have $k^{2}\lambda_{j,k}\rightarrow\mu_{j}$ as $k\rightarrow\infty$.
\item For every $s\in\mathbb{C}$ with $\text{\rm Re}(s)>1/2$ we have
$k^{-2s}\zeta_{G_{k}}(s)\rightarrow\zeta_{\Gamma}(s)$ as $k\rightarrow\infty$.
\end{enumerate}
\end{mainthm}


Note that the formula for $\Phi$ as given in \eqref{eq:mforestfunction} is a finite sum, with one term for each spanning forest of $G_0$. No roots are marked, and the number of components is unrestricted. It ranges from one component for a spanning tree up to $V_0$ components for the forest with no edges.  The number of terms can be very large.

From its definition, the function $\Phi$ also has an Euler product.  Namely, for $\text{\rm Re}\sqrt{x}$ large,
$$
\Phi(x) = \frac{\prod_{v}d_{v}}{2^{E_{0}}}\,x^{\chi/2}\,e^{L\sqrt{x}}\prod_{\gamma}
\Big(1-\alpha(\gamma)\,e^{-\ell(\gamma)\sqrt{x}}\Big),
$$
the product running over the primitive periodic orbits of $\Gamma$, whose amplitudes $\alpha(\gamma)$ are built
from the scattering coefficients $2/d_{v}$ and $2/d_{v}-1$ at the vertices; see
Remark~\ref{rem:euler-product}.  For the circle of length $L$ this reads
$\Phi(x) = e^{L\sqrt{x}}+e^{-L\sqrt{x}}-2$.

It is important to note that we prove Theorem~\ref{mthm:convergence} without appealing to a trace formula for $\Gamma$.  Instead, we study the
asymptotic behavior of the sequence of finite graphs $\{G_{k}\}$ and use results from linear algebra, complex analysis
and elementary aspects of the spectral theory of $\Delta_{\Gamma}$. In fact, Roth's trace formula can be deduced from our approach; see \cite{JKS27}.


The restriction to $\text{\rm Re}\,s>1/2$ in
Theorem~\ref{mthm:convergence}(5) is sharp.  Suppose the subdivision is exact, meaning $n_{j,k}+1 = k\ell_{j}$.  Then $G_{k}$ has
$V_{k} = kL+V_{0}-E_{0}$ vertices, so $\zeta_{G_{k}}(0) = V_{k}-1 = kL+\chi-1$, which tends to infinity, whereas $\zeta_{\Gamma}(0) = \chi/2-1$.
The discrete zeta function counts all of its roughly $kL$ eigenvalues, and to the left of $\text{\rm Re}\,s = 1/2$ their
contribution, which is proportional to the total length, dominates.

By analogy with the results of
\cite{CJK10, CJK12, FK17} for discrete tori degenerating to real tori, and the expansion put forward in \cite{RV15}, one expects an asymptotic expansion of the form
$$
k^{-2s}\zeta_{G_{k}}(s) = \zeta_{\Gamma}(s)+ L\,\zeta_{\mathbb{Z}}(s)\,k^{1-2s}+k^{-2s}\sum_{v}c_{d_{v}}(s)+\cdots,\qquad
\zeta_{\mathbb{Z}}(s) = \frac{4^{-s}\,\Gamma(1/2-s)}{\sqrt{\pi}\,\Gamma(1-s)},
$$
with a length term $L\,\zeta_{\mathbb{Z}}(s)\,k^{1-2s}$ coming from the Weyl law, a constant term $\zeta_{\Gamma}(s)$, and local
vertex terms $c_{d}(s)$.  Here $\zeta_{\mathbb{Z}}(s) = \frac{1}{2\pi}\int_{0}^{2\pi}(4\sin^{2}(\theta/2))^{-s}\,d\theta$ is the spectral zeta
function of the infinite path $\mathbb{Z}$, contributing once per unit length.  The metric graph enters through the
constant term $\zeta_{\Gamma}(s)$, and $c_{d}(s)$ are local contributions of the vertices of degree $d$, with $c_{2} = 0$.  Two
exact checks support this form.  At $s = -1$, $\zeta_{G_{k}}(-1) = \operatorname{tr}\Delta_{G_{k}} = 2E_{k} = 2kL$.  This equals the
length term, since $\zeta_{\mathbb{Z}}(-1) = 2$, and $\zeta_{\Gamma}(-1) = 0$; hence $c_{d}(-1) = 0$.  At $s = 0$, the length term is
$kL$ since $\zeta_{\mathbb{Z}}(0) = 1$, and the constant term $\chi-1$ requires $\sum_{v}c_{d_{v}}(0) = \chi/2$, which holds with
the local values $c_{d}(0) = (2-d)/4$, since $\sum_{v}d_{v} = 2E_{0}$.  Establishing this expansion, and determining the
vertex terms $c_{d}(s)$, is a natural problem which we do not pursue here except for the determinant of the Laplacian in Theorem \ref{mthm:det} below. 

We now identify the function $\Phi(x)$ 
as an expression involving the zeta-regularized determinant of
$\Delta_{\Gamma}+x$. For this purpose, define $M(x)$ to be the $V_{0}\times V_{0}$
\emph{Dirichlet-to-Neumann matrix}. It has entries
\begin{align*}
M(x)_{vw} &= -\sum_{\mathcal{E}_{j}=\{v,w\}}\frac{\sqrt{x}}{\sinh(\sqrt{x}\,\ell_{j})}\qquad(v\neq w),\\
M(x)_{vv} &= \sum_{h\text{ \rm at }v}\frac{\sqrt{x}}{\tanh(\sqrt{x}\,\ell_{e(h)})}-\sum_{\text{\rm loops at }v}\frac{2\sqrt{x}}{\sinh(\sqrt{x}\,\ell)}.
\end{align*}
The off-diagonal sum runs over the edges joining $v$ and $w$, and on the diagonal the first sum runs over the
half-edges at $v$ and the second over the loops at $v$.

\begin{mainthm}\label{mthm:det}
Let $\Delta_{D}$ be the Laplacian with Dirichlet conditions at all vertices, and let $Q_{\Gamma}$ be the period
matrix of $\Gamma$, defined below.  Then the following hold.
\begin{enumerate}
\item For every $x>0$,
$\displaystyle\text{\rm det}(\Delta_{\Gamma}+x) = \frac{2^{E_{0}}}{\prod_{v}d_{v}}\,\Phi(x)$.
\item For every $x>0$,
$\displaystyle\text{\rm det}(\Delta_{\Gamma}+x) = \frac{1}{\prod_{v}d_{v}}\;\text{\rm det}(\Delta_{D}+x)\;\det M(x)$.
\item We have
$$
\text{\rm det}^{*}\Delta_{\Gamma} = \frac{2^{E_{0}}L}{\prod_{v}d_{v}}\sum_{T}\prod_{\mathcal{E}_{j}\notin T}\ell_{j},
$$
where $T$ runs over the spanning trees of $G_{0}$.
\item Equivalently,
$\displaystyle\text{\rm det}^{*}\Delta_{\Gamma} = \prod_{v}\frac{2^{d_{v}/2}}{d_{v}}\cdot\text{\rm vol}(\Gamma)\cdot\det Q_{\Gamma}$.
\item For the graphs $G_{k}$, as $k\rightarrow\infty$,
$$
\log\text{\rm det}^{*}\Delta_{G_{k}} = (E_{0}-V_{0}+2)\log k+\log\text{\rm det}^{*}\Delta_{\Gamma}-E_{0}\log 2+\sum_{v}\log d_{v}+o(1).
$$
\end{enumerate}
\end{mainthm}

Parts (2) and (3) are known.  Part (2) is \cite[Theorem~1]{Fr06}, for the Laplacian with standard vertex
conditions, and part (3) is \cite[Theorem~2]{Fr06}.  We reprove both results from \cite{Fr06}  by studying limits of the determinants of the
discrete Laplacians of the graphs $G_{k}$, so Theorem~\ref{mthm:det} also gives an independent proof.
Parts (1), (4) and (5) are new.  Part (1) identifies $\text{\rm det}(\Delta_{\Gamma}+x)$ with $\Phi$, up to the
constant $2^{E_{0}}/\prod_{v}d_{v}$, and so transfers the forest expansion of Theorem~\ref{mthm:convergence} to
the determinant.  Part (4) turns the spanning-tree sum into a volume.  Part (5) realizes
$\text{\rm det}^{*}\Delta_{\Gamma}$ as the renormalized limit of the determinants
$\text{\rm det}^{*}\Delta_{G_{k}}$, with every correction term explicit.  A second derivation of part (3), from
the convergence of the discrete determinants, is given in Remark~\ref{rem:friedlander}.


Theorem~\ref{mthm:special} expresses the special values of $\zeta_{\Gamma}$ through the spanning forests of $G_{0}$.  The relation
can also be read in the opposite direction. The numbers of spanning forests of the discrete graphs $G_{k}$ are
governed, asymptotically as $k\rightarrow\infty$, by the regularized determinant and the special values of $\Gamma$.

\begin{mainthm}\label{mthm:forests}
Let $c_{m}(G_{k})$ be the number of rooted spanning forests of $G_{k}$ with $m$ components. 
Then, as $k\rightarrow\infty$,
$$
c_{m}(G_{k}) = k^{E_{0}-V_{0}+2m}\Big(\frac{\prod_{v}d_{v}}{2^{E_{0}}}\,\text{\rm det}^{*}\Delta_{\Gamma}\cdot
h_{m-1}\big(\zeta_{\Gamma}(1),\dots,\zeta_{\Gamma}(m-1)\big)+o(1)\Big),
$$
where $h_{n}$ is the universal polynomial expressing the elementary symmetric function of degree $n$ through the
power sums.
\end{mainthm}

The polynomials $h_{n}$ are classical.  For instance, $h_{1}(p_{1}) = p_{1}$ and $h_{2}(p_{1},p_{2}) = \frac12(p_{1}^{2}-p_{2})$,
and in general they are determined by the identity of formal power series
$$
\sum_{n\geq 0}h_{n}(p_{1},\dots,p_{n})\,x^{n} = \exp\big(\sum_{r\geq 1}(-1)^{r-1}p_{r}x^{r}/r\big),
$$
see \cite[Chapter~I, \S2]{Mac95} and Section~\ref{subsec:forest-counts}.  (Note that one should not
mistake these polynomials $h_{n}$ with the complete homogeneous symmetric functions, which \cite{Mac95} also denotes by $h_{n}$.)

When the subdivision is exact, $n_{j,k}+1 = k\ell_{j}$, the count $c_{m}(G_{k})$ is a polynomial in $k$ of degree
$E_{0}-V_{0}+2m$ with leading coefficient $P_{m}$, so the error term in Theorem~\ref{mthm:forests} can be expanded to
every order in $1/k$.  The normalized characteristic polynomials then depend analytically on $h = 1/k$, and each
simple eigenvalue satisfies $k^{2}\lambda_{j,k} = \mu_{j}(1/k)$ for a function $\mu_{j}(h)$ analytic near $h = 0$; see
Remark~\ref{rem:analytic-in-h}.


\subsection{Our approach}\label{subsec:approach}

Our method of proof involves a careful study of the metric forest function $\Phi(x)$ of
\eqref{eq:mforestfunction} from various points of view.  Let us now describe this in detail.

We begin with the definition \eqref{eq:mforestfunction} of $\Phi(x)$.  Let us use the notation
$$
S(x;s) = \frac{\sinh(\sqrt{x}\,s)}{\sqrt{x}},\qquad S_{j}(x) = S(x;\ell_{j}),\qquad
\Theta_{C}(x) = \sqrt{x}\sum_{h\text{ \rm at }C}\tanh\frac{\sqrt{x}\,\ell_{e(h)}}{2},
$$
where the last sum is over the half-edges at $C$, so that an edge with both ends in $C$ contributes twice.
With this notation, we have that
$$
\Phi(x) = \sum_{F}\ \prod_{\mathcal{E}_{j}\notin F}S_{j}(x)\ \prod_{C\in F}\Theta_{C}(x)
\qquad(x\in\mathbb{C}),
$$
This is a finite sum, with one term for each spanning forest $F$ of the combinatorial graph $G_{0}$.  The
two products run over the edges of $G_{0}$ outside $F$ and over the connected components $C$ of $F$.
In general terms, the results in this paper follow by establishing three further
descriptions of $\Phi(x)$, a determinant of fixed size, a limit of determinants of discrete Laplacians, and a
product over the spectrum of $\Delta_{\Gamma}$.  Let us discuss each of these alternate realizations.

The first description is that $\Phi(x)$ is the determinant of a weighted combinatorial Laplacian of $G_{0}$,
with vertex weights added along the diagonal, times a multiplicative factor.  Specifically, we will define
weights $w(x)$ on the edges of $G_{0}$ and weights $p(x)$ on its vertices, both depending on the lengths and
on $x$, so that
$$
\Phi(x) = \Big(\prod_{j}S_{j}(x)\Big)\det\big(L_{w(x)}+\text{\rm diag}(p(x))\big)
$$
for every $x\in\mathbb{C}$ outside the discrete set $Z$ where some $S_{j}$ vanishes, where $L_{w}$ is the
weighted Laplacian of $G_{0}$.  The construction is as follows.  The edge $\mathcal{E}_{j}$ is given the
weight $w_{j}(x) = 1/S_{j}(x)$ and the vertex $v$ the weight
$$
p_{v}(x) = \sqrt{x}\sum_{h\text{ \rm at }v}\tanh\frac{\sqrt{x}\,\ell_{e(h)}}{2},
$$
which is made so that $\sum_{v\in V(C)}p_{v}(x) = \Theta_{C}(x)$ for every component $C$.  The generalized
matrix-forest theorem of Lemma~\ref{lem:matrix-forestnew} expands the determinant of any weighted Laplacian
with an added diagonal over the spanning forests of the graph, which here gives
$$
\det\big(L_{w(x)}+\text{\rm diag}(p(x))\big)
= \sum_{F}\ \prod_{\mathcal{E}_{j}\in F}\frac{1}{S_{j}(x)}\ \prod_{C\in F}\Theta_{C}(x).
$$
Multiplying by $\prod_{j}S_{j}(x)$ turns the product over the edges of $F$ into the product over the edges
outside $F$, which gives
\begin{equation}\label{eq:Phi-detM}
\Phi(x) = \Big(\prod_{j}S_{j}(x)\Big)\det M(x),\qquad\text{with}\qquad
M(x) = L_{w(x)}+\text{\rm diag}(p(x))
\end{equation}
where $j$ ranges over the set of edges of $G_{0}$; see Propositions~\ref{prop:Mlaplace} and
\ref{prop:Phi-DtN}.

The defining sum \eqref{eq:mforestfunction} has one term for every spanning forest of $G_{0}$.  In
\eqref{eq:Phi-detM} it has become the determinant of the single $V_{0}\times V_{0}$ matrix $M(x)$,
multiplied by the single product $\prod_{j}S_{j}(x)$ of $E_{0}$ scalars.

As it turns out, the matrix $M(x)$ is the Dirichlet-to-Neumann matrix of $\Gamma$, and this is what ties
$\Phi$ to the spectrum.
Let $\varphi$ be a complex-valued function on the vertices of $G_{0}$, that is, a vector
$\varphi\in\mathbb{C}^{V_{0}}$ with one entry $\varphi_{v}$ for each vertex $v$.  For $x\notin Z$ there is a
unique extension $f_{\varphi}$ of $\varphi$ to all of $\Gamma$ which solves $f''=xf$ on every edge; see
Lemma~\ref{lem:edge-solutions}.  Explicitly, on the edge $\mathcal{E}_{j}$, parametrized by arclength
$s\in[0,\ell_{j}]$,
$$
f_{\varphi}(s) = \frac{a\,S(x;\ell_{j}-s)+b\,S(x;s)}{S(x;\ell_{j})},
$$
where $a$ and $b$ are the values of $\varphi$ at the endpoints of $\mathcal{E}_{j}$ lying at $s = 0$ and at
$s = \ell_{j}$.  Indeed
$$
f_{\varphi}(0) = a,\qquad f_{\varphi}(\ell_{j}) = b,\qquad
f_{\varphi}''(s) = x\,f_{\varphi}(s)\quad\text{for }0<s<\ell_{j},
$$
so $f_{\varphi}$ agrees with $\varphi$ at both ends of every edge and is therefore continuous on $\Gamma$.  At each vertex the
sum of its outgoing derivatives is the corresponding entry of $-M(x)\varphi$.  So $f_{\varphi}$ satisfies the Kirchhoff conditions exactly
when $M(x)\varphi = 0$, and $\varphi\mapsto f_{\varphi}$ is an isomorphism from $\ker M(x)$ onto
$\ker(\Delta_{\Gamma}+x)$; see Lemma~\ref{lem:DtN-kernel}.  The same computation gives the energy identities
$$
\varphi^{T}M(x)\varphi = \sum_{j}\int_{0}^{\ell_{j}}\big(|f_{\varphi}'|^{2}+x|f_{\varphi}|^{2}\big)\,ds,
\qquad
\varphi^{T}M'(x)\varphi = \|f_{\varphi}(x)\|_{L^{2}(\Gamma)}^{2}
$$
of Lemma~\ref{lem:energynew}, in which $M'(x)$ is the entrywise derivative $\frac{d}{dx}M(x)$.  These identities imply that
$M(x)$ and $M'(x)$ are positive definite for $x>0$.  The isomorphism $\ker M(x)\cong\ker(\Delta_{\Gamma}+x)$ implies that the
zeros of $\Phi$ are the eigenvalues of $\Delta_{\Gamma}$.  The positivity of $M'$ implies that the order of each zero is exactly the
multiplicity of the eigenvalue, and not larger.

The entries of $M$ are built from $\sqrt{x}\coth(\sqrt{x}\,\ell_{j})$ on the diagonal and
$\sqrt{x}/\sinh(\sqrt{x}\,\ell_{j})$ off it.  In a sector $|\arg x|\leq\pi-\varepsilon$ the first tends to
$\sqrt{x}$ and the second to $0$, both exponentially fast in $\text{\rm Re}\sqrt{x}$, so that each half-edge
at a vertex contributes $\sqrt{x}$ to its diagonal entry while the off-diagonal entries die.  Hence
$$
M(x) = \sqrt{x}\,\text{\rm diag}(d_{v})+O\big(|x|^{1/2}e^{-c\sqrt{|x|}}\big)
$$
in sectors; see Proposition~\ref{prop:DtN-properties}.  Together with the elementary asymptotics of the
$S_{j}$, the factorization \eqref{eq:Phi-detM} turns this into the asymptotics of $\Phi$ itself, and with
them those of $\Phi'/\Phi$.  By the product formula of the third description,
$\Phi'(x)/\Phi(x) = \sum_{j\geq 0}(x+\mu_{j})^{-1}$, which is the Laplace transform of the theta function
$Q(t) = \sum_{j\geq 0}e^{-\mu_{j}t}$.  Inverting the transform gives the small time expansion of $Q$, and
with it the constant $2^{E_{0}}/\prod_{v}d_{v}$ in Theorem~\ref{mthm:det}.

The second description is the analytic one.  It uses the same weighting of $G_{0}$ by edges and by
vertices as the first, with discrete weights in place of $w(x)$ and $p(x)$, which converge to them after
rescaling.  Let $G_{k}$ be
the graph obtained from $G_{0}$ by subdividing each edge $\mathcal{E}_{j}$ into $n_{j,k}+1$ segments, where
$n_{j,k}/k\rightarrow\ell_{j}$.  The interior vertices of a subdivided edge form a path, and the corresponding
block of $\Delta_{G_{k}}+u$ is tridiagonal with $u+2$ on the diagonal and $-1$ off it.  Its determinant is
$D_{n}(u)$, and the interior vertices act on the two ends of the path through
$\eta_{n}(u) = 1-(D_{n-1}(u)+1)/D_{n}(u)$, where, writing $u = 2\cosh\beta-2$,
$$
D_{n}(u) = \frac{\sinh((n+1)\beta)}{\sinh\beta},\qquad
\eta_{n}(u) = \frac{2\sinh(\beta/2)\sinh(n\beta/2)}{\cosh((n+1)\beta/2)}.
$$
Using Schur's determinant formula to eliminate the interior vertices reduces the roughly $kL$ vertices of
$G_{k}$ to a determinant of the fixed size $V_{0}$, and expanding that determinant over forests gives the
exact identity
\begin{equation}\label{eq:exact-Gk}
\det\big(\Delta_{G_{k}}+u\big) = \sum_{F}\ \prod_{\mathcal{E}_{j}\notin F}D_{n_{j,k}}(u)\
\prod_{C\in F}\Big(u\,|V(C)|+\sum_{h\text{ \rm at }C}\eta_{n_{e(h),k}}(u)\Big)
\end{equation}
of Proposition~\ref{prop:exact-forest}.  For every $k$, and before any limit is taken, the determinant of the
discrete Laplacian is already a sum over the spanning forests of $G_{0}$, of the same shape as
\eqref{eq:mforestfunction}.

Now set $u = x/k^{2}$.  Then $\beta = \sqrt{x}/k+O(k^{-3})$, so $k\beta\rightarrow\sqrt{x}$ and
$(n_{j,k}+1)\beta\rightarrow\sqrt{x}\,\ell_{j}$, whence
$$
\frac{D_{n_{j,k}}(x/k^{2})}{k}\longrightarrow S_{j}(x),\qquad
k\,\eta_{n_{j,k}}(x/k^{2})\longrightarrow\sqrt{x}\tanh\frac{\sqrt{x}\,\ell_{j}}{2},
$$
uniformly on compact sets; see Lemma~\ref{lem:continuum-limit}.  These two limits are the building blocks
of \eqref{eq:mforestfunction}.  The first gives the edge factors and the second the half-edge terms of
$\Theta_{C}$, while the term $u|V(C)|$ contributes nothing after rescaling.

This also explains the normalizing factor.  A spanning forest with $r$ components has $V_{0}-r$ edges, so its
term in \eqref{eq:exact-Gk} carries $E_{0}-V_{0}+r$ edge factors and $r$ component factors.  Each edge factor is of size
$k$ and each component factor of size $1/k$, so the term is of size
$k^{E_{0}-V_{0}+r}k^{-r} = k^{E_{0}-V_{0}}$.  The number of components $r$ cancels, so each term is scaled
by $k^{V_{0}-E_{0}}$.  The same factor therefore normalizes the whole sum at once.
Each forest term converges to the corresponding term of $\Phi$, and the forests are finite in number.  As a
result, we have that
$$
\lim_{k\rightarrow\infty}k^{V_{0}-E_{0}}\det\big(\Delta_{G_{k}}+x/k^{2}\big) = \Phi(x)
$$
uniformly on compact subsets of $\mathbb{C}$; see Theorem~\ref{thm:Phi-limit}, which is
Theorem~\ref{mthm:convergence}(1).  In particular $\Phi$ is entire, with no exceptional set, and
the argument uses no more than $n_{j,k}/k\rightarrow\ell_{j}$.  The same elimination shows that the discrete
Dirichlet-to-Neumann matrices converge as well, meaning that $k\,S_{k}(x/k^{2})\rightarrow M(x)$, where $S_{k}(u)$ is the
Schur complement onto the vertices of $G_{0}$; see Remark~\ref{rem:discrete-DtN}.

The third description is that $\Phi$ is a Hadamard product over the spectrum of $\Delta_{\Gamma}$.
Specifically, we will prove that
$$
\Phi(x) = P_{1}\,x\prod_{j\geq 1}\Big(1+\frac{x}{\mu_{j}}\Big),\qquad P_{1} = \Phi'(0),
$$
where $0 = \mu_{0}<\mu_{1}\leq\mu_{2}\leq\cdots$ are the eigenvalues of $\Delta_{\Gamma}$.  To prove this we need to know where $\Phi$
vanishes and a bound on the growth of $\Phi(x)$ in $x$.

The zeros come from the two previous descriptions.  For
$x_{0}\notin Z$ the product $\prod_{j}S_{j}$ is nonzero at $x_{0}$, so $\Phi$ and $\det M$ vanish there to the
same order, and Lemma~\ref{lem:DtN-kernel} ensures that this order is the multiplicity of $-x_{0}$ as an
eigenvalue of $\Delta_{\Gamma}$.  The set
$Z$ has to be treated separately, since there $M$ has poles and the extension $f_{\varphi}$ is not defined.
The points of $Z$ are the numbers $-(\pi n/\ell_{j})^{2}$, the Dirichlet eigenvalues of the individual edges, which
are exactly the $x_{0}$ at which an eigenfunction may vanish at both ends of an edge.  They are handled by
subdividing.  Inserting a vertex of degree two
in the interior of an edge leaves both the metric graph $\Gamma$ and the function $\Phi$ unchanged, by
Lemma~\ref{lem:subdivision}, while the two new lengths replace the old one in $Z$, so for any prescribed
$x_{0}$ the subdivision can be chosen to put $x_{0}$ outside the set $Z$ of the new model, where the previous
case applies.  Hence $\Phi(x_{0}) = 0$ if and only if $-x_{0}$ is an eigenvalue of $\Delta_{\Gamma}$, and the
order of the zero is its multiplicity; see Theorem~\ref{thm:zeros}, which is
Theorem~\ref{mthm:convergence}(2).

This is what the degeneration from the graphs $G_{k}$ was for.  Each
$k^{V_{0}-E_{0}}\det(\Delta_{G_{k}}+x/k^{2})$ is a polynomial, so $\Phi$ is entire, whereas individual
terms of \eqref{eq:mforestfunction} have poles on $Z$, and the divisor of $\Phi$ is pinned down at every
point, the points of $Z$ included.  Being entire with a known divisor is exactly what Hadamard's theorem
requires.

The growth is estimated from the definition \eqref{eq:mforestfunction} itself, term by term over the forests.  On an edge $|S(x;s)|\leq s\,e^{s\sqrt{|x|}}$, and
each $\Theta_{C}$ is bounded by a constant multiple of $\sqrt{|x|}$ away from the zeros of the hyperbolic
cosines, so a forest with $m$ components contributes a factor $|x|^{m/2}$ and an exponential
$e^{\sqrt{|x|}\sum_{\mathcal{E}_{j}\notin F}\ell_{j}}$.  A forest has at most $V_{0}$ components, and the
edges outside it have total length at most $L$.  This gives the bound on circles $|x| = \rho^{2}$ with
$\rho$ at a fixed distance from the bad radii, and the maximum principle carries it to every $x$, so that
$$
|\Phi(x)|\leq A\,(1+|x|)^{V_{0}/2}e^{L\sqrt{|x|}}\qquad(x\in\mathbb{C}),
$$
so $\Phi$ has order at most $1/2$; see Lemma~\ref{lem:growth}.  The total length of $\Gamma$ thus appears as the
exponential type of $\Phi$.

Since the order is less than one, Hadamard's theorem writes $\Phi$ as a monomial times its Weierstrass
product over its zeros, with no exponential factor and with the product converging absolutely.  The zeros
being known, this reads
$$
\Phi(x) = c\,x^{a}\prod_{j\geq 1}\Big(1+\frac{x}{\mu_{j}}\Big)
$$
for some integer $a\geq 0$ and some constant $c$, which are determined as follows.

The exponent is $a = 1$.  Indeed $\Gamma$ is connected, so $\mu_{0} = 0$ is a simple eigenvalue and the zero
of $\Phi$ at the origin is simple, the remaining zeros being the $-\mu_{j}$ with $j\geq 1$.  The constant is
then $c = \Phi'(0)$, and it is evaluated from the definition.  Each component factor vanishes to first order,
$\Theta_{C}(x) = x\,\ell^{*}(C)+O(x^{2})$ with $\ell^{*}(C)$ half the sum of $\ell_{e(h)}$ over the half-edges
$h$ at $C$, so a forest with $r$ components contributes at order $x^{r}$ and only the spanning trees reach
the coefficient of $x$.  A spanning tree contains every vertex, so $\ell^{*}(T) = L$, and
$$
c = P_{1} = \Phi'(0) = L\sum_{T}\prod_{\mathcal{E}_{j}\notin T}\ell_{j},
$$
the sum being over the spanning trees $T$ of $G_{0}$; see Theorem~\ref{thm:hadamard}, which is
Theorem~\ref{mthm:convergence}(3).  One entire function now has a forest expansion on one
side and the spectrum of $\Delta_{\Gamma}$ on the other.

From here the argument goes in two directions.  Taking logarithms in the Hadamard product gives
$$
\log\Big(\frac{\Phi(x)}{P_{1}x}\Big) = \sum_{r\geq 1}\frac{(-1)^{r-1}\zeta_{\Gamma}(r)}{r}\,x^{r}
\qquad(|x|<\mu_{1}),
$$
so the Taylor coefficients $P_{m}$ of $\Phi$, which are forest sums, determine the special values.  The ratios
$P_{m+1}/P_{1}$ are the elementary symmetric functions of the numbers $1/\mu_{j}$, and Newton's identities turn
them into $\zeta_{\Gamma}(n)$ for every integer $n\geq 1$; see Theorems~\ref{mthm:special} and
\ref{mthm:forests}.
In the other direction, the asymptotics of $M(x)$ in sectors give those of $\Phi$, hence the theta function, the
meromorphic continuation of $\zeta_{\Gamma}$, its values at the non-positive integers and the identity
$$
\text{\rm det}(\Delta_{\Gamma}+x) = \frac{2^{E_{0}}}{\prod_{v}d_{v}}\,\Phi(x)\qquad(x>0)
$$
of Theorem~\ref{mthm:det}(1), which transfers the forest expansion to the determinant.

\subsection{Comparison with known results}\label{subsec:comparison}
For background on quantum graphs we refer to \cite{Kuc08, BK13, Ku24}.

The spectral zeta functions of metric graphs are studied through contour integrals in \cite{HK11, HKT12}, with
applications to the determinant, the vacuum energy $\zeta_{\Gamma}(-1/2)$ and the heat coefficients, but those
representations are confined to strips such as $-\frac12<\text{\rm Re}\,s<1$ and do not reach the positive
integers.  In the equilateral case, \cite{HW18} relates the zeta functions of discrete and metric graphs through
von Below's correspondence \cite{vB85}.  To our knowledge the evaluation of $\zeta_{\Gamma}(n)$ at every positive
integer, for arbitrary edge lengths, is new.  It does not follow from expanding the formula of \cite{Fr06} at
$x = 0$.  That would require two further ingredients.  The first is the exponential factor by which the
regularized determinant and the Hadamard product differ \cite{Vo87, JoLa93, FJS19}, and the second is the
Taylor expansion of \eqref{eq:Phi-detM}, whose entries have poles that cancel only in the determinant.
Theorem~\ref{mthm:convergence}(3) supplies the first, with no exponential factor, and the spanning-forest
expansion, based on the matrix-forest theorem \cite{KC74, CS97, Ch82}, supplies the second.  For the circle,
$\zeta_{\Gamma}(s) = 2(L/2\pi)^{2s}\zeta_{\mathbb Q}(2s)$, so our evaluation extends Euler's formulas for
$\zeta_{\mathbb Q}$ at the even integers to arbitrary compact metric graphs.  Related invariants of metrized
graphs, defined with respect to the canonical measure rather than arc length, are the tau constant
of \cite{BR07} and the invariants of \cite{Ci11} in arithmetic geometry.

The convergence of the rescaled spectra of discrete models to the spectrum of a metric graph was proved in
\cite{Fa06}, and Theorem~\ref{mthm:convergence} recovers it from a stronger statement, the convergence of the
normalized characteristic polynomials to an explicit entire function.  What is new for metric graphs are the vertices, which
contribute local terms such as the constant $\chi/2$ in the theta expansion and the factor
$\prod_{v}d_{v}/2^{E_{0}}$ in the determinant.

Theorem~\ref{mthm:det} recovers, for standard vertex conditions, the formula of \cite{Fr06}, which holds for
Schr\"odinger operators with vertex conditions of $\delta$-type and is an instance of the gluing formulas of
\cite{BFK92}.  A second derivation, from Roth's trace formula, is given in \cite{Te10}.  Here the formula appears
as the limit of Schur's determinant formula \cite{Sc17} for the discrete graphs, and we write it in a
Cheeger--M\"uller form with the structure of a Quillen metric \cite{Qu85}.  In \cite{HW23} it is shown that the
determinant of an equilateral or nearly equilateral metric graph determines its number of spanning trees, and
Theorem~\ref{mthm:det} gives the length-weighted count exactly, for arbitrary lengths.

\subsection{Outline of the paper}\label{subsec:outline}

The main theorems are proved one part at a time.  Theorem~\ref{mthm:convergence} has part (1) in
Theorem~\ref{thm:Phi-limit}, part (2) in Theorem~\ref{thm:zeros} with Lemma~\ref{lem:growth}, part (3) in
Theorem~\ref{thm:hadamard}, part (4) in Corollary~\ref{cor:spectral-convergence} and part (5) in
Corollary~\ref{cor:zeta-degeneration}.  Theorem~\ref{mthm:special} has parts (1) and (2) in
Corollary~\ref{cor:zeta-continuation}, part (3) in Corollary~\ref{cor:zeta-integers} and parts (4) and (5) in
Theorem~\ref{thm:newton}.  Theorem~\ref{mthm:det} has parts (1) and (3) in Theorem~\ref{thm:det}, part (2) in
Corollary~\ref{cor:friedlander}, part (4) in Corollary~\ref{cor:CM} and part (5) in
Corollary~\ref{cor:main}.  Theorem~\ref{mthm:forests} is part (2) of Theorem~\ref{thm:forest-counts}, and
Corollary~\ref{cor:arithmetic} is proved in Section~\ref{sec:arithmetic}.

Section~\ref{sec:examples} works out the circle, the flower and dipole graphs and the molecular graph of
naphthalene, reads $\zeta_{\Gamma}(1)$ as an average effective resistance, and closes with a numerical check
on the complete graph on four vertices.

In short, Sections~\ref{sec:combinatorialgraphs} to \ref{sec:graphapproximations} assemble the tools,
Section~\ref{sec:forest-expansion} produces $\Phi$, Section~\ref{sec:zeros} attaches it to the spectrum,
Section~\ref{sec:heat-degeneration} completes Theorem~\ref{mthm:convergence}, Section~\ref{sec:det} turns
$\Phi$ into the determinant, Section~\ref{sec:explicit} makes the invariants explicit, and
Section~\ref{sec:examples} gives examples.

\section{Combinatorial graphs}\label{sec:combinatorialgraphs}

\subsection{Graphs and spanning forests}\label{subsec:graphsandspanningforests}

Let $G$ denote a finite connected graph, possibly with loops and multiple edges, with $V$ vertices and
$E$ unoriented edges $\mathcal{E}_{1},\dots,\mathcal{E}_{E}$.  The degree $d_{v}$ of a vertex $v$ is the number of
edges meeting $v$, a loop counting twice.  We will sometimes deal with oriented edges.  Every edge gives
rise to two oriented edges.  An edge has two ends, and to orient it is to choose which of the two ends is
the initial one.  A loop at $v$, both of whose ends are attached to $v$, also gives two oriented edges,
each of them originating and terminating at $v$.  The degree $d_{v}$ is then the number of oriented edges
originating at $v$, a loop once more counting twice.
The Euler characteristic of $G$ is $\chi = V-E$.

A spanning tree, usually denoted by $T$, is a maximal subgraph of $G$ containing no cycles, meaning that no further edge of $G$ can be added without creating a cycle.  Since $G$ is
connected, a spanning tree contains every vertex of $G$ and is itself connected, and it has $V-1$ edges.
It is obtained by removing $E-V+1$ edges from $G$.  A loop is a cycle, so no spanning tree contains a
loop, and of several edges joining the same pair of vertices a spanning tree contains at most one.  The
number of spanning trees of $G$ is $\tau(G)$.

More generally, for $1\leq m\leq V$, an $m$-spanning forest $\mathcal{F}$ (or denoted $F$) of $G$ is a
spanning subgraph with $m$ connected components $\mathcal{F}_{1},\dots,\mathcal{F}_{m}$, each of
which is a tree.  An $m$-spanning forest has $V-m$ edges, so it is obtained by removing
$E-V+m$ edges.  Define the weight of $\mathcal{F}$ to be
$$
w(\mathcal{F}) = \prod\limits_{i=1}^{m}V(\mathcal{F}_{i}),
$$
where $V(\mathcal{F}_{i})$ denotes the number of vertices of the component $\mathcal{F}_{i}$.
In alternative notation, let $|C|$ denote the number of vertices of a component $C$. Finally, for $1\leq m\leq V$ let
\begin{equation}\label{eq:cm}
c_{m}(G) = \sum\limits_{\mathcal{F}}w(\mathcal{F}),
\end{equation}
where the sum is over all $m$-spanning forests of $G$.  A forest with components
$\mathcal{F}_{1},\dots,\mathcal{F}_{m}$ admits exactly $\prod_{i}V(\mathcal{F}_{i})$ choices of one root
vertex in each component, so $c_{m}(G)$ is also the number of $m$-spanning forests of $G$ with a root
chosen in each component.  The only forest with $V$ components is the one with no edges, so
$c_{V}(G) = 1$.  These weighted counts are the coefficients of the polynomial $\det(\Delta_{G}+x)$, as recalled in
\S\ref{subsec:combinatoriallaplacian} below.  An end of an edge is also called a \emph{half-edge}, so
every edge has two half-edges, both half-edges of a loop at $v$ are at $v$, and $d_{v}$ is the number of
half-edges at $v$.  We write $e(h)$ for the edge containing a half-edge $h$, $v(h)$ for its vertex, and
$v_{h}$ for the other endpoint of $e(h)$.


\subsection{The combinatorial Laplacian}\label{subsec:combinatoriallaplacian}

We write $\Delta_{G}$ for the combinatorial Laplacian
$$
\Delta_{G}f(v)=\sum_{h\,:\,v(h)=v}\big(f(v)-f(v_{h})\big),
$$
where the sum is over the half-edges $h$ at $v$.  Loops do not contribute to the Laplacian.
It has eigenvalues $0 = \lambda_{0}<\lambda_{1}\leq\cdots\leq \lambda_{V-1}$.  The product of
the nonzero eigenvalues is denoted $\text{\rm det}^{*}\Delta_{G}$, and
$$\zeta_{G}(s) = \sum_{j\geq 1}\lambda_{j}^{-s}$$
is the spectral zeta function of $G$. By Kirchhoff's
matrix-tree theorem 
(see \cite[Chapter~6]{Bi93})
$$
\tau(G) =\frac{1}{V} \text{\rm det}^{*}\Delta_{G}= \frac{1}{V}\prod\limits_{j=1}^{V-1}\lambda_{j}.
$$
More generally (see \cite{KC74, CS97}), the weighted forest counts \eqref{eq:cm} are the coefficients of
the polynomial
\begin{equation}\label{eq:matrix-forest}
\det\big(\Delta_{G}+x\big) = \sum\limits_{m=1}^{V}c_{m}(G)\,x^{m},
\end{equation}
which is the matrix-forest theorem, and
$$
c_{m}(G) = \sum\limits_{0\leq j_{1}<\cdots<j_{V-m}\leq V-1}
\lambda_{j_{1}}\cdots\lambda_{j_{V-m}},
$$
the elementary symmetric function of degree $V-m$ in the eigenvalues.
In the case $m=1$, one has $c_{1}(G) = \text{\rm det}^{*}\Delta_{G}$ and $w(\mathcal{F}) = V$ for every
spanning tree, so \eqref{eq:matrix-forest} contains Kirchhoff's matrix-tree theorem.

Since $\lambda_{0}=0$, only subsets omitting $\lambda_{0}$ contribute to $c_{m}(G)$, and
factoring out the product of all nonzero eigenvalues gives
\begin{equation}\label{eq:forest-zeta}
c_{m}(G) = \text{\rm det}^{*}\Delta_{G}\cdot
e_{m-1}\!\left(\lambda_{1}^{-1},\dots,\lambda_{V-1}^{-1}\right),
\end{equation}
where $e_{m-1}$ is the elementary symmetric polynomial of degree $m-1$.  By Newton's identities,
$e_{m-1}$ is a polynomial in the power sums
$\sum_{j\geq 1}\lambda_{j}^{-n} = \zeta_{G}(n)$, $1\leq n\leq m-1$.  This is the form in which
the forest counts are related to spectral zeta values in Corollary~\ref{cor:zeta-integers}.

\subsection{A generalized matrix-forest theorem}\label{sec:generalizedmatrixtree}

Assign to the finite connected graph $G$ a weight $w_{j}$ on each edge $\mathcal{E}_{j}$ and weights
$\{p_{v}\}$ on the vertices.  Let $L_{w}$ denote the weighted Laplacian,
$$
(L_{w}f)(v) = \sum_{h\,:\,v(h)=v} w_{e(h)}\big(f(v)-f(v_{h})\big),
$$
the sum being over the half-edges $h$ at $v$, loops again contributing nothing.  For a spanning forest $F$ of $G$ with components
$C$, write $p(C) = \sum_{v\in C}p_{v}$.

\begin{lem}\label{lem:matrix-forestnew}
With the above notation,
$$
\det\big(L_{w} + \text{\rm diag}(p)\big) = \sum_{F}\ \prod_{e\in F}w_{e}\ \prod_{C\in F}p(C),
$$
where the sum is over all spanning forests $F$ of $G$, with any number of components.
\end{lem}

\begin{proof}
By multilinearity of the determinant in the diagonal entries,
$$
\det\big(L_{w} + \text{\rm diag}(p)\big) = \sum_{S\subseteq V(G)}\Big(\prod_{v\in S}p_{v}\Big)\,
\det\big(L_{w}\big)_{\overline{S}},
$$
where $(L_{w})_{\overline{S}}$ is the principal submatrix with rows and columns indexed by
$V(G)\setminus S$, and the empty determinant is $1$.  By the all-minors matrix-tree theorem
\cite{Ch82}, $\det(L_{w})_{\overline{S}}$ is the sum of $\prod_{e\in F}w_{e}$ over spanning forests
$F$ each of whose components contains exactly one vertex of $S$, the root.  For such a forest write
$\text{\rm root}(C)$ for the unique vertex of $S$ lying in the component $C$.  Regrouping the double sum
by $F$, a forest $F$ appears once for each choice of one root in each of its components, with weight
$\prod_{e\in F}w_{e}\prod_{C}p_{\text{\rm root}(C)}$.  Summing over the choices of roots gives
$$
\sum_{\text{\rm roots}}\ \prod_{C}p_{\text{\rm root}(C)} = \prod_{C}\Big(\sum_{v\in C}p_{v}\Big)
= \prod_{C}p(C),
$$
each component choosing its root independently of the others.
\end{proof}

\begin{remark}\label{rem:matrix-forest-special}
Both sides of Lemma~\ref{lem:matrix-forestnew} are polynomials in the variables $w_{j}$ and $p_{v}$, so the
identity survives any substitution for them, in particular weights that are functions of a parameter.  That
is how it will be used in Section~\ref{subsec:DtN}, where the weights are meromorphic in $x$.  Two
specializations place the lemma among the classical statements.  Taking $p_{v} = x$ for every $v$ gives
$p(C) = x\,|C|$, so a forest with $m$ components contributes $x^{m}\prod_{C}|C|$, the number of ways of
choosing a root in each of its components, and with all $w_{j} = 1$ the lemma returns the matrix-forest
theorem \eqref{eq:matrix-forest}.  Taking $p_{v_{0}} = 1$ at a single vertex $v_{0}$ and $p_{v} = 0$
elsewhere annihilates every forest with a component avoiding $v_{0}$, so only the spanning trees survive, and
since $\det L_{w} = 0$ the lemma returns the weighted matrix-tree theorem in the grounded form
$$
\det\big(L_{w}\big)_{\overline{\{v_{0}\}}} = \sum_{T}\ \prod_{\mathcal{E}_{j}\in T}w_{j},
$$
the sum being over the spanning trees of $G$.  Between these two extremes the vertex weights are
unconstrained, and it is the unconstrained case that the metric graph needs, the $p_{v}$ carrying the
hyperbolic functions of Proposition~\ref{prop:Mlaplace}.
\end{remark}

\subsection{Determinant formulas from linear algebra}\label{subsec:determinantformulas}

We need two determinant formulas.

Schur's determinant formula is an expression for the determinant of a block matrix in terms of the Schur
complement.  Let $A$ and $D$ be square matrices, not necessarily of the same size, with $A$ invertible.
Multiplying out
$$
\begin{pmatrix} I&0\\-CA^{-1}&I\end{pmatrix}\begin{pmatrix}A&B\\C&D\end{pmatrix}
=\begin{pmatrix}A&B\\0&D-CA^{-1}B\end{pmatrix}
$$
and taking determinants gives
$$
\det\begin{pmatrix}A&B\\C&D\end{pmatrix} = \det A\cdot\det(D-CA^{-1}B).
$$
The matrix $S = D-CA^{-1}B$ is called the \emph{Schur complement} of $A$ in the full matrix; we refer to \cite{Sc17}.  The Schur complement of a Laplacian is again a
Laplacian plus a diagonal matrix, which is where the spanning forests come from.

The Cauchy--Binet formula is an expansion of the determinant of a product of rectangular matrices.  For
an $m\times n$ matrix $X$ and an $n\times m$ matrix $Y$ with $m\leq n$,
$$
\det(XY) = \sum_{S}\det(X_{S})\,\det(Y^{S}),
$$
the sum being over the $m$-element subsets $S$ of $\{1,\dots,n\}$, where $X_{S}$ consists of the columns
of $X$ indexed by $S$ and $Y^{S}$ of the rows of $Y$ indexed by $S$.

\subsection{The path graph}\label{subsec:pathgraphs}
Let $\mathcal P_n$ be the path graph on the $n$ vertices $0,1,\dots,n-1$, consecutive vertices joined by
an edge.  At an inner vertex, $0<m<n-1$, the equation $(\Delta_{\mathcal P_n}+u)f = 0$ is the recurrence
\begin{equation}\label{eq:recurrence}
f_{m+1}-(2+u)\,f_{m}+f_{m-1} = 0 .
\end{equation}
For $u\in\mathbb{C}$ let $\beta$ be defined by $\cosh\beta = 1+u/2$, and let
\begin{equation}\label{eq:D-eta}
D_{n}(u) = \frac{\sinh((n+1)\beta)}{\sinh\beta}\quad(n\geq -2),\qquad
\eta_{n}(u) = \frac{2\sinh(\beta/2)\sinh(n\beta/2)}{\cosh((n+1)\beta/2)}\quad(n\geq 0).
\end{equation}
Both are even functions of $\beta$, therefore independent of the choice of $\beta$ that satisfies the given condition.  Set $D_{-2} = -1$,
$D_{-1} = 0$ and $D_{0} = 1$.  Dividing the addition formula
$\sinh((m+2)\beta)+\sinh(m\beta) = 2\cosh\beta\,\sinh((m+1)\beta)$ by $\sinh\beta$ gives the recursion
\begin{equation}\label{eq:D-recursion}
D_{m+1} = (2+u)\,D_{m}-D_{m-1},
\end{equation}
so each $D_{n}$ is a polynomial in $u$ of degree $n$.  Let $T_{n}$ be the $n\times n$ tridiagonal matrix
with $2+u$ on the diagonal and $-1$ above and below it.
In fact, $T_{n}$ represents $\Delta+u$ on a path of $n$ vertices
with Dirichlet conditions at both ends.  Expanding along the first row gives
$\det T_{n} = (2+u)\det T_{n-1}-\det T_{n-2}$, which is again \eqref{eq:D-recursion}.  Since
$\det T_{0} = 1$ and $\det T_{1} = 2+u$ we obtain $\det T_{n} = D_{n}$.

Let
$$
\mathsf{T}(u) = \begin{pmatrix}2+u & -1\\ 1 & 0\end{pmatrix}
$$
be the transfer matrix of \eqref{eq:recurrence}, so that a sequence satisfies \eqref{eq:recurrence} if and
only if $(f_{m+1},f_{m})^{T} = \mathsf{T}(u)\,(f_{m},f_{m-1})^{T}$ for all $m$.

\begin{lem}\label{lem:transfer}
For every $m\geq 0$,
$$
\mathsf{T}(u)^{m} = \begin{pmatrix}D_{m} & -D_{m-1}\\ D_{m-1} & -D_{m-2}\end{pmatrix},\qquad
\operatorname{tr}\mathsf{T}(u)^{m} = D_{m}-D_{m-2} = 2\cosh(m\beta).
$$
The eigenvalues of $\mathsf{T}(u)$ are $e^{\pm\beta}$.  If $\sinh\beta\neq 0$, equivalently if $u\neq 0$
and $u\neq -4$, every solution of \eqref{eq:recurrence} has the form
$f_{m} = a\cosh(m\beta)+b\sinh(m\beta)$.  For every $u$, the solution with $f_{0} = 0$ is
$f_{m} = f_{1}D_{m-1}$.
\end{lem}

\begin{proof}
Let $\mathsf{A}_{m}$ denote the matrix on the right.  Then $\mathsf{A}_{0} = I$, since $D_{0} = 1$,
$D_{-1} = 0$ and $D_{-2} = -1$, and \eqref{eq:D-recursion} gives
$\mathsf{T}\,\mathsf{A}_{m} = \mathsf{A}_{m+1}$, hence $\mathsf{T}^{m} = \mathsf{A}_{m}$ for all
$m\geq 0$ by induction.  Taking traces and using the addition formula for $\sinh$,
$$
D_{m}-D_{m-2} = \frac{\sinh((m+1)\beta)-\sinh((m-1)\beta)}{\sinh\beta} = 2\cosh(m\beta).
$$
The characteristic polynomial of $\mathsf{T}(u)$ is $\lambda^{2}-2\cosh\beta\,\lambda+1$, whose roots are
$e^{\pm\beta}$.  They are distinct precisely when $\sinh\beta\neq 0$, and in that case $\mathsf{T}(u)$ is
diagonalizable, so the solutions of \eqref{eq:recurrence} are exactly the combinations of $e^{m\beta}$
and $e^{-m\beta}$, that is of $\cosh(m\beta)$ and $\sinh(m\beta)$.  Finally, the sequence $(D_{m-1})_{m}$
satisfies \eqref{eq:recurrence} with $D_{-1} = 0$ and $D_{0} = 1$, so for any solution with $f_{0} = 0$
the two sequences $f_{m}$ and $f_{1}D_{m-1}$ agree at $m = 0$ and $m = 1$, hence at every $m$.
\end{proof}

\begin{remark}\label{rem:edge-functions}
    The formulas of Theorem~\ref{mthm:convergence} and of Section~\ref{subsec:not-Phi} involve the functions
$\sinh(\sqrt{x}\,\ell)$, $\cosh(\sqrt{x}\,\ell)$ and $\tanh(\sqrt{x}\,\ell/2)$.  These functions are forced by the discretization,
and let us describe why.

In $G_{k}$ an edge of length $\ell$ becomes a path with about $k\ell$ vertices.  At its interior vertices the equation
$(\Delta_{G_{k}}+u)f = 0$ is the recurrence
$$
f_{m+1}-(2+u)f_{m}+f_{m-1} = 0,
$$
the discretization with mesh $1/k$ of the differential equation $-g''+xg = 0$ when $u = x/k^{2}$.
As we saw above its solutions are $\cosh(m\beta)$ and $\sinh(m\beta)$, and
the determinant of the path is the Chebyshev-type polynomial $D_{n}(u) = \sinh((n+1)\beta)/\sinh\beta$.  As $k\rightarrow\infty$ with
$u = x/k^{2}$, we have $k\beta\rightarrow\sqrt{x}$, and a vertex $m\approx ks$ corresponds to the point at distance $s$ along the edge.  So
$m\beta\rightarrow s\sqrt{x}$, and the discrete solutions converge to
$$
\frac{D_{m-1}(x/k^{2})}{k}\longrightarrow\frac{\sinh(\sqrt{x}\,s)}{\sqrt{x}},\qquad \cosh(m\beta)\longrightarrow\cosh(\sqrt{x}\,s),
$$
the fundamental solutions of $g'' = xg$ (Lemma~\ref{lem:continuum-limit}).  Evaluated at the ends of the edges, these
limits produce the building blocks of the metric forest function.  The Dirichlet determinant of a subdivided edge
gives $\sinh(\sqrt{x}\,\ell)/\sqrt{x}$, and the flux of the discrete harmonic extension through a component gives
$\sqrt{x}\tanh(\sqrt{x}\,\ell/2)$.  At the eigenvalues, $x = -\kappa^{2}$, the same functions become $\sin(\kappa s)/\kappa$ and
$\cos(\kappa s)$, the profiles of the eigenfunctions along the edges.

\end{remark}

\section{Metric graphs}\label{sec:metricgraphs}

\subsection{Metric graphs}
Let $G_0$ be a finite, connected (combinatorial) graph with edge set $E(G_{0})$. Let $\Lambda:E(G_{0}) \rightarrow  \mathbb R_{>0}$ be the function assigning length $\ell_j$ to the edge $\mathcal E_j$. By $\Lambda^{-1}$ we mean the function which assigns length $1/\ell_j$ to $\mathcal E_j$. The total length is $L = \sum_{j}\ell_{j}$, which can be considered as the volume of $\Gamma$.

Let $\Gamma = (G_{0},\Lambda)$ be the metric graph obtained by identifying each edge $\mathcal{E}_{j}$
of $G_{0}$ with the interval $[0,\ell_{j}]$. In case we fix an orientation, then $0$ corresponds to the initial point and
$\ell_{j}$ to the terminal point of the edge $\mathcal{E}_{j}$ with the selected orientation.  A function $f$ on $\Gamma$ is a collection
$(f_{j})$ of functions $f_{j}$ on $[0,\ell_{j}]$, and
$L^{2}(\Gamma) = \bigoplus_{j}L^{2}([0,\ell_{j}])$.

\subsection{The Laplacian of metric graphs}

As in \S\ref{subsec:graphsandspanningforests}, $h$ denotes a half-edge, $e(h)$ the edge containing it and
$v(h)$ its vertex.  For $f_{j}\in C^{1}([0,\ell_{j}])$, the outgoing derivative of $f$ along $h$ is
$\partial_{h}f = f_{j}'(0)$ if $h$ is the end of $\mathcal{E}_{j}$ at which $s=0$, and
$\partial_{h}f = -f_{j}'(\ell_{j})$ if $h$ is the other end.

Here $H^{r}([0,\ell])$ is the Sobolev space of functions on $[0,\ell]$ whose derivatives up to order $r$,
taken in the weak sense, belong to $L^{2}([0,\ell])$; see \cite[Chapter~8]{Br11}.  A function $f = (f_{j})$ on $\Gamma$
belongs to $\bigoplus_{j}H^{r}([0,\ell_{j}])$ when each component $f_{j}$ belongs to $H^{r}([0,\ell_{j}])$, with no
condition imposed at the vertices; the vertex conditions are imposed separately, as follows.  Since the edges
have finite length, $H^{2}([0,\ell])$ embeds into $C^{1}([0,\ell])$ \cite[Chapter~8]{Br11}, so that the values $f_{j}(0)$,
$f_{j}(\ell_{j})$ and the derivatives $\partial_{h}f$ are defined pointwise.

By definition, the Laplacian $\Delta_{\Gamma}$ with standard (Kirchhoff) vertex conditions acts by
$(\Delta_{\Gamma}f)_{j} = -f_{j}''$ on the domain of all $f\in\bigoplus_{j}H^{2}([0,\ell_{j}])$ that
are continuous at the vertices and satisfy $\sum_{h:\,v(h)=v}\partial_{h}f = 0$ at each vertex $v$.
It is self-adjoint and non-negative with compact resolvent; see, for example, \cite[Chapter~1]{BK13},
\cite[Section~3]{Kuc08} or \cite[Chapter~3]{Ku24}.  Its
quadratic form is $\sum_{j}\int_{0}^{\ell_{j}}|f_{j}'|^{2}\,ds$ on the space of continuous functions
in $\bigoplus_{j}H^{1}([0,\ell_{j}])$.  We write its eigenvalues, repeated according to
multiplicity, as
$$
0 = \mu_{0} < \mu_{1}\leq\mu_{2}\leq\cdots\rightarrow\infty.
$$
The frequencies are $k_{n} = \sqrt{\mu_{n}}$. The eigenvalue $0$ is simple because $\Gamma$ is connected.  Indeed, if $\Delta_{\Gamma}f=0$, then the
quadratic form vanishes, so each $f_{j}$ is constant, and continuity forces $f$ to be constant.

A vertex of degree two with standard conditions imposes only the continuity of $f$ and $f'$ across
it.  Hence inserting such a vertex into the interior of an edge does not change the operator
$\Delta_{\Gamma}$.

 The spectral zeta functions are $\zeta_{\Gamma}(s) = \sum_{j\geq 1}\mu_{j}^{-s}$ and
$\zeta_{\Gamma}(s;x) = \sum_{j\geq 0}(\mu_{j}+x)^{-s}$, the zeta-regularized determinants are $\text{\rm det}(\Delta_{\Gamma}+x)$ and
$\text{\rm det}^{*}\Delta_{\Gamma}$, and $Q(t) = \sum_{j\geq 0}e^{-\mu_{j}t}$ is the theta function (Section~\ref{sec:det}).
We omit the subscript $\zeta$ throughout.  The operator says which determinant is meant, since
$\text{\rm det}^{*}\Delta_{G}$ is the ordinary determinant of the Laplacian of a finite graph $G$.

A Laplacian on a different function space is obtained by imposing Dirichlet conditions at every vertex.  This
operator $\Delta_{D}$ acts by $(\Delta_{D}f)_{j} = -f_{j}''$ on those $f\in\bigoplus_{j}H^{2}([0,\ell_{j}])$
that vanish at the vertices.  Nothing couples the edges, so $\Delta_{D}$ is the direct sum of the Dirichlet
Laplacians of the edges, with eigenvalues $(\pi n/\ell_{j})^{2}$.  In particular it depends only on the edge
lengths and not on the combinatorics of $G_{0}$.

Both operators appear here because the regularized determinant of $\Delta_{\Gamma}+x$ factors as that of
$\Delta_{D}+x$ times the determinant of the Dirichlet-to-Neumann matrix.  The first factor carries the edges and
the second the vertices.  This factorization is due to \cite{Fr06}.  We prove it independently in
Corollary~\ref{cor:friedlander}, starting from Proposition~\ref{prop:DtN-properties}(4), by the same comparison
of logarithmic derivatives that Friedlander uses.  Remark~\ref{rem:friedlander} gives a second derivation of the
formula for $\text{\rm det}^{*}\Delta_{\Gamma}$, which instead comes out of the convergence of the
discrete determinants.

\subsection{The metric forest function}\label{subsec:not-Phi}

The basic functions attached to an edge and to a component of a spanning forest of $G_0$ are
$$
S_{j}(x)=S_{\mathcal E_j}(x) = S(x;\ell_{j}) = \frac{\sinh(\sqrt{x}\,\ell_{j})}{\sqrt{x}},\qquad
\Theta_{C}(x) = \sqrt{x}\sum_{h\text{ \rm at }C}\tanh\frac{\sqrt{x}\,\ell_{e(h)}}{2},
$$
where the second sum is over the half-edges $h$ of $G_{0}$ at the vertices of $C$, so that an edge with
both ends in $C$ contributes twice.  Recall the series expansions
$$
S_{j}(x) = \sum_{n\geq 0}\frac{\ell_{j}^{2n+1}}{(2n+1)!}\,x^{n},
$$
$$
\sqrt{x}\tanh\Big(\frac{\sqrt{x}\,\ell}{2}\Big) = \sum_{n\geq 1}\frac{2(2^{2n}-1)B_{2n}}{(2n)!}\,\ell^{2n-1}x^{n},
$$
with $B_{2n}$ denoting the Bernoulli numbers. Note that the square roots disappear and the former function is entire and the latter meromorphic in $x$.

From them we build the \emph{metric forest function}
\begin{equation}\label{eq:Phi}
\Phi(x) = \Phi_{(G_{0},\Lambda)}(x) := \sum_{F}\ \prod_{\mathcal{E}_{j}\notin F}S_{j}(x)\
\prod_{C\in F}\Theta_{C}(x),
\end{equation}
where $F$ runs over all spanning forests of $G_{0}$.

For $x\in\mathbb{C}$ we also define
$$
C_{j}(x) = C(x;\ell_{j})= \cosh(\sqrt{x}\,\ell_j),
$$
which is an entire function of $x$.

\subsection{The Dirichlet-to-Neumann matrix}\label{subsec:not-DtN}\label{subsec:DtN}

A central object in the proofs is the \emph{Dirichlet-to-Neumann matrix} $M(x)$ of $\Gamma$; see \cite[Section~7.2]{Kuc08}
and \cite{BK13, Ku24}, and \cite{PM99, ACDMT00} for the physics literature. Let
$$
Z = Z(\Gamma) = \{x\in\mathbb{C}: S_{j}(x)=0\text{ for some }j\}
= \big\{-(\pi n/\ell_{j})^{2}: n\geq 1,\ 1\leq j\leq E_{0}\big\},
$$
a discrete subset of the negative real axis.

\begin{lem}\label{lem:edge-solutions}
Let $x\notin Z$.  For every $a,b\in\mathbb{C}$ there is a unique solution of $g''(s) = xg(s)$ on
$[0,\ell_{j}]$ with $g(0)=a$ and $g(\ell_{j})=b$, namely
$$
g(s) = \frac{a\,S(x;\ell_{j}-s) + b\,S(x;s)}{S_{j}(x)}.
$$
Furthermore, $g'(0) = (b-aC_{j}(x))/S_{j}(x)$ and
$-g'(\ell_{j}) = (a-bC_{j}(x))/S_{j}(x)$.
\end{lem}

\begin{proof}
The solutions of $g''=xg$ form a two-dimensional space, spanned by $S(x;\cdot)$ and $C(x;\cdot)$, whose
Wronskian is $1$.  Since $S(x;\ell_{j}-s) = S_{j}(x)C(x;s)-C_{j}(x)S(x;s)$ and $S_{j}(x)\neq 0$, the
functions $S(x;s)$ and $S(x;\ell_{j}-s)$ span it as well.  The stated $g$ is a solution with the
required boundary values.  Uniqueness holds because a solution vanishing at both ends is a multiple
of $S(x;s)$ that vanishes at $s=\ell_{j}$, hence is zero since $S_{j}(x)\neq 0$.  Since $\partial_{s}S(x;s) = C(x;s)$, we have
$g'(s) = \big(-a\,C(x;\ell_{j}-s) + b\,C(x;s)\big)/S_{j}(x)$, and evaluating at $s=0$ and $s=\ell_{j}$
gives the stated outgoing derivatives.
\end{proof}

For $x\notin Z$ and $\varphi\in\mathbb{C}^{V(G_{0})}$, let $f_{\varphi} = f_{\varphi}(x)$ be the function on
$\Gamma$ which on each edge is the solution of $g''=xg$ from Lemma~\ref{lem:edge-solutions} with the
values of $\varphi$ at the endpoints.  Define the $V_{0}\times V_{0}$ matrix $M(x)$ by
$$
\big(M(x)\varphi\big)_{w} = -\sum_{h:\,v(h)=w}\partial_{h}f_{\varphi}.
$$
By Lemma~\ref{lem:edge-solutions}, $M(x)$ is symmetric with entries
$$
M(x)_{ww} = \sum_{h:\,v(h)=w}\frac{C_{e(h)}(x)}{S_{e(h)}(x)}
-\sum_{\text{\rm loops }\mathcal{E}_{j}\text{ \rm at }w}\frac{2}{S_{j}(x)},\qquad
M(x)_{wv} = -\sum_{\mathcal{E}_{j}=\{w,v\}}\frac{1}{S_{j}(x)}\quad(w\neq v).
$$
Here the first sum is over all edges with the two orientations originating from $w$, the second sum is over the unoriented loops at $w$. The non-diagonal sum is over unoriented edges between $w$ and $v$.
Since $C_{j}/S_{j} = \sqrt{x}\coth(\sqrt{x}\ell_{j})$ and $1/S_{j} = \sqrt{x}/\sinh(\sqrt{x}\ell_{j})$,
this is the matrix $M(x)$ of Section~\ref{sec:forest-expansion}, so that
$\Phi(x) = \prod_{j}S_{j}(x)\cdot\det M(x)$ for $x\notin Z$.  The entries of $M(x)$ are
holomorphic on $\mathbb{C}\setminus Z$ and real for real $x\notin Z$.

We now consider the combinatorial graph $G_0$ and assign certain weights to its vertices and edges (\emph{cf.} section \ref{sec:generalizedmatrixtree}). Let $w_j(x)=1/S_j(x)$ be the weight on $\mathcal E_j$ and vertices $v$ get the weight
$$
p_v(x)=\sum_{h:v(h)=v}\frac {C_{e(h)}(x)-1}{S_{e(h)}(x)}=\sqrt{x}\sum_{h:v(h)=v}\tanh \frac{\sqrt{x}\ell_{e(h)}}{2}.
$$
With these weights, every component $C$ of a spanning forest satisfies
$$
\Theta_C(x)=\sum_{v\in V(C)}p_v(x).
$$
Moreover, we then have the following observation.
\begin{prop} \label{prop:Mlaplace}
   Let $x\notin Z$. Then we have that
   $$
   M(x)=L_{w(x)}+\mathrm{diag}(p_v(x)),
   $$
   where $L_w$ is the weighted Laplacian of the finite combinatorial graph $G_0$.
\end{prop}

\begin{proof}
    Take a non-loop edge $\mathcal E_j=\{ v,w\}$. Its contribution to the matrix $M(x)$ is
    $$
    \frac{1}{S_j}\begin{pmatrix}C_j & -1 \\ -1 & C_j\end{pmatrix}=\frac{1}{S_j}\begin{pmatrix}1 & -1 \\ -1 & 1\end{pmatrix}+\frac{C_j-1}{S_j}\begin{pmatrix}1 & 0 \\ 0 & 1\end{pmatrix}
    $$
    in the rows and columns corresponding to the two endpoints of the edge. Summing over the edges gives $L_w+\textrm{diag}(p)$. A loop at $v$ also fits this decomposition, since it contributes nothing to $L_w$, while its two half-edges contribute $2(C_j-1)/S_j$ to $p_v$, which is $2C_j/S_j - 2/S_j$.
\end{proof}

From this and Lemma \ref{lem:matrix-forestnew} we have the following result.

\begin{prop}\label{prop:Phi-DtN}
    With the above notation, for $x\in \mathbb C \setminus Z$,
    $$
    \Phi(x)=\Big( \prod_j S_j(x) \Big)\det M(x).
    $$
\end{prop}

\begin{proof}
    By Proposition~\ref{prop:Mlaplace}, the generalized matrix-tree formula of
    Lemma~\ref{lem:matrix-forestnew}, applied with the weights defined above, gives
$$
\det M(x)= \sum_{F}\ \prod_{\mathcal E_j\in F}\frac{1}{S_j(x)}\ \prod_{C\in F}\Theta_C (x).
$$
Multiplying both sides of this equation by $\prod_j S_j(x)$ finishes the proof in view of the definition (\ref{eq:Phi}) of the metric forest function $\Phi(x)$.
\end{proof}

\begin{remark}
    Later we will see that $\Phi(x)$ is in fact entire.
\end{remark}

\begin{lem}\label{lem:DtN-kernel}
Let $x\notin Z$.  The map $\varphi\mapsto f_{\varphi}(x)$ is a linear isomorphism from $\ker M(x)$ onto
$\ker(\Delta_{\Gamma}+x)$.  In particular, $\det M(x) = 0$ if and only if $-x$ is an eigenvalue of
$\Delta_{\Gamma}$, and then $\dim\ker M(x)$ is its multiplicity.
\end{lem}

\begin{proof}
By construction $f_{\varphi}$ is continuous, lies in $\bigoplus_{j}H^{2}$ and satisfies
$f_{\varphi}''=xf_{\varphi}$ on each edge.  It satisfies the Kirchhoff conditions if and only if
$M(x)\varphi=0$.  Hence the map sends $\ker M(x)$ into $\ker(\Delta_{\Gamma}+x)$, and it is injective
because $\varphi$ is the vector of vertex values of $f_{\varphi}$.  Conversely, if
$(\Delta_{\Gamma}+x)f=0$, let $\varphi$ be the vector of vertex values of $f$, which is well defined by
continuity.  On each edge, $f$ solves $g''=xg$ with the endpoint values of $\varphi$, so
$f = f_{\varphi}$ by the uniqueness in Lemma~\ref{lem:edge-solutions}.  Since $f$ satisfies the Kirchhoff
conditions, $M(x)\varphi=0$.
\end{proof}

\begin{lem}\label{lem:energynew}
Let $x\in\mathbb{R}\setminus Z$ and $\varphi,\psi\in\mathbb{R}^{V(G_{0})}$.  Then
$$
\psi^{T}M(x)\varphi = \sum_{j}\int_{0}^{\ell_{j}}\big(f_{\psi}'f_{\varphi}' + x f_{\psi}f_{\varphi}\big)\,ds
\qquad\text{and}\qquad
\varphi^{T}M'(x)\varphi = \|f_{\varphi}(x)\|_{L^{2}(\Gamma)}^{2},
$$
where $M'(x) = \frac{d}{dx}M(x)$.
\end{lem}

\begin{proof}
Since $f_{\varphi}'' = xf_{\varphi}$ on each edge, we have
$(f_{\psi}f_{\varphi}')' = f_{\psi}'f_{\varphi}' + xf_{\psi}f_{\varphi}$.
Integrating over each edge and summing,
$$
\sum_{j}\int_{0}^{\ell_{j}}\big(f_{\psi}'f_{\varphi}' + xf_{\psi}f_{\varphi}\big)\,ds
= \sum_{j}\Big[f_{\psi}f_{\varphi}'\Big]_{0}^{\ell_{j}}
= -\sum_{h}\psi_{v(h)}\,\partial_{h}f_{\varphi} = \psi^{T}M(x)\varphi.
$$
For the second identity, fix $\varphi$ and write $f=f_{\varphi}(x)$ and
$\dot f = \partial_{x}f_{\varphi}(x)$, which
exists by the explicit formula in Lemma~\ref{lem:edge-solutions}.  Differentiating $f''=xf$ in $x$
gives $\dot f'' = x\dot f + f$, and $\dot f$ vanishes at all vertices because the vertex values of
$f_{\varphi}(x)$ do not depend on $x$.  Differentiating the first identity with $\psi=\varphi$,
$$
\varphi^{T}M'(x)\varphi = \sum_{j}\int_{0}^{\ell_{j}}\big(2f'\dot f' + f^{2} + 2xf\dot f\big)\,ds.
$$
Integration by parts gives $\int f'\dot f' = [f'\dot f] - \int f''\dot f = -x\int f\dot f$ on each edge,
since $\dot f$ vanishes at the endpoints.  Hence
$\varphi^{T}M'(x)\varphi = \sum_{j}\int_{0}^{\ell_{j}}f^{2}\,ds$.
\end{proof}

\begin{prop}\label{prop:DtN-properties}
The Dirichlet-to-Neumann matrix $M(x)$ has the following properties.
\begin{enumerate}
\item Its entries are holomorphic on $\mathbb{C}\setminus Z$, and $M(x)$ is real symmetric for real $x\notin Z$.
\item It extends holomorphically to $x = 0$, and $M(0)$ is the weighted combinatorial Laplacian of $G_{0}$ with
edge weights $1/\ell_{j}$, with loops not contributing.
\item For real $x>0$, $M(x)$ is positive definite, and $M'(x)$ is positive definite.
\item For $x\notin Z$, $\Phi(x) = \prod_{j}S_{j}(x)\det M(x)$, and $\det M(x) = 0$ exactly when $-x$ is an eigenvalue of
$\Delta_{\Gamma}$ 
\item For every $\varepsilon>0$ there are $c,r_{0}>0$ such that, for $|\arg x|\leq\pi-\varepsilon$ and $|x|\geq r_{0}$,
$$
M(x) = \sqrt{x}\,\text{\rm diag}(d_{v})+O\big(|x|^{1/2}e^{-c\sqrt{|x|}}\big).
$$
\end{enumerate}
\end{prop}

\begin{proof}
Part (1) is clear from the explicit entries.  For (2), as $x\rightarrow 0$ we have
$\sqrt{x}\coth(\sqrt{x}\,\ell)\rightarrow 1/\ell$ and $\sqrt{x}/\sinh(\sqrt{x}\,\ell)\rightarrow 1/\ell$, and a loop contributes
$2\sqrt{x}\tanh(\sqrt{x}\,\ell/2)\rightarrow 0$.  These functions are even in $\sqrt{x}$, hence holomorphic in $x$ near $0$.
For (3), Lemma~\ref{lem:energynew} gives
$\varphi^{T}M(x)\varphi = \sum_{j}\int(|f_{\varphi}'|^{2}+x|f_{\varphi}|^{2})\,ds>0$ and
$\varphi^{T}M'(x)\varphi = \|f_{\varphi}\|^{2}>0$ for real
$\varphi\neq 0$.  Part (4) was shown above.  For (5), write $\omega = \sqrt{x}$ with $\text{\rm Re}\,\omega\geq |\omega|\sin(\varepsilon/2)$.
Then $\coth(\omega\ell) = 1+O(e^{-2\text{\rm Re}(\omega)\ell})$ and $1/\sinh(\omega\ell) = O(e^{-\text{\rm Re}(\omega)\ell})$, while a loop contributes
$2\omega\tanh(\omega\ell/2) = 2\omega+O(|\omega|e^{-\text{\rm Re}(\omega)\ell})$.  Summing over the half-edges at each vertex gives $d_{v}\sqrt{x}$ on the
diagonal, up to the stated error.
\end{proof}

  The coefficient of $x^{m}$ in $\Phi(x)$ is denoted $P_{m}$; it is the limiting
normalized $m$-spanning forest count (Theorem~\ref{thm:Phi-limit}), and the ratios $e_{m} = P_{m+1}/P_{1}$ are the elementary
symmetric functions of the numbers $1/\mu_{j}$.

In the explicit formulas of Section~\ref{sec:explicit} we use the metric size $\ell^{*}(C) = \frac12\sum_{h\text{ \rm at }C}\ell_{e(h)}$ of a
component $C$, the product $w_{F} = \prod_{\mathcal{E}_{j}\notin F}\ell_{j}$ of the lengths of the edges outside a forest $F$, and the power sums
of lengths $\sigma_{p}(F)$ and $\Lambda_{p}(C)$.  The period matrix $Q_{\Gamma}$ is the Gram matrix of a basis of $H_{1}(G_{0},\mathbb{Z})$ in the
inner product $\langle\mathcal{E}_{i},\mathcal{E}_{j}\rangle = \delta_{ij}\ell_{j}$; its determinant is $\det Q_{\Gamma} = \sum_{T}\prod_{\mathcal{E}_{j}\notin T}\ell_{j}$
(Section~\ref{subsec:CM} and Lemma~\ref{lem:period}).

\section{Combinatorial graph approximations of metric graphs}\label{sec:graphapproximations}
\subsection{Discrete graph approximations to $G_{0}$}

We now construct a sequence of graphs $\{G_{k}\}$, $k \in \Z_{\geq 0}$, starting with $G_{0}$
as above.  For each edge $\mathcal{E}_{j}$ of $G_{0}$, $1\leq j\leq E_{0}$, choose a positive real
number $\ell_{j}$, and write $\Lambda=\{\ell_{j}\}$ and $L=\sum_{j}\ell_{j}$.

Choose non-negative integers $n_{j,k}$ with $n_{j,0} = 0$ and $n_{j,k}/k \rightarrow \ell_{j}$
as $k \rightarrow \infty$.  We form $G_{k}$ by replacing each edge $\mathcal{E}_{j}$ of $G_{0}$
with the discrete interval $[0,n_{j,k}+1]$, with vertices at the integers, where $0$ is identified with the
initial point of $\mathcal{E}_{j}$ and $n_{j,k}+1$ with its terminal point.  In other words, we
subdivide $\mathcal{E}_{j}$ by inserting $n_{j,k}$ new vertices, each of degree two, so that
$\mathcal{E}_{j}$ becomes a path of $n_{j,k}+1$ edges, each of weight one.

It is convenient to view $G_{k}$ as a discretization of the metric graph $(G_{0},\Lambda)$ with
mesh $1/k$.  Set
$$
\ell_{j,k} := \frac{n_{j,k}+1}{k},\qquad \Lambda_{k}:=\{\ell_{j,k}\},\qquad L_{k}:=\sum_{j=1}^{E_{0}}\ell_{j,k},
$$
so that $\ell_{j,k}\rightarrow\ell_{j}$ and $L_{k}\rightarrow L$ as $k\rightarrow\infty$.  If each
unit edge of $G_{k}$ is given length $1/k$, then $G_{k}$ is the set of points of mesh $1/k$ on the
metric graph $(G_{0},\Lambda_{k})$.

\begin{figure}[!ht]
\centering
\begin{tikzpicture}[v/.style={circle,fill,inner sep=1.9pt}, sv/.style={circle,fill,inner sep=0.85pt},
  e/.style={thick}, lab/.style={font=\small}]
\begin{scope}
  \draw[e] (0,0) -- node[above,lab] {$\mathcal{E}_{j}$, of length $\ell_{j}$} (6,0);
  \node[v,label=below:{$u$}] at (0,0) {};
  \node[v,label=below:{$w$}] at (6,0) {};
\end{scope}
\draw[-{Stealth[length=3mm]},thick] (6.7,0) -- node[above,lab] {subdivide} (8.3,0);
\begin{scope}[shift={(9.0,0)}]
  \draw[e] (0,0) -- (6,0);
  \foreach \i in {1,...,6} \node[sv] at ({\i*6/7},0) {};
  \node[v,label=below:{$u$}] at (0,0) {};
  \node[v,label=below:{$w$}] at (6,0) {};
  \draw[decorate,decoration={brace,amplitude=4pt}] (0,0.16) -- node[above=4pt,lab] {$1/k$} ({6/7},0.16);
  \draw[decorate,decoration={brace,amplitude=4pt,mirror}] (0,-0.9) --
        node[below=4pt,lab] {$\ell_{j,k} = (n_{j,k}+1)/k$} (6,-0.9);
\end{scope}
\end{tikzpicture}
\caption{The edge $\mathcal{E}_{j}$ of $G_{0}$ and the path replacing it in $G_{k}$, here with $n_{j,k} = 6$ interior
vertices of degree two and $n_{j,k}+1$ edges.  Giving each of these edges the length $1/k$ turns the path into an
interval of length $\ell_{j,k}$, which tends to $\ell_{j}$ as $k\rightarrow\infty$.}
\label{fig:subdivision}
\end{figure}
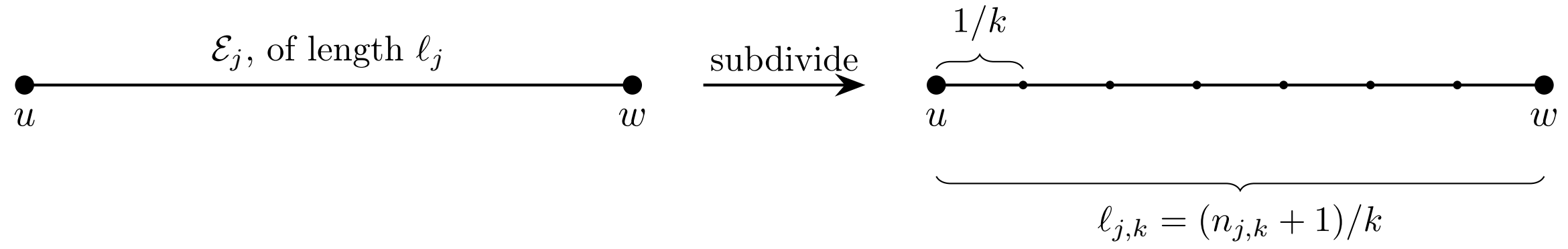

See Figure~\ref{fig:subdivision}.

Let
$$
H_{k} := \sum_{j=1}^{E_{0}}n_{j,k}.
$$
The graph $G_{k}$ has $V_{k}:=V_{0} + H_{k}$ vertices and $E_{k} := E_{0}+H_{k}$ edges.  In
particular,
$$
E_{k} - V_{k} = E_{0}-V_{0}\quad\text{for all $k$,}\qquad\text{and}\qquad
V_{k} = kL_{k} - (E_{0}-V_{0}).
$$

\subsection{Spanning trees of $G_{k}$}\label{subsec:trees-Gk}

\begin{lem}\label{lem:trees}
For all $k\geq 0$,
$$
\tau(G_{k}) = \sum\limits_{\text{\rm spanning trees $T_{0}$ of $G_{0}$}}\ \
\prod\limits_{\mathcal{E}_{j}\notin T_{0}}\left(n_{j,k}+1\right).
$$
\end{lem}

\begin{proof}
For each edge $\mathcal{E}_{j}$ of $G_{0}$, let $P_{j}$ denote the path of $n_{j,k}+1$ edges in
$G_{k}$ which replaces $\mathcal{E}_{j}$.

Let $T$ be a spanning tree of $G_{k}$.  Then $T$ omits at most one edge from each path $P_{j}$.
Indeed, if two edges of $P_{j}$ were omitted, then the vertices of $P_{j}$ lying strictly between
them would form a nonempty set of degree-two vertices of $G_{k}$ with no edge of $T$ connecting
it to the rest of $G_{k}$, contradicting connectedness.  Let $S\subseteq\{1,\dots,E_{0}\}$ be the set
of indices $j$ for which $T$ omits an edge of $P_{j}$.  Since $T$ omits
$E_{k}-V_{k}+1 = E_{0}-V_{0}+1$ edges, we have $|S| = E_{0}-V_{0}+1$, so the subgraph $T_{0}$ of
$G_{0}$ with edge set $\{\mathcal{E}_{j}:j\notin S\}$ has $V_{0}-1$ edges.  Any path in $T$
between two vertices of $G_{0}$ traverses complete paths $P_{j}$ with $j\notin S$, so $T_{0}$ is
connected.  Hence $T_{0}$ is a spanning tree of $G_{0}$.

Conversely, given a spanning tree $T_{0}$ of $G_{0}$ and, for each $\mathcal{E}_{j}\notin T_{0}$, a
choice of one of the $n_{j,k}+1$ edges of $P_{j}$ to omit, the resulting subgraph of $G_{k}$ is
connected and has $E_{k}-(E_{0}-V_{0}+1) = V_{k}-1$ edges, so it is a spanning tree.  The two
constructions are inverse to each other, which proves the lemma.
\end{proof}

Since $n_{j,k}+1 = k\ell_{j,k}$ and a spanning tree of $G_{0}$ has $V_{0}-1$ edges, so that
its complement has $E_{0}-V_{0}+1$ edges, Lemma~\ref{lem:trees} can be written as
\begin{equation}\label{eq:tau-exact}
\tau(G_{k}) = k^{E_{0}-V_{0}+1}\sum\limits_{T_{0}}\ \prod\limits_{\mathcal{E}_{j}\notin T_{0}}\ell_{j,k}
= k^{E_{0}-V_{0}+1}\Big(\prod_{j=1}^{E_{0}}\ell_{j,k}\Big)\,\tau(G_{0},\Lambda_{k}^{-1}),
\end{equation}
where, for a set of lengths $\Lambda=\{\ell_{j}\}$, we define the weighted spanning tree number
$$
\tau(G_{0},\Lambda^{-1}) := \sum\limits_{\text{\rm spanning trees $T_{0}$ of $G_{0}$}}\ \
\prod\limits_{\mathcal{E}_{j}\in T_{0}} \ell_{j}^{-1}.
$$
By the weighted matrix-tree theorem, $\tau(G_{0},\Lambda^{-1})$ is any cofactor of the Laplacian
of $G_{0}$ with edge weights $\ell_{j}^{-1}$.  Note that \eqref{eq:tau-exact} is an identity, not
only an asymptotic statement.

Combining \eqref{eq:tau-exact} with Kirchhoff's theorem
$\text{\rm det}^{*}\Delta_{G_{k}} = V_{k}\,\tau(G_{k})$ and $V_{k} = kL_{k}-(E_{0}-V_{0})$, we
obtain the exact formula
\begin{equation}\label{eq:logdet-exact}
\log (\text{\rm det}^{*}\Delta_{G_{k}}) = (E_{0}-V_{0}+1) \log k
+ \log\Big(\prod_{j}\ell_{j,k}\Big) + \log \tau(G_{0},\Lambda_{k}^{-1})
+ \log\big(kL_{k}-E_{0}+V_{0}\big).
\end{equation}

\begin{prop}\label{prop:logdet-asymptotics}
As $k\rightarrow\infty$,
$$
\log (\text{\rm det}^{*}\Delta_{G_{k}}) = (E_{0}-V_{0}+2) \log k + \log \tau(G_{0},\Lambda^{-1})
+ \log\Big(\prod_{j}\ell_{j}\Big) + \log L + o(1).
$$
If $\ell_{j,k} = \ell_{j}+O(1/k)$ for each $j$, for example if $n_{j,k}+1 = \lfloor k\ell_{j}\rfloor$,
then the error term is $O(1/k)$.
\end{prop}

\begin{proof}
Each term on the right-hand side of \eqref{eq:logdet-exact} depends continuously on
$\Lambda_{k}\in(0,\infty)^{E_{0}}$, and
$\log(kL_{k}-E_{0}+V_{0}) = \log k + \log L_{k} + O(1/k)$.  Since $\Lambda_{k}\rightarrow\Lambda$,
the first assertion follows.  The terms are smooth in $\Lambda_{k}$, which gives the second.
\end{proof}

Lemma~\ref{lem:trees} expresses the behaviour of the Kirchhoff polynomial under subdivision of edges.  A path of
$n+1$ unit edges contributes like a single edge of resistance $n+1$.  Proposition~\ref{prop:logdet-asymptotics} is
the analogue for metric graphs of the asymptotics of \cite{CJK10} for discrete tori degenerating to real tori; the
asymptotics of the number of spanning trees in that setting are treated in \cite{CJK12}.

\begin{remark}\label{rem:friedlander}
Proposition~\ref{prop:logdet-asymptotics} determines the growth of $\text{\rm det}^{*}\Delta_{G_{k}}$ by a purely combinatorial
argument, with constant term $\log\big(L\prod_{j}\ell_{j}\,\tau(G_{0},\Lambda^{-1})\big) = \log P_{1}$.  In Section~\ref{sec:det} the same
quantity is identified analytically.  The characteristic polynomials of the $G_{k}$ converge to $\Phi$
(Theorem~\ref{mthm:convergence}), and $\Phi$ equals $\det(\Delta_{\Gamma}+x)$ up to the factor $\prod_{v}d_{v}/2^{E_{0}}$
(Theorem~\ref{thm:det}).  Comparing the two descriptions proves the formula of \cite{Fr06}, namely that
$$
\text{\rm det}^{*}\Delta_{\Gamma} =
\frac{2^{E_{0}}\,L\,\prod_{j}\ell_{j}\;\tau(G_{0},\Lambda^{-1})}{\prod_{v\in V(G_{0})}d_{v}},
$$
as a corollary of the discrete approximation, and it also turns the proposition into the asymptotic formula for
$\log\text{\rm det}^{*}\Delta_{G_{k}}$ of Theorem~\ref{mthm:det}; see Corollary~\ref{cor:main}.
\end{remark}

\section{An exact forest expansion on \texorpdfstring{$G_{k}$}{Gk} and its limit}\label{sec:forest-expansion}

In this section we prove the convergence of the characteristic polynomials of the discrete graphs $G_{k}$ to the metric
forest function $\Phi$; the point of view is explained in Section~\ref{subsec:not-viewpoint}.  Throughout this section we assume
$k$ is large enough that $n_{j,k}\geq 1$ for all $j$.

\subsection{The exact identity}

We now prove the following, building on subsection \ref{subsec:pathgraphs}, where $D_n$ and $\eta_n$ are defined in (\ref{eq:D-eta}), and writing $H(G_{0})$ for the set of half-edges of $G_{0}$ as in subsection \ref{subsec:graphsandspanningforests}. Note that the second product expression in the following proposition only depends on the vertex set $V(C)$ of a subtree $C$ and not its edge set.

\begin{prop}\label{prop:exact-forest}
For all $u$,
$$
\det\big(\Delta_{G_{k}}+u\big) = \sum_{F}\ \prod_{\mathcal{E}_{j}\notin F}D_{n_{j,k}}(u)\
\prod_{C\in F}\Big(u\,|V(C)| + \sum_{\substack{h\in H(G_{0})\\ v(h)\in V(C)}}\eta_{n_{e(h),k}}(u)\Big),
$$
where $F$ runs over all spanning forests of $G_{0}$ and $C$ the connected components of $F$.
\end{prop}

\begin{proof}
Order the vertices of $G_{k}$ with the interior vertices of the paths $P_{j}$ first and the vertices
of $G_{0}$ last.
$$
\Delta_{G_{k}}+u = \begin{pmatrix}T(u)&B\\ B^t & R(u)\end{pmatrix},
$$
where $R(u)$ is the block on the vertices of $G_{0}$ and $B$ records the edges between the two sets.
The interior block, $T(u)$, is block diagonal, with one tridiagonal
block $T_{j}$ of size $n_{j,k}$ per edge, having $u+2$ on the diagonal and $-1$ off the diagonal.
As shown in subsection \ref{subsec:pathgraphs} its determinant is $\prod_{j}D_{n_{j,k}}(u)$, and by Schur's determinant formula, see subsection \ref{subsec:determinantformulas},
$$
\det\big(\Delta_{G_{k}}+u\big) = \prod_{j}D_{n_{j,k}}(u)\cdot\det S_{k}(u),
$$
where $S_{k}(u)$ is the Schur complement onto $V(G_{0})$.
By the cofactor formula for the inverse of the triangular matrices, $(T^{-1})_{11} = (T^{-1})_{nn} = D_{n-1}/D_{n}$,
because deleting the first row and column leaves $T_{n-1}$.  Similarly $(T^{-1})_{1n} = (T^{-1})_{n1} = 1/D_{n}$,
because the $(n,1)$ minor is the determinant of a lower triangular matrix with $-1$ on the diagonal,
equal to $(-1)^{n-1}$, and $(-1)^{n+1}(-1)^{n-1} = 1$.  Here $n = n_{j,k}$.  Each vertex $v$ of $G_{0}$ is joined by an edge of
weight one to the first or last interior vertex of each path at $v$, so
$$
S_{k}(u)_{vv} = u + \sum_{h\text{ \rm at }v}\Big(1-\frac{D_{n-1}}{D_{n}}\Big)
- \sum_{\text{\rm loops at }v}\frac{2}{D_{n}},\qquad
S_{k}(u)_{vw} = -\sum_{\mathcal{E}_{j}=\{v,w\}}\frac{1}{D_{n_{j,k}}}\quad (v\neq w),
$$
where in each term $n = n_{e(h),k}$.  Therefore $S_{k}(u) = L_{w}+\text{\rm diag}(p)$, with edge
weights $w_{j} = 1/D_{n_{j,k}}(u)$ and vertex weights
$$
p_{v} = u + \sum_{h\text{ \rm at }v}\Big(1-\frac{D_{n-1}+1}{D_{n}}\Big).
$$
We claim that $1-(D_{n-1}+1)/D_{n} = \eta_{n}$.  Indeed,
$$
D_{n}-D_{n-1}-1 = \frac{\sinh((n+1)\beta)-\sinh(n\beta)-\sinh\beta}{\sinh\beta}.
$$
By the identities $\sinh A-\sinh B = 2\cosh\frac{A+B}{2}\sinh\frac{A-B}{2}$ and
$\sinh\beta = 2\sinh\frac{\beta}{2}\cosh\frac{\beta}{2}$, we have
$$
\sinh((n+1)\beta)-\sinh(n\beta) = 2\cosh\frac{(2n+1)\beta}{2}\,\sinh\frac{\beta}{2}.
$$
Hence the numerator equals
$$
2\sinh\frac{\beta}{2}\Big(\cosh\frac{(2n+1)\beta}{2}-\cosh\frac{\beta}{2}\Big)
= 4\sinh\frac{\beta}{2}\,\sinh\frac{(n+1)\beta}{2}\,\sinh\frac{n\beta}{2},
$$
by $\cosh A-\cosh B = 2\sinh\frac{A+B}{2}\sinh\frac{A-B}{2}$.  Dividing by
$D_{n}\sinh\beta = \sinh((n+1)\beta) = 2\sinh\frac{(n+1)\beta}{2}\cosh\frac{(n+1)\beta}{2}$ gives
$$
\frac{D_{n}-D_{n-1}-1}{D_{n}} = \frac{2\sinh(\beta/2)\sinh(n\beta/2)}{\cosh((n+1)\beta/2)} = \eta_{n},
$$
as claimed.  Lemma~\ref{lem:matrix-forestnew} now gives
$$
\det S_{k}(u) = \sum_{F}\prod_{\mathcal{E}_{j}\in F}\frac{1}{D_{n_{j,k}}(u)}\ \prod_{C}p(C),
$$
and multiplying by $\prod_{j}D_{n_{j,k}}(u)$ completes the proof.
\end{proof}

Setting $u=0$ gives $\beta = 0$, $D_{n}(0) = n+1$ and $\eta_{n}(0) = 0$.  Since $\eta_{n}(u) = O(u)$,
only forests with one component contribute to the coefficient of $u$, and one recovers
Lemma~\ref{lem:trees} together with Kirchhoff's theorem.

\begin{lem}\label{lem:order}
Let $A(x)$ be an $N\times N$ matrix of functions holomorphic near $x_{0}\in\mathbb{R}$ which is real
symmetric for real $x$.  Let $K = \ker A(x_{0})$ have dimension $m$, and suppose that
$v^{T}A'(x_{0})v>0$ for all nonzero $v\in K$.  Then $\det A(x)$ has a zero of order exactly $m$ at
$x_{0}$.
\end{lem}

\begin{proof}
After conjugating by a real orthogonal matrix, we may assume $K$ is spanned by the first $m$
standard basis vectors.  Since $A(x_{0})$ is symmetric, it then has the block form
$\text{\rm diag}(0,B)$ with $B$ invertible.  Write
$A(x) = \begin{pmatrix}A_{11}(x) & A_{12}(x)\\ A_{21}(x) & A_{22}(x)\end{pmatrix}$ accordingly.  Near
$x_{0}$, the block $A_{22}(x)$ is invertible, and
$$
\det A(x) = \det A_{22}(x)\cdot\det\big(A_{11}(x)-A_{12}(x)A_{22}(x)^{-1}A_{21}(x)\big).
$$
Since $A_{11}(x_{0})$, $A_{12}(x_{0})$ and $A_{21}(x_{0})$ vanish, the second factor is the
determinant of $(x-x_{0})\big(A_{11}'(x_{0}) + O(x-x_{0})\big)$.  The matrix $A_{11}'(x_{0})$
represents the form $v\mapsto v^{T}A'(x_{0})v$ on $K$, so it is positive definite.  Hence
$\det A(x) = (x-x_{0})^{m}\big(\det B\cdot\det A_{11}'(x_{0}) + O(x-x_{0})\big)$ with a nonzero
leading coefficient.
\end{proof}


\subsection{Taking the scaling limit}\label{subsec:not-hyperbolic}

The next lemma makes precise the statement that, under the scaling $u = x/k^{2}$, the discrete solutions of
\eqref{eq:recurrence} converge to the hyperbolic solutions of $g'' = xg$, which are spanned by the functions
$S(x;s)$ and $C(x;s)$ of Section~\ref{subsec:not-Phi}.

\begin{lem}\label{lem:continuum-limit}
Let $u = x/k^{2}$, and let $m = m_{k}(s)$ be integers with $m_{k}(s)/k\rightarrow s$ uniformly for $s$ in compact subsets of
$[0,\infty)$.  Then, uniformly for $x$ in compact subsets of $\mathbb{C}$ and $s$ in compact subsets of $[0,\infty)$,
$$
\frac{D_{m-1}(x/k^{2})}{k}\longrightarrow S(x;s),\qquad \tfrac12\operatorname{tr}\mathsf{T}(x/k^{2})^{m} = \cosh(m\beta)\longrightarrow C(x;s).
$$
\end{lem}

\begin{proof}
Let
$$
\sigma(w) = \frac{\sinh\sqrt{w}}{\sqrt{w}} = \sum_{n\geq 0}\frac{w^{n}}{(2n+1)!},\qquad
c(w) = \cosh\sqrt{w} = \sum_{n\geq 0}\frac{w^{n}}{(2n)!},
$$
which are entire functions of $w$, with $\sigma(0) = c(0) = 1$, and note that $S(x;s) = s\,\sigma(s^{2}x)$ and $C(x;s) = c(s^{2}x)$.  All the
quantities below are even in $\beta$, so they depend only on $\beta^{2}$ and the choice of sign of $\beta$ is
irrelevant; we take the branch with $\beta\rightarrow 0$ as $u\rightarrow 0$.

Fix a compact set $K\subset\mathbb{C}$ and an interval $[0,s_{0}]$, and let $x\in K$ and $s\in[0,s_{0}]$.  By the definition of $\beta$,
$$
\cosh\beta-1 = \frac{u}{2} = \frac{x}{2k^{2}},\qquad\text{that is},\qquad \frac{\beta^{2}}{2}\,c_{1}(\beta^{2}) = \frac{x}{2k^{2}},
$$
where $c_{1}(w) = 2\big(c(w)-1\big)/w = \sum_{n\geq 0}2w^{n}/(2n+2)!$ is entire with $c_{1}(0) = 1$.  Hence $\beta^{2} = O(1/k^{2})$ uniformly for
$x\in K$, so $\beta^{2}\rightarrow 0$, and
$$
(k\beta)^{2} = \frac{x}{c_{1}(\beta^{2})} = x\big(1+O(1/k^{2})\big)\longrightarrow x
$$
uniformly for $x\in K$.  Since $m/k\rightarrow s$ uniformly, it follows that $(m\beta)^{2} = (m/k)^{2}(k\beta)^{2}\rightarrow s^{2}x$ uniformly for
$x\in K$ and $s\in[0,s_{0}]$, and these arguments remain in a fixed compact subset of $\mathbb{C}$.  As $c$ is uniformly continuous on
compact sets,
$$
\tfrac12\operatorname{tr}\mathsf{T}(x/k^{2})^{m} = \cosh(m\beta) = c\big((m\beta)^{2}\big)\longrightarrow c(s^{2}x) = C(x;s)
$$
uniformly.  Similarly, writing the quotient in terms of $\sigma$,
$$
\frac{D_{m-1}}{k} = \frac{\sinh(m\beta)}{k\sinh\beta} = \frac{m}{k}\cdot\frac{\sigma\big((m\beta)^{2}\big)}{\sigma(\beta^{2})},
$$
and the right-hand side converges uniformly to $s\,\sigma(s^{2}x) = S(x;s)$, because $m/k\rightarrow s$, $(m\beta)^{2}\rightarrow s^{2}x$, $\beta^{2}\rightarrow 0$ and
$\sigma(0) = 1$, and because $\sigma$ is uniformly continuous on compact sets and bounded away from zero near $0$.
\end{proof}

Thus the solution of \eqref{eq:recurrence} with $f_{0} = 0$ and initial slope $k(f_{1}-f_{0}) = 1$ is
$f_{m} = D_{m-1}/k$, and at the vertex $m\approx ks$ it converges to $S(x;s)$, while the solution with
$f_{0} = 1$ and symmetric initial data converges to $C(x;s)$.  This is the sense in which $\sinh$ and
$\cosh$ are the continuum limits of the transfer matrices of the subdivided edges, as described in
Remark~\ref{rem:edge-functions}, rather than functions introduced by hand.  One further limit is used
later, the ratio $k(1-D_{n-1}/D_{n})$ converging to $\sqrt{x}\coth(\sqrt{x}\,\ell_{j})$, a diagonal entry
of the Dirichlet-to-Neumann matrix (Remark~\ref{rem:discrete-DtN}).

\subsection{The completion of the proof of convergence}\label{subsec:not-viewpoint}

The convergence of the characteristic polynomials of the graphs $G_{k}$ to $\Phi$ rests on two ingredients, namely
the exact combinatorial identity of Proposition~\ref{prop:exact-forest}, valid for every $k$, and the local
convergence along each edge of Lemma~\ref{lem:continuum-limit}.  The only analytic input is that the transfer
matrices of the paths converge to those of $-g''+xg = 0$, which requires only $n_{j,k}/k\rightarrow\ell_{j}$.  In
this way the roughly $kL$ vertices of $G_{k}$ are packaged into a finite formula over the edges and the vertices
of $G_{0}$.

We will now conclude with the proof of part (1) of Theorem~\ref{mthm:convergence}.

\begin{thm}\label{thm:Phi-limit}
As $k\rightarrow\infty$,
$$
k^{V_{0}-E_{0}}\,\det\big(\Delta_{G_{k}} + x/k^{2}\big)\longrightarrow\Phi(x)
$$
uniformly for $x$ in compact subsets of $\mathbb{C}$. Consequently, $\Phi$ is entire and for every $m\geq 1$,
$$
P_{m} := \lim_{k\rightarrow\infty}\frac{c_{m}(G_{k})}{k^{E_{0}-V_{0}+2m}} = [x^{m}]\,\Phi(x),
$$
the coefficient of $x^{m}$ in the Taylor expansion of $\Phi$ at $x=0$.
\end{thm}

\begin{proof}
Put $u = x/k^{2}$.  Since $\cosh\beta = 1+u/2$, the function $\beta^{2}$ is analytic in $u$ near
$0$ with $\beta^{2} = u + O(u^{2})$, so $(k\beta)^{2} = x\,(1+O(x/k^{2}))$ uniformly for $x$ in
compact sets.  Write $n+1 = k\ell_{j,k}$ with $\ell_{j,k}\rightarrow\ell_{j}$.  The functions
$$
\frac{1}{k}D_{n}(u) = \frac{\sinh(\ell_{j,k}\,k\beta)}{k\sinh\beta},\qquad
k\,\eta_{n}(u) = 2k\sinh(\beta/2)\frac{\sinh\big((\ell_{j,k}-1/k)\,k\beta/2\big)}{\cosh(\ell_{j,k}\,k\beta/2)}
$$
are even in $k\beta$.  The first converges to $S_{j}(x)$ uniformly on compact subsets of
$\mathbb{C}$.  The second converges to $\sqrt{x}\tanh(\sqrt{x}\ell_{j}/2)$ uniformly on compact
subsets of $\mathbb{C}\setminus Z$, since the points at which some $\cosh(\sqrt{x}\ell_{j}/2)$
vanishes, namely $x = -\big((2i+1)\pi/\ell_{j}\big)^{2}$ with $i\geq 0$ and $1\leq j\leq E_{0}$, are
the points of $Z$ with odd $n$.  Individual terms of the forest sum may have poles on $Z$,
even though the sum does not; we deal with this at the end of the proof.

A spanning forest $F$ of $G_{0}$ with $r$ components has $V_{0}-r$ edges, so its complement has
$E_{0}-V_{0}+r$ edges.  Hence each term of Proposition~\ref{prop:exact-forest}, multiplied by
$k^{V_{0}-E_{0}}$, equals
$$
\prod_{\mathcal{E}_{j}\notin F}\frac{D_{n_{j,k}}(u)}{k}\ \prod_{C\in F}
\Big(\frac{x\,|C|}{k}+\sum_{h\text{ \rm at }C}k\,\eta_{n_{e(h),k}}(u)\Big)
\longrightarrow\prod_{\mathcal{E}_{j}\notin F}S_{j}(x)\prod_{C\in F}\Theta_{C}(x).
$$
Summing over the finitely many forests $F$, we conclude that
$f_{k}(x) := k^{V_{0}-E_{0}}\det(\Delta_{G_{k}}+x/k^{2})$ converges to $\Phi(x)$ uniformly on compact
subsets of $\mathbb{C}\setminus Z$.

It remains to remove the exceptional set $Z$.  Given $R>0$, choose $R'>R$ such that the circle
$|x| = R'$ does not meet $Z$; this is possible since $Z$ is discrete.  Each $f_{k}$ is a polynomial,
and $(f_{k})$ converges uniformly on $|x| = R'$.  By the maximum principle, $(f_{k})$ is uniformly
Cauchy on the disc $|x|\leq R'$, so it converges uniformly there to a holomorphic function $g$.
Since $g=\Phi$ on the disc minus $Z$, the singularities of $\Phi$ at points of $Z$ are removable.
In particular, $\Phi$ extends to an entire function, and $f_{k}\rightarrow\Phi$ uniformly on
$|x|\leq R$.  This proves the first assertion.

For the second, write $\det(\Delta_{G_{k}}+u) = \sum_{m\geq 1}c_{m}(G_{k})u^{m}$, so that
$$
k^{V_{0}-E_{0}}\det\big(\Delta_{G_{k}}+x/k^{2}\big)
= \sum_{m\geq 1}\frac{c_{m}(G_{k})}{k^{E_{0}-V_{0}+2m}}\,x^{m}.
$$
These are polynomials in $x$ converging uniformly on compact sets to the entire function $\Phi$,
so by Cauchy's estimates their coefficients converge to the Taylor coefficients of $\Phi$.
\end{proof}

\begin{lem}\label{lem:subdivision}
Let $\Gamma'$ be obtained from $\Gamma$ by inserting a vertex $w$ in the interior of an edge
$\mathcal{E}_{j}$, splitting it into two edges of lengths $a$ and $b$ with $a+b = \ell_{j}$.  Then
$\Phi_{\Gamma'} = \Phi_{\Gamma}$.
\end{lem}

\begin{proof}
Both functions are entire by Theorem~\ref{thm:Phi-limit}, so it suffices to prove the identity for
real $x>0$.  For such $x$ all the functions $S$ are positive, so $x\notin Z(\Gamma)\cup Z(\Gamma')$.
The diagonal entry of $M_{\Gamma'}(x)$ at $w$ is
$$
\frac{C(x;a)}{S(x;a)}+\frac{C(x;b)}{S(x;b)} = \frac{S_{j}(x)}{S(x;a)S(x;b)} > 0,
$$
by the addition formula $S(x;a+b) = S(x;a)C(x;b) + C(x;a)S(x;b)$.  Taking the Schur complement at
$w$,
$$
\det M_{\Gamma'}(x) = \frac{S_{j}(x)}{S(x;a)S(x;b)}\cdot\det\widetilde{M}(x),
$$
where $\widetilde{M}(x)$ is the Schur complement of the $(w,w)$ entry.  We claim
$\widetilde{M}(x) = M_{\Gamma}(x)$.  For $\varphi\in\mathbb{C}^{V(G_{0})}$, the Schur complement is
computed by extending $\varphi$ to a vector $\tilde\varphi$ on $V(\Gamma')$ whose value at $w$ is chosen so
that $(M_{\Gamma'}(x)\tilde\varphi)_{w} = 0$; then $\widetilde{M}(x)\varphi$ is the restriction of
$M_{\Gamma'}(x)\tilde\varphi$ to $V(G_{0})$.  The condition at $w$ says that $f_{\tilde\varphi}$ is $C^{1}$
across $w$, so it solves $g''=xg$ on all of $\mathcal{E}_{j}$ with the endpoint values of $\varphi$.  By
uniqueness, $f_{\tilde\varphi} = f_{\varphi}$ on $\Gamma$, so the outgoing derivatives at the vertices of
$G_{0}$ agree, which proves the claim.  Finally, $\prod S$ over the edges of $\Gamma'$ equals
$S(x;a)S(x;b)/S_{j}(x)$ times the corresponding product for $\Gamma$, and the factors cancel.
\end{proof}

\begin{remark}\label{rem:discrete-DtN}
The Schur complement $S_{k}(u)$ in the proof of Proposition~\ref{prop:exact-forest} is a discrete
Dirichlet-to-Neumann matrix.  Given values at the vertices of $G_{0}$, it returns the net flux at these vertices
of the unique extension that satisfies $(\Delta_{G_{k}}+u)f = 0$ at all subdivision vertices.  It converges to
the metric Dirichlet-to-Neumann matrix
$$
k\,S_{k}(x/k^{2})\longrightarrow M(x)
$$
uniformly for $x$ in compact subsets of $\mathbb{C}\setminus Z$.  Indeed, with $u = x/k^{2}$ and $n = n_{j,k}$, we have
$k/D_{n}(u)\rightarrow 1/S_{j}(x) = \sqrt{x}/\sinh(\sqrt{x}\,\ell_{j})$ and $ku\rightarrow 0$.  Moreover, since
$D_{n}-D_{n-1} = \cosh((2n+1)\beta/2)/\cosh(\beta/2)$,
$$
k\Big(1-\frac{D_{n-1}}{D_{n}}\Big) = \frac{\cosh\big((\ell_{j,k}-\frac{1}{2k})\,k\beta\big)}{\cosh(\beta/2)\,D_{n}/k}
\longrightarrow\frac{\cosh(\sqrt{x}\,\ell_{j})}{S_{j}(x)} = \sqrt{x}\coth(\sqrt{x}\,\ell_{j}),
$$
and the loop terms $-2k/D_{n}$ converge to $-2\sqrt{x}/\sinh(\sqrt{x}\,\ell_{j})$.  Comparing with the formula for
$S_{k}(u)$ gives the claim.  We have also confirmed this convergence numerically, at rate $O(1/k)$, for a graph
with a loop and a double edge.
\end{remark}

\begin{remark}\label{rem:harrison-weyand}
Relations between discrete and metric graphs have been used before to study spectral zeta functions.
In \cite{HW18}, the spectral zeta function of an \emph{equilateral}
metric graph is expressed in terms of the spectral zeta function of the normalized Laplacian on the same
combinatorial graph.  The method there rests on von Below's correspondence between the two spectra
\cite{vB85}, which holds only when all edges have the same length, and therefore applies
only to equilateral graphs.  Theorem~\ref{thm:Phi-limit} is
considerably more general.  It holds for arbitrary positive edge lengths, including rationally
independent ones, and for graphs with loops and multiple edges.  It is based on an exact identity,
Proposition~\ref{prop:exact-forest}, valid for every $k$.  It identifies the limit of the
characteristic polynomials $\det(\Delta_{G_{k}}+x/k^{2})$ as an explicit entire function, not only
the limit of a single spectral invariant.  As a result, it yields simultaneously the regularized
determinant, all the forest invariants $P_{m}$, and all the values $\zeta_{\Gamma}(r)$ at positive
integers.  For equilateral graphs our results imply those of \cite{HW18}.

\end{remark}

\begin{corollary}\label{cor:kirchhoff-type}
For a component $C$ of a spanning forest of $G_{0}$, let
$\ell^{*}(C) = \frac12\sum_{h\text{ \rm at }C}\ell_{e(h)}$, the total length of the edges of
$G_{0}$ inside $C$ plus half the length of each edge with exactly one end in $C$.  Then
$$
P_{m} = \sum_{F:\ m\text{ \rm components}}\ \prod_{\mathcal{E}_{j}\notin F}\ell_{j}\
\prod_{C\in F}\ell^{*}(C)\ +\ R_{m},
$$
where $R_{m}$ collects the terms from forests with fewer than $m$ components,
obtained from \eqref{eq:Phi} using the series expansions of $\sinh$ and $\tanh$.
In particular, $P_{1} = L\sum_{T}\prod_{\mathcal{E}_{j}\notin T}\ell_{j}$,
which is in agreement with Proposition~\ref{prop:logdet-asymptotics}.
\end{corollary}

\begin{proof}
Since $S_{j}(x) = \ell_{j}+O(x)$ and $\Theta_{C}(x) = x\,\ell^{*}(C)+O(x^{2})$, a forest with $r$
components contributes to $\Phi$ only at orders $x^{r}$ and higher, and its contribution at
order $x^{r}$ is
$$
\prod_{\mathcal{E}_{j}\notin F}\ell_{j}\prod_{C}\ell^{*}(C).
$$
\end{proof}

\begin{remark}
The formula of Corollary~\ref{cor:kirchhoff-type} is the metric analogue of the matrix-forest
theorem \eqref{eq:matrix-forest}.  The number of vertices $V(\mathcal{F}_{i})$ of a component is
replaced by its metric size $\ell^{*}(C)$, and each deleted edge is weighted by its length.  Unlike
the discrete case, forests with fewer than $m$ components also contribute to $P_{m}$.  Indeed, such forests
account for the $m$-forests of $G_{k}$ in which some edge of $G_{0}$ is cut more than once.  For example, the
circle ($V_{0}=E_{0}=1$) has no $2$-forests, yet $P_{2} = \ell^{4}/12$.
\end{remark}

\section{Zeros of \texorpdfstring{$\Phi$}{Phi} and the spectrum of the metric graph}\label{sec:zeros}

In this section we show that the zeros of the metric forest function $\Phi$ are exactly the
negatives of the eigenvalues of the Laplacian on the metric graph, counted with multiplicity.
As consequences we obtain a Hadamard product for $\Phi$, formulas for the spectral zeta function
at all positive integers, hence the convergence of the rescaled spectrum of $G_{k}$ to that of the
metric graph $\Gamma = (G_{0},\Lambda)$.

\subsection{The functions \texorpdfstring{$S_{j}$}{Sj}}

\begin{remark}\label{rem:motivation-S}
Besides being the renormalized limit $D_{n_{j,k}}(x/k^{2})/k\rightarrow S_{j}(x)$ of the determinant of the
subdivided edge, the function $S_{j}(x)$ is essentially the zeta-regularized determinant of the Dirichlet
Laplacian of that edge.  The
Dirichlet eigenvalues of $-d^{2}/ds^{2}$ on $[0,\ell_{j}]$ are $(\pi n/\ell_{j})^{2}$, $n\geq 1$, and
$$
S_{j}(x) = \ell_{j}\prod_{n\geq 1}\Big(1+\frac{x\,\ell_{j}^{2}}{\pi^{2}n^{2}}\Big),\qquad
\text{\rm det}\Big(-\frac{d^{2}}{ds^{2}}+x\Big)_{\text{\rm Dirichlet on }[0,\ell_{j}]} = 2\,S_{j}(x).
$$
The second identity is classical; see \cite{BFK95}.  Hence $\prod_{j}S_{j}(x) = 2^{-E_{0}}\text{\rm det}(\Delta_{D}+x)$, where $\Delta_{D}$ is the Laplacian on $\Gamma$ with
Dirichlet conditions at all vertices.  This is the source of the factor $2^{E_{0}}$ in Theorem~\ref{thm:det}, and
the reason why the formula $\Phi = \prod_{j}S_{j}\det M$ matches the formula of \cite{Fr06}
$\det(\Delta_{\Gamma}+x) = \prod_{v}d_{v}^{-1}\det(\Delta_{D}+x)\det M(x)$.
\end{remark}

\subsection{Zeros of \texorpdfstring{$\Phi$}{Phi}}

\begin{thm}\label{thm:zeros}
For $x_{0}\in\mathbb{C}$, $\Phi(x_{0}) = 0$ if and only if $-x_{0}$ is an eigenvalue of
$\Delta_{\Gamma}$, and the order of the zero is the multiplicity of the eigenvalue.
\end{thm}

\begin{proof}
First suppose $x_{0}\notin Z$.  Then $\prod_{j}S_{j}$ does not vanish at $x_{0}$, so the order of
vanishing of $\Phi$ at $x_{0}$ equals that of $\det M(x)$.  If $-x_{0}$ is not an eigenvalue, then
$\det M(x_{0})\neq 0$ by Lemma~\ref{lem:DtN-kernel}.  If $-x_{0}$ is an eigenvalue, then $x_{0}$ is
real because $\Delta_{\Gamma}$ is self-adjoint, and $M(x)$ is holomorphic near $x_{0}$ and real
symmetric for real $x$.  For nonzero $\varphi\in\ker M(x_{0})$, the function $f_{\varphi}(x_{0})$ is nonzero
because its vertex values are $\varphi$, so $\varphi^{T}M'(x_{0})\varphi>0$ by Lemma~\ref{lem:energynew}.  By
Lemma~\ref{lem:order}, $\det M$ vanishes at $x_{0}$ to order $\dim\ker M(x_{0})$, which is the
multiplicity of $-x_{0}$ by Lemma~\ref{lem:DtN-kernel}.

Now suppose $x_{0}\in Z$, so $x_{0}<0$.  For each edge $\mathcal{E}_{j}$ with $S_{j}(x_{0})=0$, insert
a vertex splitting $\mathcal{E}_{j}$ into edges of lengths $a_{j}$ and $b_{j}$, where $a_{j}$ is
chosen so that neither $\sqrt{|x_{0}|}\,a_{j}$ nor $\sqrt{|x_{0}|}\,b_{j}$ lies in $\pi\mathbb{Z}$.
This excludes only countably many values of $a_{j}$.  The resulting metric graph $\Gamma'$ has
the same Laplacian as $\Gamma$, and $x_{0}\notin Z(\Gamma')$, since the edges not subdivided
already satisfy $S_{j}(x_{0})\neq 0$.  By Lemma~\ref{lem:subdivision}, $\Phi_{\Gamma'} = \Phi_{\Gamma}$,
and the first case applied to $\Gamma'$ completes the proof.
\end{proof}

\begin{lem}\label{lem:growth}
There is a constant $A>0$ such that $|\Phi(x)|\leq A\,(1+|x|)^{V_{0}/2}\,e^{L\sqrt{|x|}}$ for all
$x\in\mathbb{C}$.  In particular, $\Phi$ is an entire function of order at most $1/2$.
\end{lem}

\begin{proof}
We use the forest formula \eqref{eq:Phi}, choosing the branch of $\sqrt{x}$ with
$\text{\rm Re}\sqrt{x}\geq 0$.  From the power series, $|S(x;s)|\leq s\,e^{s\sqrt{|x|}}$.  For
the factors $\Theta_{C}$, we use the following elementary estimate.  For each $\delta>0$ there is
$c_{\delta}>0$ such that $|\cosh w|\geq c_{\delta}e^{|\text{\rm Re}\,w|}$ whenever $w$ has distance
at least $\delta$ from the zeros $i\pi(\mathbb{Z}+\frac12)$ of $\cosh$.  Indeed, for
$|\text{\rm Re}\,w|\geq 1$ this follows from $|\cosh w|\geq\sinh|\text{\rm Re}\,w|$, and for
$|\text{\rm Re}\,w|\leq 1$ it follows from periodicity and compactness.  Since
$|\sinh w|\leq e^{|\text{\rm Re}\,w|}$, we get $|\tanh w|\leq 1/c_{\delta}$ for such $w$.

Call $r>0$ good if $|\ell_{j}r/2 - \pi(i+\frac12)|\geq\delta$ for all $j$ and all integers $i$.  If
$|x| = r^{2}$ with $r$ good, then each $w = \sqrt{x}\,\ell_{j}/2$ has $|w| = \ell_{j}r/2$, so its distance
from every zero of $\cosh$ is at least $\delta$.  Indeed the zeros lie on the imaginary axis, so
$|w-i\pi(k+\frac12)|\geq\big||w|-\pi|k+\tfrac12|\big|$, and the condition defining a good $r$ runs over all
integers $i$, hence over every value of $\pi|k+\frac12|$.  Hence
$|\Theta_{C}(x)|\leq 2E_{0}\,r/c_{\delta}$, and every term of \eqref{eq:Phi} is bounded by a
constant times $(1+r)^{V_{0}}e^{Lr}$, since a forest has at most $V_{0}$ components.  For each $j$,
the set of $r$ which are not good for $j$ is a union of intervals of length $4\delta/\ell_{j}$ with
centres spaced $2\pi/\ell_{j}$ apart.  A unit interval meets at most $\ell_{j}/(2\pi)+2$ of these intervals, so the
set of $r\in[r_{0},r_{0}+1]$ that are not good has total length at most $\sum_{j}(\ell_{j}/2\pi+2)\,4\delta/\ell_{j}$.  For
$\delta$ small enough, depending only on $\Lambda$, this is less than $1$, so every interval $[r_{0},r_{0}+1]$
contains a good $r$.  Given $x$, choose a good
$r\in[\sqrt{|x|},\sqrt{|x|}+1]$.  By the maximum principle on the disc of radius $r^{2}$,
$|\Phi(x)|$ is at most a constant times $(1+r)^{V_{0}}e^{Lr}$, which gives the stated bound.
\end{proof}

\begin{thm}\label{thm:hadamard}
The series $\sum_{j\geq 1}\mu_{j}^{-s}$ converges for $s>1/2$, and
$$
\Phi(x) = P_{1}\,x\prod_{j\geq 1}\Big(1+\frac{x}{\mu_{j}}\Big),\qquad
P_{1} = L\sum_{T}\prod_{\mathcal{E}_{j}\notin T}\ell_{j},
$$
where the product converges absolutely and uniformly on compact subsets of $\mathbb{C}$.
\end{thm}

\begin{proof}
By Lemma~\ref{lem:growth}, $\Phi$ has order at most $1/2$.  For an entire function of finite order
$\rho$, the series $\sum|z|^{-s}$ over its nonzero zeros converges for every $s>\rho$; see, e.g.,
\cite{Ti39}.  By Theorem~\ref{thm:zeros} the nonzero zeros are the $-\mu_{j}$, repeated according
to multiplicity, and $x=0$ is a simple zero because $\mu_{0}=0$ is simple.  Since the order is
less than one, Hadamard's factorization theorem gives
$\Phi(x) = c\,x\prod_{j\geq 1}(1+x/\mu_{j})$ for a constant $c$.  Comparing coefficients of $x$ and
using Theorem~\ref{thm:Phi-limit} and Corollary~\ref{cor:kirchhoff-type}, $c = \Phi'(0) = P_{1}$.
\end{proof}

\subsection{Consequences}

By Theorem~\ref{thm:hadamard}, the spectral zeta function $\zeta_{\Gamma}$ converges for
$\text{\rm Re}(s)>1/2$.

\begin{corollary}\label{cor:zeta-integers}
For $|x|<\mu_{1}$,
$$
\log\Big(\frac{\Phi(x)}{P_{1}\,x}\Big) = \sum_{r\geq 1}\frac{(-1)^{r-1}\zeta_{\Gamma}(r)}{r}\,x^{r}.
$$
Consequently, for each integer $r\geq 1$, $\zeta_{\Gamma}(r)$ is a polynomial with rational
coefficients in the ratios $P_{2}/P_{1},\dots,P_{r+1}/P_{1}$, given by Newton's identities; see
Theorem~\ref{thm:newton} for the explicit formulas.  For
example,
$$
\zeta_{\Gamma}(1) = \frac{P_{2}}{P_{1}},\qquad
\zeta_{\Gamma}(2) = \Big(\frac{P_{2}}{P_{1}}\Big)^{2} - 2\,\frac{P_{3}}{P_{1}}.
$$
\end{corollary}

\begin{proof}
Take logarithms in Theorem~\ref{thm:hadamard} and expand
$\log(1+x/\mu_{j}) = \sum_{r\geq 1}(-1)^{r-1}x^{r}/(r\mu_{j}^{r})$.  The double series converges
absolutely for $|x|<\mu_{1}$ since $\sum_{j}\mu_{j}^{-1}<\infty$.  On the other side,
$\Phi(x)/(P_{1}x) = 1+\sum_{m\geq 1}(P_{m+1}/P_{1})x^{m}$ by Theorem~\ref{thm:Phi-limit}, and
comparing coefficients gives the formulas.
\end{proof}

\begin{corollary}\label{cor:spectral-convergence}
Let $0 = \lambda_{0,k}<\lambda_{1,k}\leq\cdots$ be the eigenvalues of $\Delta_{G_{k}}$.  For each
$j\geq 0$, $k^{2}\lambda_{j,k}\rightarrow\mu_{j}$ as $k\rightarrow\infty$.
\end{corollary}

\begin{proof}
The polynomial $f_{k}(x) = k^{V_{0}-E_{0}}\det(\Delta_{G_{k}}+x/k^{2})$ vanishes exactly at the points
$-k^{2}\lambda_{j,k}$, with multiplicity, and $f_{k}\rightarrow\Phi$ uniformly on compact sets by
Theorem~\ref{thm:Phi-limit}.  Let $R>0$ be distinct from all $\mu_{j}$.  By Hurwitz's theorem
applied to the disc $|x|<R$ and to small discs around each $-\mu_{j}$ with $\mu_{j}<R$, for
large $k$ the number of $j$ with $k^{2}\lambda_{j,k}<R$ equals the number of $j$ with $\mu_{j}<R$,
and these rescaled eigenvalues converge to the corresponding $\mu_{j}$, with multiplicity.  Since
all zeros are real, the convergence respects the ordering.
\end{proof}

\begin{remark}
Corollary~\ref{cor:spectral-convergence} recovers a theorem of \cite{Fa06}, which gives the
convergence of the rescaled eigenvalues of discrete models of a metrized graph, together with
convergence of the interpolated eigenvectors.  Our proof is different.  Here the result is a
consequence of the exact forest expansion, without separate estimates on eigenfunctions.
\end{remark}

\section{Degeneration of the discrete theta functions}\label{sec:heat-degeneration}

In this section we prove part (5) of Theorem~\ref{mthm:convergence}, the convergence
$k^{-2s}\zeta_{G_{k}}(s)\rightarrow\zeta_{\Gamma}(s)$ for $\text{\rm Re}(s)>1/2$.  At the positive integers
$s = r$ this also follows by combinatorics, since $k^{-2r}\zeta_{G_{k}}(r)$ is a polynomial in the normalized
forest counts of $G_{k}$, which converge by part (1), and that route is taken in
Corollary~\ref{cor:forest-hierarchy}.  Here we use theta functions to interpolate between the integer values.  This
requires uniform control of all the eigenvalues at once, and it comes from the convergence of the individual
eigenvalues (Corollary~\ref{cor:spectral-convergence}) together with a lower bound on the eigenvalues that is
uniform in $k$ (Lemma~\ref{lem:uniform-lower}).  The theta functions used here are sums over eigenvalues, and no heat
kernel is constructed.

The method stops at $\text{\rm Re}(s) = 1/2$, the uniform bound $\theta_{G_{k}}(k^{2}t)\leq A(1+t^{-1/2})$ below
reflecting the one-dimensional bulk.  The special values at $s = 0$ are nevertheless obtained in
Section~\ref{sec:det}, namely $\zeta_{\Gamma}(0)$ from the expansion of the theta function of $\Gamma$
(Corollary~\ref{cor:zeta-continuation}), and $\zeta_{\Gamma}'(0)$ from the determinants of the $G_{k}$, after the
divergent terms are removed explicitly (Corollary~\ref{cor:main}).  The behavior of $k^{-2s}\zeta_{G_{k}}(s)$ for
$\text{\rm Re}(s)<1/2$ is discussed, but not proved, after Theorem~\ref{mthm:convergence}.

We now show that the rescaled theta functions of the discrete graphs $G_{k}$ converge to the theta
function of the metric graph $\Gamma$, with bounds uniform in $k$.  These are the analogues for metric graphs
of Proposition~5.2 and Lemma~5.3 of \cite{CJK10}.  Let
$$
\theta_{G_{k}}(t) = \sum_{j=0}^{V_{k}-1}e^{-\lambda_{j,k}t}
$$
denote the theta function of $\Delta_{G_{k}}$, where
$0=\lambda_{0,k}<\lambda_{1,k}\leq\cdots\leq\lambda_{V_{k}-1,k}$ are its eigenvalues.

\begin{lem}\label{lem:uniform-lower}
For all $k$ sufficiently large and all $j\geq V_{0}+E_{0}$,
$$
k^{2}\lambda_{j,k}\geq\frac{(j-V_{0}-E_{0})^{2}}{L^{2}}.
$$
\end{lem}

\begin{proof}
Let $W$ be the graph obtained from $G_{k}$ by deleting every edge incident to a vertex of $G_{0}$.
It is the disjoint union of the $V_{0}$ vertices of $G_{0}$, now isolated, and $E_{0}$ paths, the
$j$-th having $n_{j,k}$ vertices.  Since $\Delta_{G_{k}}-\Delta_{W}$ is the Laplacian of the deleted
edges, it is positive semi-definite.  By the min-max principle, see for example \cite[Theorem~4.5.3]{HJ13},
$\lambda_{j,k}\geq\lambda_{j}(W)$ for all $j$, where the eigenvalues of $W$ are also listed in
increasing order.  The eigenvalues of a
path with $n$ vertices are $4\sin^{2}(\pi i/2n)$ for $0\leq i\leq n-1$.  For $0<\tau<4$ and
$y = \sqrt{\tau}/2$, the number of these that are less than $\tau$ is the number of $i$ with
$i<(2n/\pi)\arcsin y$, which is at most $1+ny$ because $\arcsin y\leq\pi y/2$.  Hence the number
$N_{W}(\tau)$ of eigenvalues of $W$ less than $\tau$ satisfies
$$
N_{W}(\tau)\leq V_{0}+E_{0}+\frac{H_{k}\sqrt{\tau}}{2}.
$$
For $j\geq V_{0}+E_{0}$ and $\sqrt{\tau} = 2(j-V_{0}-E_{0})/H_{k}$, this gives $N_{W}(\tau)\leq j$, and
therefore $\lambda_{j}(W)\geq\tau$.  The case $\tau\geq 4$ is empty, since then
$j\geq V_{0}+E_{0}+H_{k}$, which exceeds $V_{k}-1$.  Since
$H_{k}\leq kL_{k}\leq kL(1+o(1))$, we obtain
$k^{2}\lambda_{j,k}\geq 4(j-V_{0}-E_{0})^{2}/(L_{k}^{2}(1+o(1)))$, which implies the stated bound for
$k$ large.
\end{proof}

\begin{prop}\label{prop:heat-degeneration}
For every $t>0$,
$$
\lim_{k\rightarrow\infty}\theta_{G_{k}}(k^{2}t) = Q(t) = \sum_{j\geq 0}e^{-\mu_{j}t}.
$$
Moreover, there are constants $A$ and $a>0$, independent of $k$ and $t$, such that for all $k$
sufficiently large,
$$
\theta_{G_{k}}(k^{2}t)\leq A\big(1+t^{-1/2}\big)\quad(t>0),\qquad
0\leq\theta_{G_{k}}(k^{2}t)-1\leq A\,e^{-at}\quad(t\geq 1).
$$
\end{prop}

\begin{proof}
By Lemma~\ref{lem:uniform-lower}, for $k$ large,
$$
\theta_{G_{k}}(k^{2}t)\leq V_{0}+E_{0}+\sum_{m\geq 0}e^{-m^{2}t/L^{2}}
\leq V_{0}+E_{0}+1+\frac{L}{2}\sqrt{\frac{\pi}{t}},
$$
which is the first bound.  The same estimate shows that for fixed $t>0$, the terms of
$\sum_{j}e^{-k^{2}\lambda_{j,k}t}$ are dominated by a summable sequence independent of $k$, namely
$1$ for $j<V_{0}+E_{0}$ and $e^{-(j-V_{0}-E_{0})^{2}t/L^{2}}$ otherwise.  Since
$k^{2}\lambda_{j,k}\rightarrow\mu_{j}$ for each $j$ by Corollary~\ref{cor:spectral-convergence},
dominated convergence gives $\theta_{G_{k}}(k^{2}t)\rightarrow Q(t)$.

For the second bound, Corollary~\ref{cor:spectral-convergence} gives
$k^{2}\lambda_{1,k}\geq\mu_{1}/2 =: a$ for $k$ large.  For $t\geq 1$ we then have
$$
\theta_{G_{k}}(k^{2}t)-1 = \sum_{j\geq 1}e^{-k^{2}\lambda_{j,k}t}
\leq e^{-a(t-1)}\sum_{j\geq 1}e^{-k^{2}\lambda_{j,k}}\leq e^{a}\,\theta_{G_{k}}(k^{2})\,e^{-at},
$$
and $\theta_{G_{k}}(k^{2})$ is bounded by the first estimate.
\end{proof}

\begin{remark}
In \cite{CJK10}, the pointwise convergence of theta functions was proved from the expression of the
heat kernel on a discrete torus as a periodized sum of I-Bessel functions, and the long-time bound
from explicit knowledge of the spectrum.  Here both come from the spectral convergence of
Corollary~\ref{cor:spectral-convergence}, which in turn comes from the exact forest expansion.  In
place of the explicit spectrum, we use the monotonicity argument of Lemma~\ref{lem:uniform-lower}.
\end{remark}

\begin{corollary}\label{cor:zeta-degeneration}
Let $\zeta_{G_{k}}(s) = \sum_{j\geq 1}\lambda_{j,k}^{-s}$.  For every $s\in\mathbb{C}$ with
$\text{\rm Re}(s)>1/2$,
$$
\lim_{k\rightarrow\infty}k^{-2s}\zeta_{G_{k}}(s) = \zeta_{\Gamma}(s).
$$
\end{corollary}

\begin{proof}
We have $k^{-2s}\zeta_{G_{k}}(s) = \sum_{j\geq 1}(k^{2}\lambda_{j,k})^{-s}$, and each term converges to
$\mu_{j}^{-s}$ by Corollary~\ref{cor:spectral-convergence}.  For $k$ large, $k^{2}\lambda_{j,k}\geq\mu_{1}/2$
for all $j\geq 1$, and $k^{2}\lambda_{j,k}\geq (j-V_{0}-E_{0})^{2}/L^{2}$ for $j\geq V_{0}+E_{0}$ by
Lemma~\ref{lem:uniform-lower}.  So the terms are dominated by a constant for finitely many $j$ and
by a constant times $(j-V_{0}-E_{0})^{-2\,\text{\rm Re}(s)}$ otherwise, which is summable for
$\text{\rm Re}(s)>1/2$.  Dominated convergence completes the proof.
\end{proof}

\begin{remark}\label{rem:CJK}
The results of this section and of Section~\ref{sec:graphapproximations} correspond closely to those of
\cite{CJK10} for discrete tori degenerating to real tori.  Proposition~\ref{prop:heat-degeneration}, the convergence of the
rescaled theta functions, corresponds to the convergence of the theta functions of discrete tori.
Corollary~\ref{cor:zeta-degeneration} is the analogue of \cite[(8)]{CJK10}, with the threshold $d/2$ for $d$-dimensional
tori replaced by $1/2$.  Proposition~\ref{prop:logdet-asymptotics} and Corollary~\ref{cor:main} correspond to the
asymptotics of the number of spanning trees of discrete tori, whose constant term is the regularized determinant of the
limiting torus.  There are two differences.  First, a metric graph is one-dimensional but singular at its vertices, and
the vertices produce local contributions, namely the constant $\chi/2$ in the theta expansion, and the factor $\prod_{v}d_{v}/2^{E_{0}}$ in
the determinant.  For tori there are no such terms.  Second, the convergence of the characteristic polynomials in
Theorem~\ref{mthm:convergence} gives more than the convergence of theta functions and zeta functions.  It also gives the
spanning forest counts of all orders (Theorem~\ref{mthm:forests}), and through them all the special values.  In particular,
at positive integers $s = r$, Corollary~\ref{cor:zeta-degeneration} follows from the forest counts, since
$k^{-2r}\zeta_{G_{k}}(r)$ is a polynomial in the ratios $c_{m}(G_{k})/c_{1}(G_{k})$ with $m\leq r+1$ (Corollary~\ref{cor:forest-hierarchy}).
In the region $\text{\rm Re}\,s<1/2$ one expects an asymptotic expansion analogous to those of \cite{CJK10} and \cite{FK17}, with
local vertex contributions; see the discussion following Theorem~\ref{mthm:convergence}.
\end{remark}

\section{Theta function asymptotics and the zeta-regularized determinant}\label{sec:det}

In this section we derive the small-time asymptotics of the theta function of $\Delta_{\Gamma}$ from the
function $\Phi$, and then evaluate the zeta-regularized determinant of $\Delta_{\Gamma}+x$ for all
$x\geq 0$.  The methods are analytic and use neither the heat kernel nor Roth's trace formula.
As a corollary of the small-time asymptotics of the theta function,
we obtain the formula of \cite{Fr06}, together with the values of
$\zeta_{\Gamma}$ at zero and at the negative integers.  Throughout, $\chi = V_{0}-E_{0}$ denotes the
Euler characteristic of $G_{0}$, and $d_{v}$ denotes the degree of the vertex $v$, where a loop
contributes $2$.  For $x\in\mathbb{C}\setminus(-\infty,0]$ we use the principal branches of
$\sqrt{x}$ and $\log x$.

\subsection{Asymptotics of $\Phi$ in sectors}

For $0<\varepsilon<\pi$, let
$\Sigma_{\varepsilon} = \{x\in\mathbb{C}\setminus\{0\}:|\arg x|\leq\pi-\varepsilon\}$.

\begin{lem}\label{lem:Phi-sector}
Fix $\varepsilon\in(0,\pi)$, and let $c = \sin(\varepsilon/2)\min_{j}\ell_{j}$.  There is $r_{0}>0$
such that $\Phi$ has no zeros in $\{x\in\Sigma_{\varepsilon}:|x|\geq r_{0}\}$, and on this set
$$
\log\Phi(x) = L\sqrt{x} + \frac{\chi}{2}\log x + \log\Big(\frac{\prod_{v}d_{v}}{2^{E_{0}}}\Big)
+ O\big(e^{-c\sqrt{|x|}}\big),
$$
for a suitable branch of $\log\Phi$.
\end{lem}

\begin{proof}
Put $\rho = \sqrt{x}$.  For $x\in\Sigma_{\varepsilon}$ we have
$\text{\rm Re}\,\rho\geq|\rho|\sin(\varepsilon/2)$, so $|e^{-\rho\ell_{j}}|\leq e^{-c|\rho|}$ for each $j$.  Hence, for $|x|$ large,
$$
S_{j}(x) = \frac{e^{\rho\ell_{j}}}{2\rho}\big(1-e^{-2\rho\ell_{j}}\big),\qquad
\frac{C_{j}(x)}{S_{j}(x)} = \rho\big(1+O(e^{-2c|\rho|})\big),\qquad
\frac{1}{S_{j}(x)} = O\big(|\rho|e^{-c|\rho|}\big).
$$
Hence each diagonal entry of $M(x)$ is $\rho\,d_{v}(1+O(e^{-c|\rho|}))$, where a loop at $v$ contributes
$2C_{j}/S_{j}-2/S_{j} = 2\rho\tanh(\rho\ell_{j}/2) = 2\rho(1+O(e^{-c|\rho|}))$.  Each off-diagonal entry
is $O(|\rho|e^{-c|\rho|})$.  Therefore
$\det M(x) = \rho^{V_{0}}\prod_{v}d_{v}\,(1+O(e^{-c|\rho|}))$ and
$$
\Phi(x) = \prod_{j}S_{j}(x)\cdot\det M(x)
= \frac{e^{L\rho}}{2^{E_{0}}}\,\rho^{\chi}\prod_{v}d_{v}\,\big(1+O(e^{-c|\rho|})\big),
$$
from which both assertions follow.
\end{proof}

\subsection{The theta function}

Let $Q(t) = \sum_{j\geq 0}e^{-\mu_{j}t}$ for $t>0$ be the theta function of
$\Delta_{\Gamma}$.  The series converges, because $\sum_{j}\mu_{j}^{-1}<\infty$ by
Theorem~\ref{thm:hadamard} and the $\mu_{j}$ are increasing, so that $j/\mu_{j}\rightarrow 0$.

\begin{thm}\label{thm:heat-asymptotics}
There is $c'>0$ such that
$$
Q(t) = \frac{L}{\sqrt{4\pi t}} + \frac{\chi}{2} + O\big(e^{-c'/t}\big)\qquad\text{as }t\rightarrow 0^{+}.
$$
\end{thm}

\begin{proof}
For $x\notin(-\infty,0]$, let $R(x) = \sum_{j\geq 0}(x+\mu_{j})^{-1}$, the harmonic series of the sequence
$\{\mu_{j}\}$ in the sense of \cite{JoLa93}, which converges since
$\sum_{j}\mu_{j}^{-1}<\infty$.  Taking the logarithmic derivative of the Hadamard product in
Theorem~\ref{thm:hadamard}, which converges uniformly on compact sets, gives
$R(x) = \Phi'(x)/\Phi(x)$.  For $\text{\rm Re}\,x>0$, Fubini's theorem gives
$R(x) = \int_{0}^{\infty}e^{-xt}Q(t)\,dt$.

Fix $\varepsilon\in(0,\pi/4)$, so that $\Sigma_{2\varepsilon}$ contains the closed right half-plane minus $0$ and
$\cos 2\varepsilon>0$.  Applying Cauchy's estimate to the remainder in
Lemma~\ref{lem:Phi-sector} on discs of radius $\delta|x|$, with $\delta$ small, shows that for
$x\in\Sigma_{2\varepsilon}$ with $|x|$ large,
$$
R(x) = \frac{L}{2\sqrt{x}} + \frac{\chi}{2x} + O\big(e^{-c_{1}\sqrt{|x|}}\big)
$$
for some $c_{1}>0$.  Let
$$
q(t) = Q(t)-\frac{L}{\sqrt{4\pi t}}-\frac{\chi}{2},\qquad
r(x) = R(x)-\frac{L}{2\sqrt{x}}-\frac{\chi}{2x}.
$$
Since $\int_{0}^{\infty}e^{-xt}(4\pi t)^{-1/2}\,dt = (4x)^{-1/2}$, the function $r$ is the Laplace
transform of $q$ for $\text{\rm Re}\,x>0$.  Moreover, $r$ is holomorphic on
$\mathbb{C}\setminus(-\infty,0]$, and $r(x) = O(e^{-c_{1}\sqrt{|x|}})$ in $\Sigma_{2\varepsilon}$ as
$|x|\rightarrow\infty$.  In particular $r$ is absolutely integrable on each vertical line
$\text{\rm Re}\,x = \sigma>0$.  Since $q$ is continuous on $(0,\infty)$, Laplace inversion gives
$$
q(t) = \frac{1}{2\pi i}\int_{\sigma-i\infty}^{\sigma+i\infty}e^{xt}\,r(x)\,dx\qquad(t>0).
$$
Given $\rho_{0}$ large, deform the contour into the path consisting of the ray from
$\infty e^{-i(\pi-2\varepsilon)}$ to $\rho_{0}e^{-i(\pi-2\varepsilon)}$, the arc
$|x|=\rho_{0}$ through the positive real axis, and the ray from $\rho_{0}e^{i(\pi-2\varepsilon)}$ to
$\infty e^{i(\pi-2\varepsilon)}$.  The region between the vertical line and the new path lies in $\Sigma_{2\varepsilon}$ together
with a bounded subset of the half-plane $\text{\rm Re}\,x>0$.  In particular it does not meet
$(-\infty,0]$, where the singularities of $r$ lie, since the closed curve formed by the two paths
has winding number zero about every point of $(-\infty,0]$.  The contributions of the arcs at
infinity connecting the two paths tend to zero, because $|e^{xt}|\leq e^{\sigma t}$ there and $r$
decays in $\Sigma_{2\varepsilon}$.  Hence the deformation is justified.  On the rays $|e^{xt}| = e^{-|x|t\cos 2\varepsilon}$, and on the arc
$|e^{xt}|\leq e^{\rho_{0}t}$.  Hence
$$
|q(t)|\leq C\Big(\rho_{0}\,e^{\rho_{0}t-c_{1}\sqrt{\rho_{0}}}
+\int_{\rho_{0}}^{\infty}e^{-\rho t\cos 2\varepsilon-c_{1}\sqrt{\rho}}\,d\rho\Big)
\leq C\Big(\rho_{0}\,e^{\rho_{0}t-c_{1}\sqrt{\rho_{0}}}+\frac{e^{-c_{1}\sqrt{\rho_{0}}}}{t\cos 2\varepsilon}\Big).
$$
Choosing $\rho_{0} = (c_{1}/2t)^{2}$ gives $|q(t)| = O(t^{-2}e^{-c_{1}^{2}/4t})$ as $t\rightarrow 0^{+}$.
\end{proof}

\begin{remark}
The asymptotics in Theorem~\ref{thm:heat-asymptotics} are known.  They follow from the expression
of the heat kernel as a sum over closed walks given in \cite{Ro84, KPS07}, and are stated in this form in
\cite{BHJ23}.  Here we derive them from the forest expansion of
Section~\ref{sec:forest-expansion}, without constructing the heat kernel.  The constant term
$\chi/2 = \sum_{v}(2-d_{v})/4$ is a sum of local vertex contributions, which is the origin of the
$(E_{0}-V_{0})\log k$ term in Proposition~\ref{prop:logdet-asymptotics}; see Remark~\ref{rem:friedlander}.
\end{remark}

\begin{remark}\label{rem:euler-product}
The following identity, valid for arbitrary edge lengths, is known; it appears in
\cite[(50)--(51)]{Te10}, building on \cite{ACDMT00}.  We record it because it connects the
forest expansion with periodic orbits.  Theorem~\ref{thm:heat-asymptotics} has an exact version, Roth's trace formula \cite{Ro84} (see
also \cite{KPS07}, \cite[Theorem~5.1]{BHJ23} and \cite[Chapter~8]{Ku24}).  It expresses $q(t)$ as a sum over the
primitive oriented periodic orbits $\gamma$ of $\Gamma$,
$$
q(t) = (4\pi t)^{-1/2}\sum_{\gamma}\sum_{k\geq 1}\alpha(\gamma)^{k}\ell(\gamma)\,e^{-k^{2}\ell(\gamma)^{2}/4t},
$$
where $\alpha(\gamma)$ is the product of the vertex
scattering coefficients along $\gamma$, namely $2/d_{v}$ for a transfer and $2/d_{v}-1$ for a bounce.  The
Laplace transform of $(4\pi t)^{-1/2}e^{-\ell^{2}/4t}$ is $e^{-\ell\sqrt{x}}/(2\sqrt{x})$, and
$R = \Phi'/\Phi$.  Integrating in $x$ and fixing the constant by Lemma~\ref{lem:Phi-sector} gives the
product formula
$$
\Phi(x) = \frac{\prod_{v}d_{v}}{2^{E_{0}}}\;x^{\chi/2}\,e^{L\sqrt{x}}\prod_{\gamma}
\Big(1-\alpha(\gamma)\,e^{-\ell(\gamma)\sqrt{x}}\Big),\qquad x>0\text{ large},
$$
an Euler product over primitive periodic orbits.  Equivalently, with $U(x)$ the $2E_{0}\times 2E_{0}$ matrix
indexed by the directed edges, whose entry from a directed edge $b$ ending at $v$ to a directed edge $b'$
starting at $v$ is $(2/d_{v}-\delta_{b',\bar b})\,e^{-\ell(b')\sqrt{x}}$, where $\bar b$ is the reverse
of $b$,
\begin{equation}\label{eq:secular-det}
\Phi(x) = \frac{\prod_{v}d_{v}}{2^{E_{0}}}\;x^{\chi/2}\,e^{L\sqrt{x}}\,\det\big(I-U(x)\big),
\end{equation}
which extends by analytic continuation to all $x\in\mathbb{C}\setminus(-\infty,0]$.  Deriving a trace formula by computing a characteristic function in two ways is also the approach taken in
\cite{CHJSV23} for regular graphs.
\end{remark}

\subsection{Special values of $\zeta_{\Gamma}$}

\begin{corollary}\label{cor:zeta-continuation}
The spectral zeta function $\zeta_{\Gamma}(s) = \sum_{j\geq 1}\mu_{j}^{-s}$ has a meromorphic
continuation to $\mathbb{C}$ whose only pole is a simple pole at $s=1/2$ with residue
$L/(2\pi)$.  Moreover,
$$
\zeta_{\Gamma}(0) = \frac{\chi}{2}-1,\qquad \zeta_{\Gamma}(-n) = 0\quad (n\geq 1).
$$
\end{corollary}

\begin{proof}
For $\text{\rm Re}\,s>1/2$,
$\zeta_{\Gamma}(s) = \Gamma(s)^{-1}\int_{0}^{\infty}(Q(t)-1)t^{s}\,\frac{dt}{t}$.  Split the integral at
$t=1$.  The integral over $(1,\infty)$ is entire in $s$, since $Q(t)-1 = O(e^{-\mu_{1}t})$.  On
$(0,1)$, write $Q(t)-1 = L(4\pi t)^{-1/2} + (\chi/2-1) + q(t)$.  By
Theorem~\ref{thm:heat-asymptotics}, the integral of $q(t)t^{s-1}$ over $(0,1)$ is entire, and
$$
\zeta_{\Gamma}(s) = \frac{1}{\Gamma(s)}\Big(\frac{L}{\sqrt{4\pi}\,(s-\frac12)} + \frac{\chi/2-1}{s}
+ \text{\rm (entire)}\Big).
$$
The residue at $s=1/2$ is $L/(\sqrt{4\pi}\,\Gamma(1/2)) = L/(2\pi)$.  At $s=0$, $1/(s\Gamma(s))\rightarrow 1$
gives $\zeta_{\Gamma}(0) = \chi/2-1$.  At $s=-n$ with $n\geq 1$, the function $1/\Gamma(s)$ vanishes
while the bracket is finite, so $\zeta_{\Gamma}(-n) = 0$.
\end{proof}

\subsection{The zeta-regularized determinant}

For $x>0$ let $\zeta_{\Gamma}(s;x) = \sum_{j\geq 0}(\mu_{j}+x)^{-s}$, which converges for
$\text{\rm Re}\,s>1/2$, and define
$$
\text{\rm det}(\Delta_{\Gamma}+x) = \exp\big(-\partial_{s}\zeta_{\Gamma}(s;x)|_{s=0}\big),\qquad
\text{\rm det}^{*}\Delta_{\Gamma} = \exp\big(-\zeta_{\Gamma}'(0)\big),
$$
using the continuations established below and in Corollary~\ref{cor:zeta-continuation}.

\begin{thm}\label{thm:det}
For every $x>0$,
$$
\text{\rm det}(\Delta_{\Gamma}+x) = \frac{2^{E_{0}}}{\prod_{v}d_{v}}\,\Phi(x).
$$
Consequently, $\text{\rm det}(\Delta_{\Gamma}+x)$ extends to an entire function of $x$, and
$$
\text{\rm det}^{*}\Delta_{\Gamma} = \frac{2^{E_{0}}\,P_{1}}{\prod_{v}d_{v}}
= \frac{2^{E_{0}}\,L}{\prod_{v}d_{v}}\sum_{T}\prod_{\mathcal{E}_{j}\notin T}\ell_{j}.
$$
\end{thm}

\begin{proof}
For $x>0$ and $\text{\rm Re}\,s>1/2$,
$\zeta_{\Gamma}(s;x) = \Gamma(s)^{-1}\int_{0}^{\infty}e^{-xt}Q(t)t^{s}\,dt/t$.  Writing
$Q = L(4\pi t)^{-1/2}+\chi/2+q$ and integrating the first two terms explicitly,
$$
\zeta_{\Gamma}(s;x) = \frac{L}{\sqrt{4\pi}}\frac{\Gamma(s-\frac12)}{\Gamma(s)}x^{\frac12-s}
+\frac{\chi}{2}x^{-s} + \frac{1}{\Gamma(s)}\int_{0}^{\infty}e^{-xt}q(t)\,t^{s}\frac{dt}{t}.
$$
The last integral is entire in $s$ by Theorem~\ref{thm:heat-asymptotics}, since $q(t)$ is bounded
as $t\rightarrow\infty$.  This gives the meromorphic continuation, which is regular at $s=0$.
Differentiating at $s=0$, using $1/\Gamma(s) = s+O(s^{2})$ and $\Gamma(-\frac12) = -2\sqrt{\pi}$,
$$
F(x) := -\partial_{s}\zeta_{\Gamma}(s;x)|_{s=0}
= L\sqrt{x}+\frac{\chi}{2}\log x-\int_{0}^{\infty}e^{-xt}q(t)\,\frac{dt}{t}.
$$
By dominated convergence the last integral tends to $0$ as $x\rightarrow+\infty$.  Hence
$F(x) = L\sqrt{x}+\frac{\chi}{2}\log x+o(1)$ as $x\rightarrow+\infty$.

Next we compute $F''$.  The representation above shows that $\zeta_{\Gamma}(s;x)$ is jointly holomorphic in
$(s,x)$ for $s$ away from $\{\frac12,-\frac12,-\frac32,\dots\}$ and $\text{\rm Re}\,x>0$, since the
first two terms are explicit, and the integral may be differentiated under the integral sign.
For $\text{\rm Re}\,s>1/2$, termwise differentiation of the series gives
$\partial_{x}\zeta_{\Gamma}(s;x) = -s\,\zeta_{\Gamma}(s+1;x)$.  For fixed $x$, both sides are
meromorphic in $s$, so the identity persists after continuation.  Hence
$\partial_{x}^{2}\zeta_{\Gamma}(s;x) = s(s+1)\zeta_{\Gamma}(s+2;x)$, where $\zeta_{\Gamma}(s+2;x)$ is
holomorphic near $s=0$.  Differentiating in $s$ at $s=0$, and using that $\zeta_{\Gamma}(s;x)$ is
jointly holomorphic so the derivatives may be interchanged,
$$
F''(x) = -\zeta_{\Gamma}(2;x) = -\sum_{j\geq 0}(\mu_{j}+x)^{-2}.
$$
On the other hand, by Theorem~\ref{thm:hadamard},
$(\log\Phi)''(x) = -\sum_{j\geq 0}(\mu_{j}+x)^{-2}$ for $x>0$.  Hence $F(x)-\log\Phi(x) = a+bx$ for
constants $a$ and $b$.  Comparing the asymptotics of $F$ with Lemma~\ref{lem:Phi-sector} as
$x\rightarrow+\infty$ gives $b=0$ and $a = \log(2^{E_{0}}/\prod_{v}d_{v})$, which proves the first
assertion.

For the second, write $\zeta_{\Gamma}(s;x) = x^{-s}+\zeta^{*}_{\Gamma}(s;x)$ with
$\zeta^{*}_{\Gamma}(s;x) = \sum_{j\geq 1}(\mu_{j}+x)^{-s}$.  Let $0<x<\mu_{1}$ and $\text{\rm Re}\,s>1/2$.  Expanding
$(\mu_{j}+x)^{-s} = \mu_{j}^{-s}(1+x/\mu_{j})^{-s}$ by the binomial series and summing over $j$, which is permitted by
absolute convergence, gives
$$
\zeta^{*}_{\Gamma}(s;x) = \zeta_{\Gamma}(s)+\sum_{k\geq 1}\binom{-s}{k}x^{k}\,\zeta_{\Gamma}(s+k).
$$
For $|s|\leq\frac14$ and $k\geq 1$ we have $\mu_{j}^{-\text{\rm Re}\,s-k}\leq\mu_{1}^{-(k-1)}\max(\mu_{j}^{-3/4},\mu_{j}^{-5/4})$, hence
$|\zeta_{\Gamma}(s+k)|\leq\mu_{1}^{-(k-1)}\big(\zeta_{\Gamma}(\frac34)+\zeta_{\Gamma}(\frac54)\big)$; also $|\binom{-s}{k}|\leq 1$.  So for $|x|\leq\rho<\mu_{1}$ the series converges uniformly, and its sum is holomorphic in
$(s,x)$ on $\{|s|<\frac14\}\times\{|x|<\rho\}$.  By Corollary~\ref{cor:zeta-continuation}, $\zeta_{\Gamma}(s)$ is holomorphic
near $s=0$.  For fixed $x>0$ both sides of the identity are meromorphic in $s$, so the identity persists after
continuation.  Since $\frac{d}{ds}\binom{-s}{k}\big|_{s=0} = (-1)^{k}/k$ for $k\geq 1$, we obtain
$$
\partial_{s}\zeta^{*}_{\Gamma}(s;x)|_{s=0} = \zeta_{\Gamma}'(0)+\sum_{k\geq 1}\frac{(-1)^{k}\zeta_{\Gamma}(k)}{k}\,x^{k},
$$
which is continuous at $x = 0$ with value $\zeta_{\Gamma}'(0)$.  Since $F(x) = \log x-\partial_{s}\zeta^{*}_{\Gamma}(s;x)|_{s=0}$, this gives
$$
-\zeta_{\Gamma}'(0) = \lim_{x\rightarrow 0^{+}}\big(F(x)-\log x\big)
= \lim_{x\rightarrow 0^{+}}\log\Big(\frac{2^{E_{0}}}{\prod_{v}d_{v}}\frac{\Phi(x)}{x}\Big)
= \log\Big(\frac{2^{E_{0}}P_{1}}{\prod_{v}d_{v}}\Big),
$$
since $\Phi(x)/x\rightarrow\Phi'(0) = P_{1}$.  Comparing the Taylor expansions of $F(x)-\log x$ and of
$\log(\Phi(x)/x)$ at $x=0$ also recovers Corollary~\ref{cor:zeta-integers}.
\end{proof}

\begin{corollary}\label{cor:friedlander}
Let $\Delta_{D}$ be the Laplacian on $\Gamma$ with Dirichlet conditions at all vertices, so that $\Delta_{D}$ is the direct sum
of the Dirichlet Laplacians of the edges.  Then, for every $x>0$,
$$
\text{\rm det}(\Delta_{\Gamma}+x) = \frac{1}{\prod_{v}d_{v}}\;\text{\rm det}(\Delta_{D}+x)\;\det M(x),
$$
where $M(x)$ is the Dirichlet-to-Neumann matrix of Section~\ref{subsec:DtN}.  Moreover,
$$
\text{\rm det}^{*}\Delta_{\Gamma} = \frac{2^{E_{0}}\,L\prod_{j}\ell_{j}\;\text{\rm det}^{*}M(0)}{V_{0}\prod_{v}d_{v}},
$$
where $\text{\rm det}^{*}M(0)$ is the product of the nonzero eigenvalues of $M(0)$, the Laplacian of $G_{0}$ with edge weights
$1/\ell_{j}$.
\end{corollary}

\begin{proof}
The zeta function of a direct sum of operators is the sum of their zeta functions, so the zeta-regularized
determinant of $\Delta_{D}+x$ is the product of those of the edges.  By Remark~\ref{rem:motivation-S}, the determinant of
the Dirichlet Laplacian of $\mathcal{E}_{j}$ is $2S_{j}(x)$, so $\text{\rm det}(\Delta_{D}+x) = 2^{E_{0}}\prod_{j}S_{j}(x)$.  By
Proposition~\ref{prop:DtN-properties}(4), $\Phi(x) = \prod_{j}S_{j}(x)\det M(x)$ for $x>0$.  Substituting both into
Theorem~\ref{thm:det} gives the first formula.  For the second, Theorem~\ref{thm:det} gives
$\text{\rm det}^{*}\Delta_{\Gamma} = 2^{E_{0}}P_{1}/\prod_{v}d_{v}$ with $P_{1} = L\prod_{j}\ell_{j}\,\tau(G_{0},\Lambda^{-1})$, by
Corollary~\ref{cor:kirchhoff-type}.  By Proposition~\ref{prop:DtN-properties}(2) and the weighted matrix-tree
theorem, $\text{\rm det}^{*}M(0) = V_{0}\,\tau(G_{0},\Lambda^{-1})$.
\end{proof}

\begin{remark}\label{rem:friedlander-thm1}
Corollary~\ref{cor:friedlander} is the case $q=0$, $A=0$ of \cite[Theorem~1]{Fr06}, the
Dirichlet-to-Neumann matrix $R(\lambda)$ there corresponding to our $M(\lambda)$ and the Dirichlet
determinant of each edge being $2S_{j}(\lambda)$.  It agrees with \cite[Theorem~2]{Fr06} because
$\text{\rm det}'R(0) = V_{0}\,\tau(G_{0},\Lambda^{-1})$ by the weighted matrix-tree theorem.  Both proofs
compare logarithmic derivatives and then fix the constant from the large-$\lambda$ asymptotics.  In \cite[(51)]{Te10} the product $\prod_{j}S_{j}\det M$ appears as the vertex determinant of
\cite{PM99}.  The new ingredient here is that $\Phi$ arises as the limit of
discrete determinants on $G_{k}$, through the exact forest expansion of
Section~\ref{sec:forest-expansion}.
\end{remark}

\begin{remark}\label{rem:regularized-products}
The arguments of this section are an instance of the theory of regularized products and series of
\cite{JoLa93}, with $\{\mu_{j}\}$ a sequence whose theta function $Q(t)$ has the expansion of
Theorem~\ref{thm:heat-asymptotics}, $R(x)$ its harmonic series and
$\text{\rm det}(\Delta_{\Gamma}+x)$ its regularized product.  That theory expresses the
regularized product as a Weierstrass product times the exponential of a polynomial; here the polynomial
is the constant $\log(2^{E_{0}}/\prod_{v}d_{v})$, because $\Phi$ has order $1/2$.
\end{remark}

\subsection{A Cheeger--M\"uller type formulation}\label{subsec:CM}

The spanning-tree sum in Theorem~\ref{thm:det} has a cohomological meaning.  Choose a basis
$c_{1},\dots,c_{\beta}$ of the integral cycle space $H_{1}(G_{0},\mathbb{Z})$, where $\beta = E_{0}-V_{0}+1$ is
the first Betti number.  Each $c_{i}$ is an integer combination of oriented edges.  The period matrix of
$\Gamma$ is the $\beta\times\beta$ matrix
$$
(Q_{\Gamma})_{ik} = \sum_{j=1}^{E_{0}}\ell_{j}\,\langle c_{i},\mathcal{E}_{j}\rangle\langle c_{k},\mathcal{E}_{j}\rangle,
$$
the Gram matrix of the cycle basis in the inner product defined by the edge lengths.  The torus
$H_{1}(\Gamma,\mathbb{R})/H_{1}(\Gamma,\mathbb{Z})$ with this inner product is the tropical Jacobian of $\Gamma$
\cite{KS00, MZ08}, and its volume is $(\det Q_{\Gamma})^{1/2}$.  If $G_{0}$ is a tree, then $\beta = 0$, $Q_{\Gamma}$ is the empty
matrix, and we set $\det Q_{\Gamma} = 1$; this agrees with the lemma below, since the only spanning tree is $G_{0}$
itself and the empty product equals $1$.

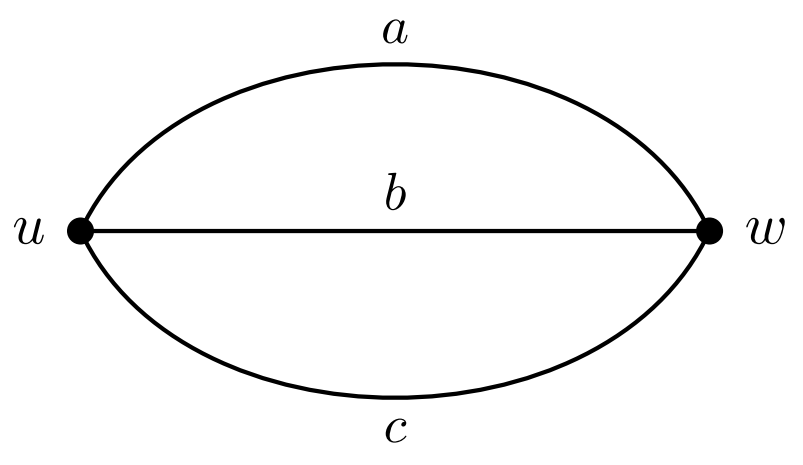
\begin{figure}[ht]
\centering
\begin{tikzpicture}[scale=1.3, v/.style={circle,fill,inner sep=1.8pt}, e/.style={thick}, lab/.style={font=\small}]
  \coordinate (U) at (-1.6,0); \coordinate (W) at (1.6,0);
  \draw[e] (U) to[bend left=65] node[lab,above] {$a$} (W);
  \draw[e] (U) -- node[lab,above] {$b$} (W);
  \draw[e] (U) to[bend right=65] node[lab,below] {$c$} (W);
  \node[v,label=left:{$u$}] at (U) {};
  \node[v,label=right:{$w$}] at (W) {};
\end{tikzpicture}
\caption{The theta graph $D_3$. Two vertices $u$ and $w$ are joined by three edges of lengths $a$, $b$ and $c$.}
\label{fig:theta}
\end{figure}

As an example, for the theta graph with lengths $a,b,c$ and cycles $\mathcal{E}_{1}-\mathcal{E}_{2}$, $\mathcal{E}_{2}-\mathcal{E}_{3}$,
one gets that
$Q_{\Gamma} = \begin{pmatrix}a+b & -b\\ -b & b+c\end{pmatrix}$, with determinant $ab+bc+ca$.

\begin{lem}\label{lem:period} We have that
$\displaystyle\det Q_{\Gamma} = \sum_{T}\prod_{\mathcal{E}_{j}\notin T}\ell_{j}$, where the sum is over spanning trees of
$G_{0}$.  In particular, $\det Q_{\Gamma}$ does not depend on the choice of integral cycle basis.
\end{lem}

\begin{proof}
Let $C$ be the $\beta\times E_{0}$ matrix with entries $\langle c_{i},\mathcal{E}_{j}\rangle$, so that
$Q_{\Gamma} = C\,\text{\rm diag}(\ell)\,C^{T}$.  By the Cauchy--Binet formula,
$\det Q_{\Gamma} = \sum_{S}(\det C_{S})^{2}\prod_{j\in S}\ell_{j}$, where $S$ runs over $\beta$-element sets of edges
and $C_{S}$ is the corresponding $\beta\times\beta$ minor.

We first show that $C_{S}$ is nonsingular if and only if the complement $T = E(G_{0})\setminus S$ is a spanning
tree.  The rows of $C$ are a basis of the real cycle space $Z = H_{1}(G_{0},\mathbb{R})$, which consists of
the real combinations of oriented edges with zero boundary.  So $C_{S}$ is nonsingular if and only if the
restriction map $Z\rightarrow\mathbb{R}^{S}$, which records the coefficients on the edges of $S$, is
injective; both spaces have dimension $\beta$.  The kernel of this map consists of the cycles supported on
$T$.  So $C_{S}$ is nonsingular if and only if $T$ contains no cycle, where a loop counts as a cycle.  Since $T$ has
$E_{0}-\beta = V_{0}-1$ edges, it is acyclic if and only if it is a spanning tree.

Now let $T$ be a spanning tree.  For each $e\in S$, the graph $T\cup\{e\}$ contains a unique cycle $f_{e}$,
the fundamental cycle, in which $e$ has coefficient $1$.  The $f_{e}$ with $e\in S$ form a
$\mathbb{Z}$-basis of $H_{1}(G_{0},\mathbb{Z})$.  Indeed, for an integral cycle $z$ with coefficients
$z_{e}$ on the edges $e\in S$, the cycle $z-\sum_{e\in S}z_{e}f_{e}$ is supported on $T$, and is therefore zero.
Relative to this basis, the minor on $S$ is the identity matrix, since $f_{e}$ contains no edge of $S$
other than $e$.  Since $c_{1},\dots,c_{\beta}$ is another $\mathbb{Z}$-basis, the two bases differ by a matrix
$U\in GL_{\beta}(\mathbb{Z})$, and $\det C_{S} = \det U = \pm 1$.  Substituting into the Cauchy--Binet formula,
with $T$ the complement of $S$, proves the lemma.
\end{proof}

Since $\text{\rm vol}(\Gamma) = L$, Theorem~\ref{thm:det} can therefore be written as follows.

\begin{corollary}\label{cor:CM}
For every compact connected metric graph $\Gamma$ with standard vertex conditions,
$$
\text{\rm det}^{*}\Delta_{\Gamma} \;=\; \prod_{v\in V(G_{0})}\frac{2^{d_{v}/2}}{d_{v}}\;\cdot\;
\text{\rm vol}(\Gamma)\;\cdot\;\det Q_{\Gamma},
$$
where the first factor is a product of local contributions from the vertices, and
$\text{\rm vol}(\Gamma)\det Q_{\Gamma}$ is the product of the squared volumes of $H^{0}$ and $H^{1}$ of $\Gamma$ in
the metrics induced by the edge lengths.
\end{corollary}

\begin{proof}
Combine Theorem~\ref{thm:det} with Lemma~\ref{lem:period}, and use
$2^{E_{0}} = \prod_{v}2^{d_{v}/2}$, which holds because $\sum_{v}d_{v} = 2E_{0}$.
\end{proof}

\begin{remark}\label{rem:HW-trees}
Corollary~\ref{cor:CM} can be read as a spectral form of Kirchhoff's theorem for metric graphs
$$
\frac{\prod_{v}d_{v}}{2^{E_{0}}\,L}\,\text{\rm det}^{*}\Delta_{\Gamma} = \sum_{T}\prod_{\mathcal{E}_{j}\notin T}\ell_{j},
$$
the length-weighted count of spanning trees, for arbitrary edge lengths.  For an equilateral graph with edge
length $\ell$ the right-hand side is $\ell^{\beta}$ times the number of spanning trees, and $L = E_{0}\ell$.  This recovers the
theorem \cite[Theorem~3]{HW23} that the number of spanning trees equals
$\prod_{v}d_{v}\,\text{\rm det}^{*}\Delta_{\Gamma}/(E_{0}2^{E_{0}}\ell^{\beta+1})$.  It is also shown in \cite{HW23}, by a perturbation
argument, that the nearest integer to this expression is still the number of spanning trees when all edge
lengths lie in a sufficiently short interval $[\ell,\ell+\delta)$ \cite[Theorem~2]{HW23}.  The displayed identity gives the
weighted count exactly, without any restriction on the lengths.  The unweighted count is recovered precisely
when the weights $\prod_{\mathcal{E}_{j}\notin T}\ell_{j}$ are nearly equal.
\end{remark}

Corollary~\ref{cor:CM} can be written as
$$
\frac{\text{\rm det}^{*}\Delta_{\Gamma}}{\text{\rm vol}(\Gamma)\cdot\det Q_{\Gamma}} = \prod_{v}\frac{2^{d_{v}/2}}{d_{v}}.
$$
Thus the left-hand side has the structure of a Quillen metric, an $L^{2}$-metric
on the determinant of cohomology corrected by a regularized determinant \cite{Qu85}.

\begin{remark}\label{rem:CM}
In this form Corollary~\ref{cor:CM} has the shape of a Cheeger--M\"uller theorem \cite{RS71, Ch79, Mu78}, with the
vertex factors $2^{d_{v}/2}/d_{v}$ in the role of the anomaly terms for singular spaces.  These factors equal $1$
exactly when $d_{v} = 2$, in accordance with the subdivision invariance of Lemma~\ref{lem:subdivision}.  Our approximation by the
refinements $G_{k}$ parallels the approximation of the Hodge Laplacian by combinatorial Laplacians of refined
triangulations in \cite{Do76}.  The Cheeger--M\"uller theorem concerns torsion, but on a one-dimensional space the Laplacians on
functions and on $1$-forms have the same nonzero spectrum, so the torsion reduces essentially to
$\text{\rm det}^{*}\Delta_{\Gamma}$ and the cohomology terms; we do not match the normalizations of Ray--Singer torsion exactly.
\end{remark}

\begin{corollary}\label{cor:main}
As $k\rightarrow\infty$,
$$
k^{V_{0}-E_{0}}\det\big(\Delta_{G_{k}}+x/k^{2}\big)\longrightarrow
\frac{\prod_{v}d_{v}}{2^{E_{0}}}\,\text{\rm det}(\Delta_{\Gamma}+x)
$$
uniformly for $x$ in compact subsets of $\mathbb{C}$.  In particular,
$$
\log\text{\rm det}^{*}\Delta_{G_{k}} = (E_{0}-V_{0}+2)\log k + \log\text{\rm det}^{*}\Delta_{\Gamma}
- E_{0}\log 2 + \sum_{v}\log d_{v} + o(1).
$$
\end{corollary}

\begin{proof}
The first statement combines Theorems~\ref{thm:Phi-limit} and~\ref{thm:det}.  The second combines
Proposition~\ref{prop:logdet-asymptotics} with Theorem~\ref{thm:det}.
\end{proof}

\section{Explicit formulas for the special values}\label{sec:explicit}

In this section we make the results of Sections~\ref{sec:forest-expansion}--\ref{sec:det} completely
explicit.  We first write each forest invariant $P_{m}$ as a finite sum over spanning forests of
$G_{0}$, with coefficients that are explicit monomials in the edge lengths.  We then give the formulas,
from Newton's identities, that express $\zeta_{\Gamma}(r)$ in terms of $P_{1},\dots,P_{r+1}$.  Finally we
summarize the procedure that produces every special value of $\zeta_{\Gamma}$ at the integers.

\subsection{Notation}

Let $F$ be a spanning forest of $G_{0}$, with components $C$.  We write
$$
w_{F} = \prod_{\mathcal{E}_{j}\notin F}\ell_{j},\qquad
\sigma_{p}(F) = \sum_{\mathcal{E}_{j}\notin F}\ell_{j}^{\,p},\qquad
\Lambda_{p}(C) = \sum_{h\text{ \rm at }C}\ell_{e(h)}^{\,p},\qquad
\ell^{*}(C) = \tfrac12\Lambda_{1}(C).
$$
Thus $w_{F}$ is the product of the lengths of the edges not in $F$, and $\sigma_{p}(F)$ is the sum of
their $p$-th powers.  The quantity $\Lambda_{p}(C)$ is the sum of $\ell^{p}$ over the half-edges at
vertices of $C$, so an edge with both ends in $C$ is counted twice.  Finally, $\ell^{*}(C)$ is the
metric size of $C$, as in Corollary~\ref{cor:kirchhoff-type}.  If $T$ is a spanning tree, it has a
single component, and $\Lambda_{p}(T) = 2\sum_{j=1}^{E_{0}}\ell_{j}^{\,p}$ and $\ell^{*}(T) = L$.  We
also use the rational numbers
$$
a_{n} = \frac{1}{(2n+1)!}\quad(n\geq 0),\qquad
b_{n} = \frac{2(2^{2n}-1)B_{2n}}{(2n)!}\quad(n\geq 1),
$$
where $B_{2n}$ are the Bernoulli numbers.  They are the Taylor coefficients in
$\sinh(z)/z = \sum_{n}a_{n}z^{2n}$ and $z\tanh(z/2) = \sum_{n\geq 1}b_{n}z^{2n}$; for example,
$b_{1} = \frac12$, $b_{2} = -\frac{1}{24}$ and $b_{3} = \frac{1}{240}$.

\subsection{The forest invariants in terms of the edge lengths}

\begin{thm}\label{thm:explicit-P}
For every $m\geq 1$,
$$
P_{m} = \sum_{F}\;w_{F}\sum\;
\prod_{\mathcal{E}_{j}\notin F}a_{n_{j}}\,\ell_{j}^{\,2n_{j}}\;
\prod_{C\in F}b_{q_{C}}\,\ell_{e(h_{C})}^{\,2q_{C}-1}.
$$
Here the outer sum is over all spanning forests $F$ of $G_{0}$.  The inner sum is over the following data.
Each edge $\mathcal{E}_{j}$ not in $F$ is given an integer $n_{j}\geq 0$.  Each component $C$ of $F$ is given an integer
$q_{C}\geq 1$ and a half-edge $h_{C}$ at one of its vertices.  These choices are subject to the single condition
$$
\sum_{\mathcal{E}_{j}\notin F}n_{j}+\sum_{C\in F}q_{C} = m.
$$
In particular, only forests with at most $m$ components contribute to $P_{m}$.  Moreover, $P_{m}$ is a
homogeneous polynomial of degree $E_{0}-V_{0}+2m$ in the edge lengths, with rational coefficients.
\end{thm}

\begin{proof}
By Theorem~\ref{thm:Phi-limit}, $P_{m}$ is the coefficient of $x^{m}$ in the expansion \eqref{eq:Phi}
of $\Phi$.  Insert $S_{j}(x) = \ell_{j}\sum_{n\geq 0}a_{n}\ell_{j}^{2n}x^{n}$ and
$\Theta_{C}(x) = \sum_{h\text{ \rm at }C}\sum_{q\geq 1}b_{q}\ell_{e(h)}^{2q-1}x^{q}$ into \eqref{eq:Phi}, and collect the
coefficient of $x^{m}$.  A forest with $r$ components contributes only in degrees $x^{m}$ with
$m\geq r$, since each $\Theta_{C}$ is divisible by $x$.  For homogeneity, note that a forest with $r$
components has $E_{0}-V_{0}+r$ edges outside it.  Each term therefore has degree
$(E_{0}-V_{0}+r)+\sum_{j}2n_{j}+\sum_{C}(2q_{C}-1) = E_{0}-V_{0}+2m$.
\end{proof}

\begin{corollary}\label{cor:first-three}
For the first three invariants, Theorem~\ref{thm:explicit-P} reads as follows.  Here $T$ runs over
spanning trees, $F$ over spanning forests with two components $C_{1},C_{2}$, and $F'$ over spanning
forests with three components $C_{1}',C_{2}',C_{3}'$.  We also write
$\epsilon_{2}(T) = \sum\ell_{i}^{2}\ell_{j}^{2}$, summed over pairs of distinct edges not in $T$.
\begin{align*}
P_{1} &= L\sum_{T}w_{T},\\[4pt]
P_{2} &= \sum_{T}w_{T}\Big(\frac{L}{6}\,\sigma_{2}(T)-\frac{1}{12}\sum_{j=1}^{E_{0}}\ell_{j}^{3}\Big)
+\sum_{F}w_{F}\,\ell^{*}(C_{1})\,\ell^{*}(C_{2}),\\[4pt]
P_{3} &= \sum_{T}w_{T}\Big(L\Big(\frac{\sigma_{4}(T)}{120}+\frac{\epsilon_{2}(T)}{36}\Big)
-\frac{\sigma_{2}(T)}{72}\sum_{j=1}^{E_{0}}\ell_{j}^{3}+\frac{1}{120}\sum_{j=1}^{E_{0}}\ell_{j}^{5}\Big)\\
&\quad+\sum_{F}w_{F}\Big(\frac{\sigma_{2}(F)}{6}\,\ell^{*}(C_{1})\,\ell^{*}(C_{2})
-\frac{1}{24}\big(\ell^{*}(C_{1})\Lambda_{3}(C_{2})+\ell^{*}(C_{2})\Lambda_{3}(C_{1})\big)\Big)\\
&\quad+\sum_{F'}w_{F'}\,\ell^{*}(C_{1}')\,\ell^{*}(C_{2}')\,\ell^{*}(C_{3}').
\end{align*}
\end{corollary}

\begin{proof}
In Theorem~\ref{thm:explicit-P}, the sum over half-edges at $C$ of $b_{q}\ell_{e(h)}^{2q-1}$ equals
$b_{q}\Lambda_{2q-1}(C)$.  In particular, $b_{1}\Lambda_{1}(C) = \ell^{*}(C)$.  For a forest with $r$ components
we list the choices of the $n_{j}\geq 0$ and $q_{C}\geq 1$ with $\sum_{j}n_{j}+\sum_{C}q_{C} = m$.

For $m=1$, only spanning trees contribute, with all $n_{j} = 0$ and $q = 1$.  This gives
$w_{T}\,\ell^{*}(T) = w_{T}L$.

For $m=2$ and a spanning tree $T$ there are two choices.  Either $n_{j} = 1$ for a single edge $j\notin T$
and $q=1$, contributing $w_{T}\,a_{1}\sigma_{2}(T)\,\ell^{*}(T) = w_{T}\,L\,\sigma_{2}(T)/6$.  Or all
$n_{j} = 0$ and $q = 2$, contributing $w_{T}\,b_{2}\Lambda_{3}(T) = -w_{T}\sum_{j}\ell_{j}^{3}/12$.  For a
two-component forest the only choice is $n_{j} = 0$ and $q_{C_{1}} = q_{C_{2}} = 1$, contributing
$w_{F}\,\ell^{*}(C_{1})\,\ell^{*}(C_{2})$.

For $m=3$ and a spanning tree there are three choices.  With $q=1$ and $\sum n_{j} = 2$, either one
$n_{j} = 2$, contributing $a_{2}\sigma_{4}(T) = \sigma_{4}(T)/120$, or two $n_{j} = 1$, contributing
$a_{1}^{2}\epsilon_{2}(T) = \epsilon_{2}(T)/36$; both are multiplied by $w_{T}L$.  With $q=2$ and one $n_{j} = 1$,
the contribution is
$$
w_{T}\,b_{2}\Lambda_{3}(T)\,a_{1}\sigma_{2}(T) = -\frac{1}{72}\,w_{T}\,\sigma_{2}(T)\sum_{j}\ell_{j}^{3}.
$$
With $q=3$
and all $n_{j} = 0$, it is $w_{T}\,b_{3}\Lambda_{5}(T) = w_{T}\sum_{j}\ell_{j}^{5}/120$.  For a two-component forest,
either one $n_{j} = 1$ and $q_{C_{1}} = q_{C_{2}} = 1$, contributing $w_{F}\,\sigma_{2}(F)\,\ell^{*}(C_{1})\ell^{*}(C_{2})/6$, or all
$n_{j} = 0$ and $(q_{C_{1}},q_{C_{2}})\in\{(2,1),(1,2)\}$, contributing
$$
w_{F}\,b_{2}\big(\Lambda_{3}(C_{1})\ell^{*}(C_{2})+\Lambda_{3}(C_{2})\ell^{*}(C_{1})\big).
$$
For a three-component forest the
contribution is $w_{F'}\prod_{i}\ell^{*}(C_{i}')$.  Here we used $a_{1} = \frac16$, $a_{2} = \frac{1}{120}$,
$b_{2} = -\frac{1}{24}$, $b_{3} = \frac{1}{240}$, and $\Lambda_{p}(T) = 2\sum_{j}\ell_{j}^{p}$.
\end{proof}

\subsection{From the forest invariants to the zeta values by Newton's identities}

Set $e_{0} = 1$ and
$$
e_{m} = \frac{P_{m+1}}{P_{1}}\qquad(m\geq 1).
$$
By Theorem~\ref{thm:hadamard}, $\Phi(x)/(P_{1}x) = \prod_{j\geq 1}(1+x/\mu_{j}) = \sum_{m\geq 0}e_{m}x^{m}$.  So
$e_{m}$ is the elementary symmetric function of degree $m$ in the reciprocals $1/\mu_{j}$
$$
e_{m} = \sum_{1\leq j_{1}<j_{2}<\cdots<j_{m}}\frac{1}{\mu_{j_{1}}\mu_{j_{2}}\cdots\mu_{j_{m}}},
$$
while $\zeta_{\Gamma}(r) = \sum_{j}\mu_{j}^{-r}$ is the $r$-th power sum of the same numbers.  The classical
relations between elementary symmetric functions and power sums, which remain valid here because all
the series converge absolutely, give the following.

\begin{thm}\label{thm:newton}
For every integer $r\geq 1$ the following hold.
\begin{enumerate}
\item (Newton's recursion)
$\displaystyle \zeta_{\Gamma}(r) = \sum_{i=1}^{r-1}(-1)^{i-1}e_{i}\,\zeta_{\Gamma}(r-i)+(-1)^{r-1}\,r\,e_{r}$.
\item (Girard--Waring formula)
$$
\zeta_{\Gamma}(r) = \sum_{m_{1}+2m_{2}+\cdots+rm_{r} = r}
(-1)^{r+m_{1}+\cdots+m_{r}}\,\frac{r\,(m_{1}+\cdots+m_{r}-1)!}{m_{1}!\cdots m_{r}!}\;
e_{1}^{m_{1}}\cdots e_{r}^{m_{r}},
$$
the sum being over all $(m_{1},\dots,m_{r})\in\mathbb{Z}_{\geq 0}^{r}$ with $\sum_{i}i\,m_{i} = r$.
\item (Determinant formula)
$$
\zeta_{\Gamma}(r) = \det\begin{pmatrix}
e_{1} & 1 & 0 & \cdots & 0\\
2e_{2} & e_{1} & 1 & \cdots & 0\\
3e_{3} & e_{2} & e_{1} & \ddots & \vdots\\
\vdots & \vdots & & \ddots & 1\\
r\,e_{r} & e_{r-1} & \cdots & e_{2} & e_{1}
\end{pmatrix}.
$$
\item (Generating function) For $|x|<\mu_{1}$,
$\displaystyle \sum_{r\geq 1}\frac{(-1)^{r-1}\zeta_{\Gamma}(r)}{r}\,x^{r} = \log\Big(\sum_{m\geq 0}e_{m}x^{m}\Big)$.
\end{enumerate}
\end{thm}

The formula in part (2) is classical.  Its first cases go back to Girard, in his \emph{Invention nouvelle en
l'alg\`ebre} of 1629, and the general case to Waring, in the \emph{Meditationes Algebraicae} of 1770, whence the name.
Its history, and its use for linear recurrences, are discussed in \cite{Go99}.

\begin{proof}
Part (4) is Corollary~\ref{cor:zeta-integers}.  Differentiating (4) and comparing coefficients gives
(1).  Parts (2) and (3) are the standard solutions of the recursion (1); see, for example,
\cite{Mac95}.  Parts (1) and (2) here are parts (4) and (5) of Theorem~\ref{mthm:special}, whose part (3) is
Corollary~\ref{cor:zeta-integers} and whose parts (1) and (2) are Corollary~\ref{cor:zeta-continuation}.  This
completes the proof of Theorem~\ref{mthm:special}.
\end{proof}

Explicitly, for the first four values,
\begin{align*}
\zeta_{\Gamma}(1) &= e_{1},\\
\zeta_{\Gamma}(2) &= e_{1}^{2}-2e_{2},\\
\zeta_{\Gamma}(3) &= e_{1}^{3}-3e_{1}e_{2}+3e_{3},\\
\zeta_{\Gamma}(4) &= e_{1}^{4}-4e_{1}^{2}e_{2}+2e_{2}^{2}+4e_{1}e_{3}-4e_{4}.
\end{align*}
Combining with Theorem~\ref{thm:explicit-P}, we obtain $\zeta_{\Gamma}(1)$ directly in terms of the edge
lengths and the spanning trees and two-component forests of $G_{0}$
$$
\zeta_{\Gamma}(1) = \frac{\displaystyle\sum_{T}w_{T}\Big(\frac{L}{6}\,\sigma_{2}(T)-\frac{1}{12}\sum_{j}\ell_{j}^{3}\Big)
+\sum_{F}w_{F}\,\ell^{*}(C_{1})\,\ell^{*}(C_{2})}{\displaystyle L\sum_{T}w_{T}}.
$$
Similarly, $\zeta_{\Gamma}(2) = (P_{2}/P_{1})^{2}-2P_{3}/P_{1}$ is given in terms of spanning forests with at
most three components by the formulas above.  In general, $\zeta_{\Gamma}(r)$ involves exactly the
spanning forests of $G_{0}$ with at most $r+1$ components.

\subsection{Spanning forest counts in terms of special values}\label{subsec:forest-counts}

Theorem~\ref{thm:newton} expresses the special values through the forest invariants.  The relation can be
inverted.  The forest invariants, and hence the spanning forest counts, are polynomials in the special values.
For $n\geq 0$, let $h_{n}(p_{1},\dots,p_{n})$ be the polynomial expressing the elementary symmetric function of
degree $n$ through the power sums,
$$
h_{n}(p_{1},\dots,p_{n}) = \sum_{(m_{1},\dots,m_{n})}\ \prod_{i=1}^{n}
\frac{(-1)^{(i-1)m_{i}}\,p_{i}^{m_{i}}}{i^{m_{i}}\,m_{i}!},
$$
the sum being over all $(m_{1},\dots,m_{n})\in\mathbb{Z}_{\geq 0}^{n}$ with $m_{1}+2m_{2}+\cdots+nm_{n} = n$.  Thus $h_{0} = 1$ and
\begin{align*}
h_{1} &= p_{1}, \qquad h_{2} = \tfrac12\big(p_{1}^{2}-p_{2}\big), \qquad
h_{3} = \tfrac16\big(p_{1}^{3}-3p_{1}p_{2}+2p_{3}\big),\\
h_{4} &= \tfrac{1}{24}\big(p_{1}^{4}-6p_{1}^{2}p_{2}+3p_{2}^{2}+8p_{1}p_{3}-6p_{4}\big).
\end{align*}
Equivalently, $\sum_{n\geq 0}h_{n}x^{n} = \exp\big(\sum_{r\geq 1}(-1)^{r-1}p_{r}x^{r}/r\big)$ as formal power series.  These
formulas are classical; see \cite[Chapter~I, \S2]{Mac95}, where our $h_{n}$ is the elementary symmetric function
$e_{n}$ written in terms of power sums.

\begin{thm}\label{thm:forest-counts}
\begin{enumerate}
\item For every finite connected graph $G$, possibly with loops and multiple edges, and every $m\geq 1$,
$$
c_{m}(G) = \text{\rm det}^{*}\Delta_{G}\cdot h_{m-1}\big(\zeta_{G}(1),\dots,\zeta_{G}(m-1)\big).
$$
\item For every compact metric graph $\Gamma = (G_{0},\Lambda)$ and every $m\geq 1$,
$$
P_{m} = \frac{\prod_{v}d_{v}}{2^{E_{0}}}\,\text{\rm det}^{*}\Delta_{\Gamma}\cdot h_{m-1}\big(\zeta_{\Gamma}(1),\dots,\zeta_{\Gamma}(m-1)\big).
$$
Consequently, the weighted numbers of $m$-spanning forests of the discrete approximations satisfy
$$
c_{m}(G_{k}) = k^{E_{0}-V_{0}+2m}\Big(\frac{\prod_{v}d_{v}}{2^{E_{0}}}\,\text{\rm det}^{*}\Delta_{\Gamma}\cdot
h_{m-1}\big(\zeta_{\Gamma}(1),\dots,\zeta_{\Gamma}(m-1)\big)+o(1)\Big)\qquad(k\rightarrow\infty).
$$
\item Equivalently, for $|x|<\mu_{1}$,
$$
\Phi(x) = \frac{\prod_{v}d_{v}}{2^{E_{0}}}\,\text{\rm det}^{*}\Delta_{\Gamma}\cdot x\,
\exp\Big(\sum_{r\geq 1}\frac{(-1)^{r-1}\zeta_{\Gamma}(r)}{r}\,x^{r}\Big).
$$
\end{enumerate}
\end{thm}

\begin{proof}
For (1), by \eqref{eq:forest-zeta}, $c_{m}(G) = \text{\rm det}^{*}\Delta_{G}\cdot e_{m-1}(\lambda_{1}^{-1},\dots,\lambda_{V-1}^{-1})$.
The elementary symmetric functions of finitely many numbers $y_{h}$ are the coefficients of
$$
\prod_{h}(1+y_{h}x) = \exp\Big(\sum_{r\geq 1}\frac{(-1)^{r-1}p_{r}}{r}\,x^{r}\Big),\qquad p_{r} = \sum_{h}y_{h}^{r}.
$$
With
$y_{h} = \lambda_{h}^{-1}$ we have $p_{r} = \zeta_{G}(r)$, so $e_{m-1} = h_{m-1}(\zeta_{G}(1),\dots,\zeta_{G}(m-1))$.

For (2) and (3), by Theorem~\ref{thm:hadamard} and Corollary~\ref{cor:zeta-integers}, for $|x|<\mu_{1}$,
$$
\frac{\Phi(x)}{P_{1}x} = \prod_{j\geq 1}\Big(1+\frac{x}{\mu_{j}}\Big)
= \exp\Big(\sum_{r\geq 1}\frac{(-1)^{r-1}\zeta_{\Gamma}(r)}{r}\,x^{r}\Big).
$$
Comparing the coefficients of $x^{m}$ gives $P_{m} = P_{1}\,h_{m-1}(\zeta_{\Gamma}(1),\dots,\zeta_{\Gamma}(m-1))$.  By
Theorem~\ref{thm:det}, $P_{1} = \prod_{v}d_{v}\,\text{\rm det}^{*}\Delta_{\Gamma}/2^{E_{0}}$.  The statement about
$c_{m}(G_{k})$ is then Theorem~\ref{thm:Phi-limit}.  In particular, part~(2) is Theorem~\ref{mthm:forests}.
\end{proof}

\begin{remark}\label{rem:forest-asymptotics}\label{rem:analytic-in-h}
Parts (1) and (2) say that the leading asymptotics of the weighted number of $m$-spanning forests of $G_{k}$ is governed
by the regularized determinant of $\Gamma$ and by the special values $\zeta_{\Gamma}(1),\dots,\zeta_{\Gamma}(m-1)$, the spanning trees
being the case $m = 1$.  For degenerating families of discrete tori, the corresponding asymptotics of the number of
spanning trees, the complexity, are obtained in \cite{CJK12}, where they are related to the heights of the associated
real tori.  The two parts are consistent, since $k^{-2r}\zeta_{G_{k}}(r)\rightarrow\zeta_{\Gamma}(r)$ by
Corollary~\ref{cor:zeta-degeneration} and the growth of $\text{\rm det}^{*}\Delta_{G_{k}}$ is given by
Proposition~\ref{prop:logdet-asymptotics}.

The error term $o(1)$ cannot be improved for arbitrary subdivisions, since the lower-order coefficients depend on the
offsets $n_{j,k}+1-k\ell_{j}$.  For exact subdivisions it can be expanded to every order, as asserted after
Theorem~\ref{mthm:forests}, and we now outline the proof.  Write $N_{e} = n_{e,k}+1$ for the number of segments of the edge
$e$ in $G_{k}$, and $d = E_{0}-V_{0}+2m$.

Since $E_{k}-V_{k} = E_{0}-V_{0}$, an $m$-component forest of $G_{k}$ omits exactly $E_{0}-V_{0}+m$ segments.  Such forests
therefore fall into finitely many combinatorial types, independent of the $N_{e}$, and within a type they are
parametrized by the positions of the omitted segments.  The weight $\prod_{C}|C|$ is a polynomial in these positions and
in the $N_{e}$, and summing a polynomial over the lattice points of a simplex again gives a polynomial.  Hence
$c_{m}(G_{k}) = Q_{m}(N_{1},\dots,N_{E_{0}})$ for a polynomial $Q_{m}$ with rational coefficients, of degree at most $d$.  Taking
$N_{e} = k\ell_{e}$ for rational lengths and letting $k\rightarrow\infty$, Theorem~\ref{thm:forest-counts} identifies the homogeneous
part of $Q_{m}$ of degree $d$ with $P_{m}$, which is nonzero.  So an exact subdivision gives
$c_{m}(G_{k}) = \sum_{i=0}^{d}Q_{m}^{(d-i)}(\ell)\,k^{d-i}$, a finite expansion whose leading term is that of
Theorem~\ref{thm:forest-counts}.  For $m = 1$ it is the identity $c_{1}(G_{k}) = V_{k}\,\tau(G_{k})$ of Lemma~\ref{lem:trees}.

\end{remark}

\begin{corollary}\label{cor:forest-hierarchy}
For $n\geq 1$ let $N_{n}$ be the polynomial of Theorem~\ref{thm:newton} expressing $\zeta(n)$ through $e_{1},\dots,e_{n}$.
\begin{enumerate}
\item For every finite connected graph $G$, possibly with loops and multiple edges, and every $n\geq 1$,
$$
\zeta_{G}(n) = N_{n}\Big(\frac{c_{2}(G)}{c_{1}(G)},\dots,\frac{c_{n+1}(G)}{c_{1}(G)}\Big),\qquad c_{1}(G) = \text{\rm det}^{*}\Delta_{G}.
$$
For example, $\zeta_{G}(1) = c_{2}/c_{1}$ and $\zeta_{G}(2) = (c_{2}/c_{1})^{2}-2c_{3}/c_{1}$.
\item For every compact metric graph $\Gamma$ and every $n\geq 1$,
$$
\zeta_{\Gamma}(n) = N_{n}\Big(\frac{P_{2}}{P_{1}},\dots,\frac{P_{n+1}}{P_{1}}\Big),\qquad \text{\rm det}^{*}\Delta_{\Gamma} = \frac{2^{E_{0}}}{\prod_{v}d_{v}}\,P_{1}.
$$
\item As $k\rightarrow\infty$, $\zeta_{G_{k}}(n) = k^{2n}\big(\zeta_{\Gamma}(n)+o(1)\big)$ for every $n\geq 1$.
\end{enumerate}
\end{corollary}

\begin{proof}
For (1), the coefficient of $u^{m}$ in $\det(\Delta_{G}+u) = u\prod_{h\geq 1}(u+\lambda_{h})$ is $c_{m}(G) = \text{\rm det}^{*}\Delta_{G}\cdot e_{m-1}$,
where $e_{m-1}$ is the elementary symmetric function of the $\lambda_{h}^{-1}$; see \eqref{eq:forest-zeta}.  In particular
$c_{1} = \text{\rm det}^{*}\Delta_{G}$, so $e_{m-1} = c_{m}/c_{1}$, and Newton's identities give the claim.  Part (2) is
Theorem~\ref{thm:newton} together with Theorem~\ref{thm:det}.  For (3), Theorem~\ref{thm:Phi-limit} gives
$c_{m}(G_{k}) = k^{E_{0}-V_{0}+2m}(P_{m}+o(1))$, hence $c_{m}(G_{k})/c_{1}(G_{k}) = k^{2(m-1)}\big(P_{m}/P_{1}+o(1)\big)$.  The polynomial $N_{n}$ is
isobaric of weight $n$ when $e_{i}$ is given weight $i$, that is, each monomial $e_{1}^{a_{1}}e_{2}^{a_{2}}\cdots$ has $\sum_{i}ia_{i} = n$.  Hence
$N_{n}(k^{2}y_{1},k^{4}y_{2},\dots,k^{2n}y_{n}) = k^{2n}N_{n}(y_{1},\dots,y_{n})$, and (3) follows from (1) and (2).  This is consistent with
Corollary~\ref{cor:zeta-degeneration}, which was obtained from the theta functions.
\end{proof}

\subsection{Computing the special values}\label{subsec:procedure}

Given a finite connected graph $G_{0}$, possibly with loops and multiple edges, positive edge
lengths $\ell_{1},\dots,\ell_{E_{0}}$, and an integer $R\geq 1$, the special values of $\zeta_{\Gamma}$ on
$\Gamma = (G_{0},\Lambda)$ at the arguments $r\leq R$ are obtained as follows.
\begin{enumerate}
\item \emph{Forests.}  List the spanning forests $F$ of $G_{0}$ with at most $R+1$ components.  For
each, record $w_{F}$, $\sigma_{p}(F)$, and for each component $C$ the quantities $\ell^{*}(C)$ and
$\Lambda_{p}(C)$.
\item \emph{Forest invariants.}  Compute $P_{1},\dots,P_{R+1}$ from Theorem~\ref{thm:explicit-P}.
Equivalently, $P_{m}$ is the coefficient of $x^{m}$ in $\prod_{j}S_{j}(x)\det M(x)$.
\item \emph{Elementary symmetric functions.}  Set $e_{m} = P_{m+1}/P_{1}$ for $1\leq m\leq R$.
\item \emph{Positive integers.}  Compute $\zeta_{\Gamma}(r)$ for $1\leq r\leq R$ by
Theorem~\ref{thm:newton}.
\item \emph{Zero and negative integers.}  $\zeta_{\Gamma}(0) = \chi/2-1$ and $\zeta_{\Gamma}(-n) = 0$ for
$n\geq 1$, by Corollary~\ref{cor:zeta-continuation}.
\item \emph{Determinant.}  $\text{\rm det}^{*}\Delta_{\Gamma} = \exp(-\zeta_{\Gamma}'(0)) = 2^{E_{0}}P_{1}/\prod_{v}d_{v}$,
by Theorem~\ref{thm:det}.
\item \emph{Reduced spectrum.}  Subtract the contributions of the loops from $\zeta_{\Gamma}(r)$ to
obtain the power sums over the reduced spectrum, by Theorem~\ref{thm:reduced}.
\end{enumerate}
Every step is a finite computation with rational coefficients in the edge lengths.  This is the
source of Theorem~\ref{thm:field}.

Explicit examples of these computations are given in Section~\ref{sec:examples}.

\subsection{Arithmetic consequences}\label{sec:arithmetic}

The results of Sections~\ref{sec:forest-expansion}--\ref{sec:det} determine every special value of
$\zeta_{\Gamma}$ at the integers, and the regularized determinant, in terms of the edge lengths and
the combinatorics of $G_{0}$.  This has arithmetic consequences, which we now draw.  Their interest
comes from what is known about the spectrum itself, where the results of \cite{KS20, KS23} and
\cite[Chapter~10]{Ku24} are in part conditional on Schanuel's conjecture, but nothing below depends on
that conjecture or on any hypothesis beyond the compactness and connectedness of $\Gamma$.

\begin{thm}\label{thm:field}
Let $K = \mathbb{Q}(\ell_{1},\dots,\ell_{E_{0}})$.  Then the following hold.
\begin{enumerate}
\item each $P_{m}$ is a polynomial in $\ell_{1},\dots,\ell_{E_{0}}$ with rational coefficients;
\item $\zeta_{\Gamma}(r)\in K$ for every integer $r\geq 1$; more precisely, $\zeta_{\Gamma}(r)$ is the
value at $(\ell_{1},\dots,\ell_{E_{0}})$ of a homogeneous rational function of degree $2r$ with rational
coefficients, which depends only on the combinatorial graph $G_{0}$;
\item $\text{\rm det}^{*}\Delta_{\Gamma}\in K$, and all Taylor coefficients of
$\text{\rm det}(\Delta_{\Gamma}+x)$ at $x=0$ lie in $K$;
\item $\zeta_{\Gamma}(0) = \chi/2-1\in\mathbb{Q}$ and $\zeta_{\Gamma}(-n) = 0$ for $n\geq 1$.
\end{enumerate}
\end{thm}

\begin{proof}
Part (1) is Theorem~\ref{thm:explicit-P}, whose coefficients $a_{n}$ and $b_{n}$ are rational.  Part (2) follows from (1) and Corollary~\ref{cor:zeta-integers}.  Part (3) follows from (1) and
Theorem~\ref{thm:det}, and part (4) is Corollary~\ref{cor:zeta-continuation}.
\end{proof}

Theorem~\ref{thm:field} contains Corollary~\ref{cor:arithmetic}, which we can now prove.

\begin{proof}[Proof of Corollary~\ref{cor:arithmetic}]
Parts (2) and (4) of Theorem~\ref{thm:field} give $\zeta_{\Gamma}(n)\in K$ for every integer $n$.  If
every $\ell_{j}$ is algebraic, then $K$ is an algebraic extension of $\mathbb{Q}$, so every element of
$K$ is an algebraic number, which is part (1).  If $\ell_{j} = q_{j}a$ with $q_{j}\in\mathbb{Q}$ for
every $j$, then for $n\geq 1$ Theorem~\ref{thm:field}(2) writes $\zeta_{\Gamma}(n) = R(\ell_{1},\dots,\ell_{E_{0}})$
for a homogeneous rational function $R$ of degree $2n$ with rational coefficients, whence
$$
\zeta_{\Gamma}(n) = R(q_{1}a,\dots,q_{E_{0}}a) = a^{2n}\,R(q_{1},\dots,q_{E_{0}})\in\mathbb{Q}\,a^{2n}.
$$
For $n\leq 0$ the values are $\chi/2-1$ and $0$ by Theorem~\ref{thm:field}(4), which lie in
$\mathbb{Q}\,a^{2n}$ as well.  This is part (2).
\end{proof}

The values $\zeta_{\Gamma}(r)$ are sums over the whole spectrum, part of which is arithmetically
transparent.  A loop of length $\ell$ always carries the frequencies $2\pi n/\ell$ with $n\geq 1$, and
the contribution of such a progression is Euler's.  Let $J\subseteq\{1,\dots,E_{0}\}$ index the loops
of $\Gamma$, and write $\text{\rm Spec}^{*}(\Gamma)$ for the frequencies of $\Gamma$, counted with
multiplicity, from which one copy of the progression $\{2\pi n/\ell_{j}: n\geq 1\}$ has been removed
for each $j\in J$.

\begin{thm}\label{thm:reduced}
Let $\Gamma$ be a compact metric graph.  Then, for every integer $r\geq 1$,
$$
\sum_{\kappa\in\text{\rm Spec}^{*}(\Gamma)}\kappa^{-2r}
= \zeta_{\Gamma}(r)-\sum_{j\in J}\frac{(-1)^{r+1}B_{2r}}{2\,(2r)!}\,\ell_{j}^{2r}.
$$
In particular this sum lies in $K$, and it is an algebraic number whenever the edge lengths are
algebraic.  When $\Gamma$ has no loops the set $J$ is empty, $\text{\rm Spec}^{*}(\Gamma)$ is the whole
set of frequencies, and the identity is the definition of $\zeta_{\Gamma}(r)$.
\end{thm}

\begin{proof}
Let $\mathcal{E}_{j}$ be a loop at the vertex $v$, and for $n\geq 1$ let $f$ be the function equal to
$\sin(2\pi ns/\ell_{j})$ on $\mathcal{E}_{j}$ in the arclength parameter $s$ and equal to zero on every
other edge.  Then $f$ vanishes at every vertex, so it is continuous on $\Gamma$, and
$-f'' = (2\pi n/\ell_{j})^{2}f$ on each edge.  The two half-edges of $\mathcal{E}_{j}$ contribute
$2\pi n/\ell_{j}$ and $-2\pi n/\ell_{j}$ to the sum of the outgoing derivatives at $v$, and every other
half-edge contributes nothing, so $f$ satisfies the Kirchhoff conditions at every vertex.  Hence
$2\pi n/\ell_{j}$ is a frequency of $\Gamma$ for every $n\geq 1$ and every $j\in J$, and the functions
attached to distinct loops have disjoint supports, so one copy of each progression can indeed be
removed and $\text{\rm Spec}^{*}(\Gamma)$ is well defined.  By Euler's formula the progression removed
for $\mathcal{E}_{j}$ accounts for
$$
\sum_{n\geq 1}\Big(\frac{\ell_{j}}{2\pi n}\Big)^{2r}
= \frac{\ell_{j}^{2r}\,\zeta_{\mathbb{Q}}(2r)}{(2\pi)^{2r}}
= \frac{(-1)^{r+1}B_{2r}}{2\,(2r)!}\,\ell_{j}^{2r}
$$
of $\zeta_{\Gamma}(r)$, which gives the identity.  Both terms on the right lie in $K$, the first by
Theorem~\ref{thm:field}(2) and the second because the Bernoulli numbers are rational.
\end{proof}

The focus is therefore on the graphs that do have loops, where an arithmetically transparent part of
the spectrum is removed before the power sums are taken, and on what the remaining frequencies are.
When the edge lengths are
rationally independent and $\Gamma$ is neither a segment, nor a cycle, nor the figure-eight graph,
$\text{\rm Spec}^{*}(\Gamma)$ is the reduced spectrum of \cite[Theorem~10.1 and (10.11)]{Ku24} and
\cite{KS20}, which contains no arithmetic progression and whose $\mathbb{Q}$-linear span is infinite
dimensional, and which is expected on Schanuel's conjecture to be algebraically independent
\cite{KS23}.  In this way Theorem~\ref{thm:reduced} stands to those results as Rayleigh's evaluation
of the sums $\sum_{k}j_{n,k}^{-2r}$ in \eqref{eq:Rayleigh} stands to Siegel's theorem \cite{Si29} on
the transcendence of the zeros $j_{n,k}$.

\section{Examples}\label{sec:families}\label{sec:examples}

In all the examples of this section the function $\Phi$ is computed from the Dirichlet-to-Neumann form
$\Phi = \prod_{j}S_{j}\cdot\det M$ of Proposition~\ref{prop:DtN-properties}(4), and Theorems~\ref{thm:hadamard}
and~\ref{thm:det} and Corollary~\ref{cor:zeta-integers} give the spectrum, the determinant and the
special values of the zeta function.

\subsection{The circle}\label{subsec:circle}

The simplest metric graph is the circle of length $L$.  As a graph, it is a single vertex $v$ with a loop of length $L$, so
$V_{0} = E_{0} = 1$, $\chi = 0$ and $d_{v} = 2$.  Every quantity in Theorems~\ref{mthm:convergence}--\ref{mthm:forests} can be computed in
closed form, and compared with classical formulas.

The only spanning forest is the vertex itself, with one component, and the loop lies
outside it.  The formula of Theorem~\ref{mthm:convergence} therefore has a single term
$$
\Phi(x) = \frac{\sinh(\sqrt{x}\,L)}{\sqrt{x}}\cdot 2\sqrt{x}\tanh\frac{\sqrt{x}\,L}{2} = 4\sinh^{2}\frac{\sqrt{x}\,L}{2} = 2\cosh(\sqrt{x}\,L)-2,
$$
the factor $2$ coming from the two half-edges of the loop.  The zeros of $\Phi$ are $x = -(2\pi j/L)^{2}$, $j\geq 0$, that at $x=0$ simple and the others of multiplicity
two, so the nonzero eigenvalues of the circle are $\mu = (2\pi j/L)^{2}$ with $j\geq 1$, each with multiplicity
two, as expected.

If $L$ is chosen so that the subdivision is exact, the graph $G_{k}$ is the cycle $C_{n}$ with
$n = kL$ vertices, whose Laplacian has the eigenvalues $4\sin^{2}(\pi j/n)$, $0\leq j<n$.  Since $V_{0} = E_{0}$, no normalization is needed in
Theorem~\ref{mthm:convergence}, which states that
$$
\det\big(\Delta_{C_{n}}+x/k^{2}\big) = \prod_{j=0}^{n-1}\Big(\frac{x}{k^{2}}+4\sin^{2}\frac{\pi j}{n}\Big) = 2\cosh(n\beta)-2\;\longrightarrow\;2\cosh(L\sqrt{x})-2,
$$
where $\cosh\beta = 1+x/2k^{2}$, so that $n\beta\rightarrow L\sqrt{x}$.  This is the classical product formula for $\cosh$.  The convergence
of the eigenvalues in Theorem~\ref{mthm:convergence} is equally explicit, since $k^{2}\cdot 4\sin^{2}(\pi j/n)\rightarrow(2\pi j/L)^{2}$ for each fixed $j$.

Since the eigenvalues are $(2\pi j/L)^{2}$ with multiplicity two,
$$
\zeta_{\Gamma}(s) = 2\Big(\frac{L}{2\pi}\Big)^{2s}\zeta_{\mathbb Q}(2s),
$$
where $\zeta_{\mathbb Q}$ is the Riemann zeta function.  Every statement of Theorems~\ref{mthm:det} and~\ref{mthm:special} can be checked directly.  The only pole is at $s = 1/2$, with residue $L/2\pi$.  At $s = 0$, $\zeta_{\Gamma}(0) = 2\zeta_{\mathbb Q}(0) = -1 = \chi/2-1$, and the trivial zeros
of $\zeta_{\mathbb Q}$ give $\zeta_{\Gamma}(-n) = 2(L/2\pi)^{-2n}\zeta_{\mathbb Q}(-2n) = 0$ for $n\geq 1$.  Using $\zeta_{\mathbb Q}(0) = -\frac12$ and $\zeta_{\mathbb Q}'(0) = -\frac12\log 2\pi$, one
finds $\zeta_{\Gamma}'(0) = -2\log L$, that is,
$$
\text{\rm det}^{*}\Delta_{\Gamma} = L^{2},
$$
in agreement with Theorem~\ref{mthm:det}, where $2^{E_{0}}/\prod_{v}d_{v} = 1$ and $L\sum_{T}\prod_{\mathcal{E}_{j}\notin T}\ell_{j} = L\cdot L$.  At the positive integers,
$\zeta_{\Gamma}(n) = 2(L/2\pi)^{2n}\zeta_{\mathbb Q}(2n)$.  Theorem~\ref{mthm:special} computes the left-hand side
without using this identity.  Indeed $\Phi(x) = 2\sum_{m\geq 1}(L^{2}x)^{m}/(2m)!$, so the forest expansion gives
$P_{m} = 2L^{2m}/(2m)!$ and $e_{m} = P_{m+1}/P_{1} = 2L^{2m}/(2m+2)!$, and Newton's recursion returns
$$
\zeta_{\Gamma}(1) = \frac{L^{2}}{12},\qquad \zeta_{\Gamma}(2) = \frac{L^{4}}{720},\qquad
\zeta_{\Gamma}(3) = \frac{L^{6}}{30240},\qquad \zeta_{\Gamma}(4) = \frac{L^{8}}{1209600},
$$
and so on, with no input beyond $\Phi$.  Comparing with the identity above gives
$\zeta_{\mathbb Q}(2n) = (2\pi)^{2n}\zeta_{\Gamma}(n)/(2L^{2n})$, that is, $\zeta_{\mathbb Q}(2) = \pi^{2}/6$,
$\zeta_{\mathbb Q}(4) = \pi^{4}/90$, $\zeta_{\mathbb Q}(6) = \pi^{6}/945$ and in general Euler's evaluation of
$\zeta_{\mathbb Q}(2n)$.  In this sense Theorem~\ref{mthm:special} contains Euler's theorem as the case of a
single loop.  Finally, the convergence $k^{-2s}\zeta_{C_{n}}(s)\rightarrow 2(L/2\pi)^{2s}\zeta_{\mathbb Q}(2s)$ for
$\text{\rm Re}(s)>1/2$ in Theorem~\ref{mthm:convergence} recovers the Riemann zeta function as a limit of the spectral zeta functions
of the discrete cycles.  The spectral zeta functions of the cycles are studied in \cite{FK17} in the critical strip,
where the asymptotic expansion involves both the Riemann zeta function and the zeta function $\zeta_{\mathbb{Z}}$ of the
infinite path.

The weighted numbers of spanning forests of the cycle $C_{n}$ are the coefficients of its
characteristic polynomial, namely
$$
c_{m}(C_{n}) = \frac{n}{m}\binom{n+m-1}{2m-1},\qquad\text{so that}\qquad c_{1}(C_{n}) = n^{2},\quad c_{2}(C_{n}) = \frac{n^{2}(n^{2}-1)}{12}.
$$
As $n = kL\rightarrow\infty$, $c_{m}(C_{n})\sim k^{2m}\cdot 2L^{2m}/(2m)!$, in accordance with Theorem~\ref{mthm:forests}.  Here
$\text{\rm det}^{*}\Delta_{\Gamma} = L^{2}$, so the theorem reads $c_{1}(C_{n})\sim k^{2}L^{2}$ and $c_{2}(C_{n})\sim k^{4}L^{2}\zeta_{\Gamma}(1) = k^{4}L^{4}/12$.  In
general, $h_{m-1}\big(\zeta_{\Gamma}(1),\dots,\zeta_{\Gamma}(m-1)\big) = 2L^{2m-2}/(2m)!$, which is Euler's product
$\sinh z/z = \prod_{n\geq 1}(1+z^{2}/\pi^{2}n^{2})$ written as an identity between forest counts and special values.

\subsection{Flower graphs}

Let $F_{m}$ be the flower graph, a single vertex carrying $m$ loops of lengths
$\ell_{1},\dots,\ell_{m}$, so
that $d_{v} = 2m$, $E_{0} = m$ and $\chi = 1-m$.  The only spanning forest is the single vertex, and
\eqref{eq:Phi} gives
$$
\Phi_{F_{m}}(x) = 2\sqrt{x}\,\sum_{i=1}^{m}\tanh\Big(\frac{\sqrt{x}\,\ell_{i}}{2}\Big)\prod_{i=1}^{m}S_{i}(x),
\qquad P_{1} = L\prod_{i}\ell_{i},\qquad
\text{\rm det}^{*}\Delta_{F_{m}} = \frac{2^{m-1}}{m}\,L\prod_{i}\ell_{i}.
$$
For $m=1$ this is the circle of Section~\ref{subsec:circle}.

If all the lengths are equal to $a$,
$$
\Phi_{F_{m}}(x) = 2m\,(\cosh a\sqrt{x}-1)\Big(\frac{\sinh a\sqrt{x}}{\sqrt{x}}\Big)^{m-1}.
$$
The nonzero eigenvalues are therefore $(2\pi j/a)^{2}$, $j\geq 1$, with multiplicity $m+1$, and
$((2j+1)\pi/a)^{2}$, $j\geq 0$, with multiplicity $m-1$.  For $r\geq 1$,
$$
\zeta_{F_{m}}(r) = \big(2\,(2\pi)^{-2r} + (m-1)\,\pi^{-2r}\big)\,a^{2r}\,\zeta_{\mathbb Q}(2r).
$$

\subsection{Dipole graphs}\label{subsec:dipole}

Let $D_{m}$ consist of two vertices joined by $m$ edges of lengths
$\ell_{1},\dots,\ell_{m}$, with $L = \sum_{i}\ell_{i}$; see Figure~\ref{fig:watermelon-intro}.  Its Dirichlet-to-Neumann matrix is
$\begin{pmatrix}\alpha & -\beta\\ -\beta & \alpha\end{pmatrix}$, with $\alpha = \sum_{i}\sqrt{x}\coth(\sqrt{x}\,\ell_{i})$ and
$\beta = \sum_{i}\sqrt{x}/\sinh(\sqrt{x}\,\ell_{i})$.  Hence $\det M = (\alpha-\beta)(\alpha+\beta)$ and
$$
\Phi_{D_{m}}(x) = x\Big(\sum_{i}\tanh\frac{\sqrt{x}\,\ell_{i}}{2}\Big)\Big(\sum_{i}\coth\frac{\sqrt{x}\,\ell_{i}}{2}\Big)\prod_{i}\frac{\sinh(\sqrt{x}\,\ell_{i})}{\sqrt{x}}.
$$

The spanning trees are the single edges, and the only
two-component forest is the empty one.  So Theorems~\ref{mthm:det} and~\ref{mthm:special} give
$$
\text{\rm det}^{*}\Delta_{D_{m}} = \frac{2^{m}}{m^{2}}\,L\sum_{i}\prod_{j\neq i}\ell_{j},
$$
$$
\zeta_{D_{m}}(1) = \frac{\displaystyle\sum_{i}\prod_{j\neq i}\ell_{j}\Big(\frac{L}{6}\sum_{j\neq i}\ell_{j}^{2}-\frac{1}{12}\sum_{j}\ell_{j}^{3}\Big)
+\frac{L^{2}}{4}\prod_{j}\ell_{j}}{\displaystyle L\sum_{i}\prod_{j\neq i}\ell_{j}}.
$$
More generally, every $\zeta_{D_{m}}(r)$ is read off from the Taylor expansion of $\Phi_{D_{m}}$ at $x = 0$, since by
Theorem~\ref{mthm:convergence},
$$
\sum_{r\geq 1}\frac{(-1)^{r-1}\zeta_{D_{m}}(r)}{r}\,x^{r} = \log\frac{\Phi_{D_{m}}(x)}{P_{1}x},\qquad P_{1} = \Phi_{D_{m}}'(0).
$$

If $\ell_{i} = a$ for all $i$, this gives $\Phi_{D_{m}}(x) = m^{2}x\big(\sinh(a\sqrt{x})/\sqrt{x}\big)^{m}$.  Using
$\sinh z/z = \prod_{n\geq 1}(1+z^{2}/\pi^{2}n^{2})$, the eigenvalues are $(\pi n/a)^{2}$, $n\geq 1$, each with multiplicity $m$,
and
$$
\zeta_{D_{m}}(r) = m\Big(\frac{a}{\pi}\Big)^{2r}\zeta_{\mathbb Q}(2r),\qquad \text{\rm det}^{*}\Delta_{D_{m}} = 2^{m}a^{m},
$$
as in the circle.  For $m=1$ the dipole graph is the interval, with $\Phi(x) = \sqrt{x}\sinh(\ell\sqrt{x})$ and $\text{\rm det}^{*}\Delta = 2\ell$.  For
$m=2$ it is a circle of length $L$ with two marked vertices at arbitrary positions, and $\Phi_{D_{2}}$ is the function $\Phi$ of
the circle in Section~\ref{subsec:circle}, independently of the positions of the vertices, as
Lemma~\ref{lem:subdivision} predicts.

\subsection{Naphthalene}\label{subsec:naphthalene}

\begin{figure}[ht]
\centering
\begin{tikzpicture}[scale=0.95, v/.style={circle,fill,inner sep=1.7pt}, e/.style={thick}]
  \foreach \i in {0,...,5} \coordinate (L\i) at ({30+60*\i}:1.4);
  \foreach \i in {0,...,5} \coordinate (R\i) at ($({30+60*\i}:1.4)+(2.425,0)$);
  \draw[e] (L0) -- (L1) -- (L2) -- (L3) -- (L4) -- (L5) -- cycle;
  \draw[e] (R0) -- (R1) -- (R2) -- (R3) -- (R4) -- (R5) -- cycle;
  \foreach \i in {0,...,5} \node[v] at (L\i) {};
  \foreach \i in {0,1,2,5} \node[v] at (R\i) {};
\end{tikzpicture}
\caption{Naphthalene as a metric graph, two hexagons fused along an edge.  It has $V_{0} = 10$ vertices, $E_{0} = 11$ edges
and $\chi = -1$; the two vertices shared by the hexagons have degree three and the other eight have degree two.}
\label{fig:naphthalene}
\end{figure}
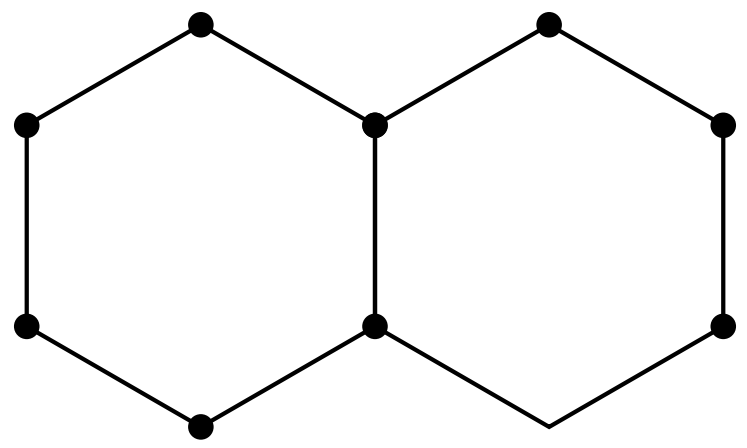

Metric graphs were used to model conjugated molecules long before they were called quantum graphs.  In the
free-electron network model of \cite{RS53}, the $\pi$-electrons of an aromatic hydrocarbon move along the bonds of
the molecule, with the Kirchhoff conditions at the atoms, so that the molecule is exactly a metric graph in the
sense of Section~\ref{sec:metricgraphs}; see \cite{Kuc08} for the later history.  Benzene is then the
circle of Section~\ref{subsec:circle}, of length six times the bond length.  We treat the next case, naphthalene,
which consists of two hexagons fused along an edge; see Figure~\ref{fig:naphthalene}.

Let all bonds have the common length $a$, so that $L = 11a$.  The graph has
$\tau = 35$ spanning trees, and $\prod_{v}d_{v} = 9\cdot 2^{8}$, so Theorem~\ref{mthm:det} gives
$$
\text{\rm det}^{*}\Delta_{\Gamma} = \frac{2^{11}}{9\cdot 2^{8}}\cdot 11a\cdot 35\,a^{2} = \frac{3080}{9}\,a^{3}.
$$

When all edges have a common length the spectrum can be
read off from a finite matrix.  Let $\Delta_{N} = I-D^{-1/2}AD^{-1/2}$ be the normalized Laplacian of
$G_{0}$, with eigenvalues $\lambda_{1},\dots,\lambda_{V_{0}}$, and put $Z(s) = \sum_{\lambda_{i}\neq 0}\lambda_{i}^{-s}$.
By von Below's correspondence \cite{vB85}, see also \cite[Chapter~4]{Ku24} and \cite{HW18}, the eigenvalues
of $\Delta_{\Gamma}$ with $\kappa a\notin\pi\mathbb{Z}$ are the numbers $\kappa^{2}$ with
$\cos(\kappa a) = 1-\lambda_{i}$, with the same multiplicities, while the numbers $(\pi n/a)^{2}$ with $n\geq 1$
are eigenvalues with a multiplicity that depends only on $\beta = E_{0}-V_{0}+1$ and on whether $G_{0}$ is
bipartite.  Expanding $\Phi$ at $x = 0$ and applying
Theorem~\ref{thm:newton}, the first two special values of an equilateral metric graph are
$$
\zeta_{\Gamma}(1) = a^{2}\Big(\frac{2(E_{0}-V_{0})+1}{12}+\frac{Z(1)}{2}\Big),\qquad
\zeta_{\Gamma}(2) = a^{4}\Big(\frac{E_{0}-V_{0}}{90}+\frac{1}{720}-\frac{Z(1)}{12}+\frac{Z(2)}{4}\Big).
$$
For the dipole graph $D_{m}$ of Section~\ref{subsec:dipole} these give $m\,a^{2}/6$ and $m\,a^{4}/90$, in
agreement with the values computed above.

For naphthalene the normalized Laplacian has the
eigenvalues
$$
0,\quad 2,\quad \frac{5\pm\sqrt5}{4},\quad \frac{3\pm\sqrt5}{4},\quad \frac{13\pm\sqrt{37}}{12},\quad
\frac{11\pm\sqrt{37}}{12},
$$
each simple.  The eigenvalues of $\Delta_{\Gamma}$ are therefore the numbers $\kappa^{2}$ with
$\cos(\kappa a) = 1-\lambda_{i}$ and $\kappa a\notin\pi\mathbb{Z}$, together with the numbers
$(\pi n/a)^{2}$ with $n\geq 1$, each of multiplicity $\beta+1 = 3$, since the graph is bipartite with
$\beta = E_{0}-V_{0}+1 = 2$.  The multiplicity at these points replaces the contributions of the
eigenvalues $\lambda_{i} = 0$ and $\lambda_{i} = 2$, which lie there.  One checks that the two families
together have the Weyl density $L/\pi$.

Here $Z(1) = \frac{2157}{154}$ and $Z(2) = \frac{4802957}{118580}$, so that the formulas above give
$$
\zeta_{\Gamma}(1) = \frac{1117}{154}\,a^{2} = 7.2532\ldots a^{2},\qquad
\zeta_{\Gamma}(2) = \frac{531907}{59290}\,a^{4} = 8.9712\ldots a^{4}.
$$
The frequencies involve $\sqrt5$ and $\sqrt{37}$, while the special values are rational multiples of $a^{2r}$, in
accordance with Corollary~\ref{cor:arithmetic}.  The value $\zeta_{\Gamma}(1)$ is the Kirchhoff index of the metric
graph, read as in Section~\ref{sec:numerics}.  With the carbon--carbon bond length $a = 1.4$ \AA, it
equals $14.216\ldots\,\text{\AA}^{2}$ for naphthalene.

\subsection{The Kirchhoff index, and a numerical check on $K_{4}$}\label{sec:numerics}

Our focus here is to compute $\zeta_{\Gamma}(1)$, and the reason to do so is the following.  For a finite graph
$G$ on $V$ vertices, the Kirchhoff index from \cite{KR93} is the sum of the effective resistances between all pairs
of vertices, and it equals $V\,\zeta_{G}(1)$ by a theorem in \cite{GM96}.  The same result holds on a metric
graph in integrated form.  If $\mathcal{G}(x,y)$ is the Green's function of $\Delta_{\Gamma}$ with
$\int_{\Gamma}\mathcal{G}(x,y)\,dy = 0$, and
$R(x,y) = \mathcal{G}(x,x)+\mathcal{G}(y,y)-2\mathcal{G}(x,y)$ is the effective resistance between the points
$x$ and $y$, then
$$
\zeta_{\Gamma}(1) = \int_{\Gamma}\mathcal{G}(x,x)\,dx = \frac{1}{2L}\int_{\Gamma}\int_{\Gamma}R(x,y)\,dx\,dy.
$$
In other words, $\zeta_{\Gamma}(1)$ is the average effective resistance of $\Gamma$, times $L/2$.  For the circle of length $L$ this gives $L^{2}/12$, and for
the interval of length $\ell$ it gives $\ell^{2}/6$.

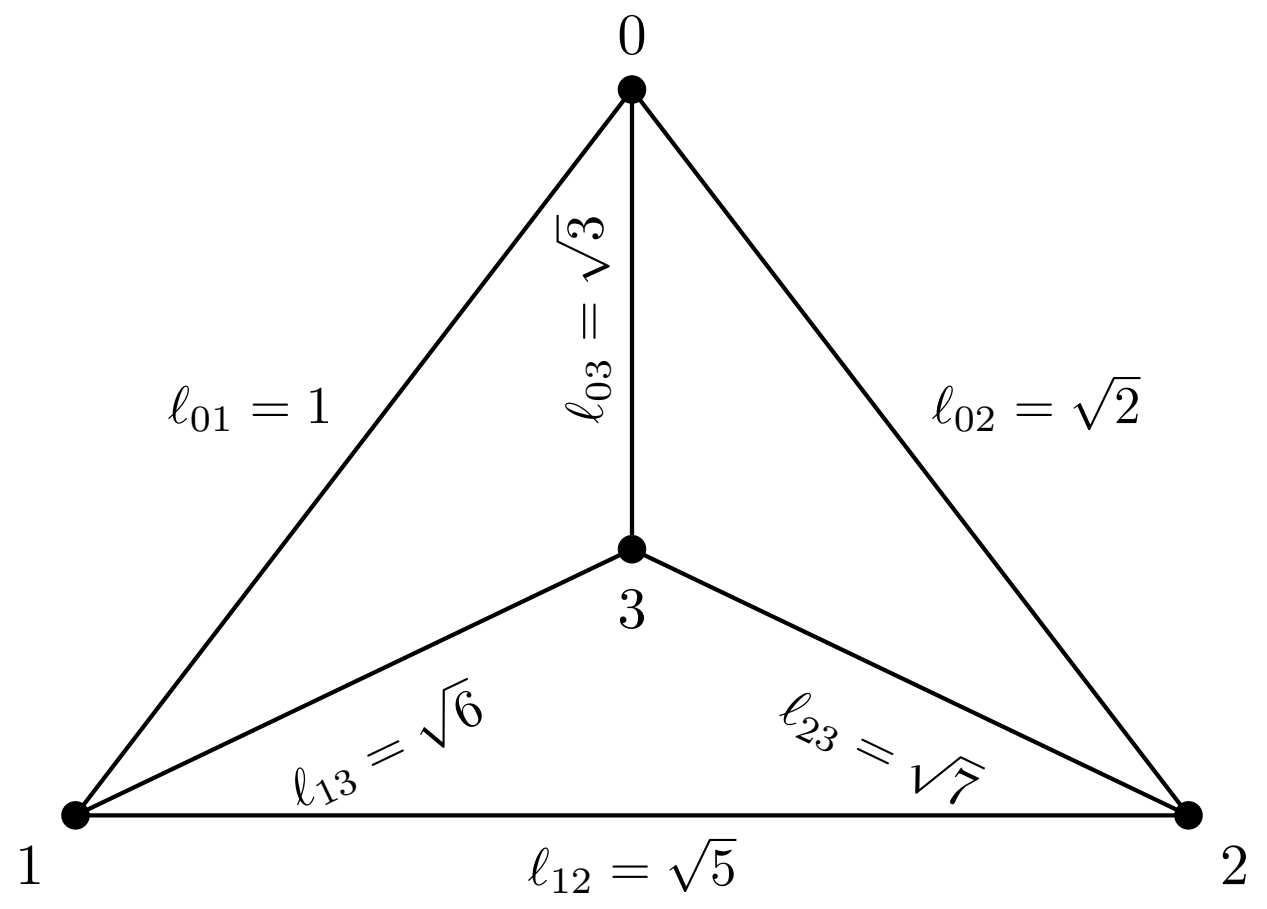
\begin{figure}[ht]
\centering
\begin{tikzpicture}[scale=1.6, v/.style={circle,fill,inner sep=1.9pt}, e/.style={thick}, lab/.style={font=\small}]
  \coordinate (A) at (0,2.0); \coordinate (B) at (-2.3,-1.0); \coordinate (C) at (2.3,-1.0); \coordinate (D) at (0,0.1);
  \draw[e] (A) -- node[lab,above left] {$\ell_{01} = 1$} (B);
  \draw[e] (A) -- node[lab,above right] {$\ell_{02} = \sqrt2$} (C);
  \draw[e] (D) -- node[lab,above,sloped,pos=0.5] {$\ell_{03} = \sqrt3$} (A);
  \draw[e] (B) -- node[lab,below] {$\ell_{12} = \sqrt5$} (C);
  \draw[e] (B) -- node[lab,below=4pt,sloped,pos=0.5] {$\ell_{13} = \sqrt6$} (D);
  \draw[e] (D) -- node[lab,below=4pt,sloped,pos=0.5] {$\ell_{23} = \sqrt7$} (C);
  \node[v,label=above:{$0$}] at (A) {};
  \node[v,label=below left:{$1$}] at (B) {};
  \node[v,label=below right:{$2$}] at (C) {};
  \node[v,label={[label distance=1pt]below:{$3$}}] at (D) {};
\end{tikzpicture}
\caption{The metric graph $K_{4}$  with edge
lengths $1,\sqrt2,\sqrt3,\sqrt5,\sqrt6,\sqrt7$.}
\label{fig:K4-numerics}
\end{figure}

Consider $K_{4}$, the complete graph with four vertices, which we label as $0,1,2,3$.  Assign the edge lengths
$$
\ell_{01} = 1,\quad \ell_{02} = \sqrt2,\quad \ell_{03} = \sqrt3,\quad \ell_{12} = \sqrt5,\quad \ell_{13} = \sqrt6,\quad \ell_{23} = \sqrt7,
$$
which are rationally independent, and with total length $L = 1+\sqrt2+\sqrt3+\sqrt5+\sqrt6+\sqrt7$; see Figure~\ref{fig:K4-numerics}.

We compute the special values $\zeta_{K_{4}}(n)$ for $n = 1,\dots,5$ from Theorems~\ref{mthm:convergence}
and~\ref{mthm:special}, and compare them with the power sums of numerically computed eigenvalues.

We computed the Taylor coefficients $P_{1},\dots,P_{6}$ of
$\Phi = \prod_{j}S_{j}\det M$ to $50$ digits, and obtained $\zeta(1),\dots,\zeta(5)$ from Theorem~\ref{thm:newton}.
As we have shown, there are closed algebraic expressions for $\zeta(1),\dots,\zeta(5)$, although they can be
unwieldy.  The first coefficient, however, is given
exactly by Theorem~\ref{thm:explicit-P} as
\begin{align*}
P_{1} = {}& 18+42\sqrt2+50\sqrt3+36\sqrt5+30\sqrt6+34\sqrt7+19\sqrt{10}+19\sqrt{14}\\
&+22\sqrt{15}+16\sqrt{21}+15\sqrt{30}+10\sqrt{35}+11\sqrt{42}+12\sqrt{70}+10\sqrt{105}+6\sqrt{210},
\end{align*}
which agrees with the numerically computed coefficient to all $25$ digits compared.  Since every vertex has degree $3$ and $E_{0} = 6$,
Theorem~\ref{mthm:det} gives
$$
\text{\rm det}^{*}\Delta = \frac{2^{6}}{3^{4}}\,P_{1}.
$$
By Corollary~\ref{cor:arithmetic}, all the values $\zeta(n)$ lie in $\mathbb{Q}(\sqrt2,\sqrt3,\sqrt5,\sqrt7)$.
We computed the eigenvalues $\mu = \kappa^{2}$ with $\kappa < K = 1000$ from the secular determinant
\eqref{eq:secular-det} of Remark~\ref{rem:euler-product}, which involves only the scattering amplitudes and the
edge lengths and so is independent of the forest expansion.  Putting $\sqrt{x} = i\kappa$ there, the function
$$
F(\kappa) = e^{i\kappa L}\,(i\kappa)^{\chi}\,\det\big(I-U(\kappa)\big) = \frac{2^{E_{0}}}{\prod_{v}d_{v}}\,\Phi(-\kappa^{2})
$$
is real for real $\kappa$, and its positive zeros are the frequencies $k_{j} = \sqrt{\mu_{j}}$, with
multiplicities, by Theorem~\ref{thm:zeros}.  In floating point the imaginary part of the computed $F(\kappa)$ is
of relative size $10^{-14}$, which checks the implementation.  We evaluated $F$ on a uniform grid of step
$2\cdot10^{-4}$ in $(0,K]$ and refined each sign change by Brent's method to an accuracy of $10^{-14}$ in
$\kappa$.  Such a search can miss only a zero of even multiplicity, which does not change sign and does not
occur for generic edge lengths \cite{Fr05}, or two zeros inside one grid interval.  Against the second, we found
$3651$ eigenvalues where the Weyl term gives $LK/\pi\approx 3653.4$, and the smallest gap between consecutive
zeros is $0.0043$, about twenty times the grid step.  Comparing $\zeta(n)$ with the power sum $\sum_{k_{j}<K}\mu_{j}^{-n}$, corrected by
the Weyl estimate $\frac{L}{\pi}\frac{K^{1-2n}}{2n-1}$ of the tail, the two agree to $4\cdot10^{-7}$ for
$n = 1$, where the tail estimate limits the accuracy, to $4\cdot10^{-13}$ for $n = 2$, and to the double
precision of the computed eigenvalues for $3\leq n\leq 5$.  The agreement confirms three results at once, the
identification of the zeros of $\Phi$ with the spectrum, with multiplicities (Theorem~\ref{thm:zeros}), Newton's
identities for the special values (Theorem~\ref{thm:newton}), and the identity \eqref{eq:secular-det}, since the
two sides were evaluated independently.

\section*{Acknowledgments}
This project began when the first-named author (J.J.) visited the second-named author (A.K.) at his home
university in Geneva in December 2023, and it continued in earnest when both of them met the third-named author (L.S.) in
Sarajevo in February 2026.  We are all grateful to both institutions for providing the opportunity to collaborate on this
research.  This article developed from numerous visits by J.J. to the home institutions of A.K. and L.S., and J.J. is very grateful for their hospitality.  All authors thank Gautam Chinta for mathematical conversations and prior collaborations.

A large language model (Anthropic) was used as an assistive software tool for numerical investigation, Python/SageMath and Pari/GP script generation, and logical error checking. All mathematical results were devised and obtained 
by the authors, who assume full responsibility for all manuscript content.

\enlargethispage{3\baselineskip}
\vspace{1mm}\noindent
\begin{minipage}[t]{0.32\textwidth}
\raggedright\footnotesize
Jay Jorgenson\\
Department of Mathematics\\
The City College of New York\\
Convent Avenue at 138th Street\\
New York, NY 10031\\
U.S.A.\\
jjorgenson@mindspring.com
\end{minipage}\hfill
\begin{minipage}[t]{0.32\textwidth}
\raggedright\footnotesize
Anders Karlsson\\
Section de math\'ematiques\\
Universit\'e de Gen\`eve\\
2-4 Rue du Li\'evre\\
Case Postale 64, 1211\\
Gen\`eve 4, Suisse\\
anders.karlsson@unige.ch\\[1ex]
and\\[1ex]
Matematiska institutionen\\
Uppsala universitet\\
Box 256, 751 05\\
Uppsala, Sweden\\
anders.karlsson@math.uu.se
\end{minipage}\hfill
\begin{minipage}[t]{0.32\textwidth}
\raggedright\footnotesize
Lejla Smajlovi\'c\\
Department of Mathematics\\
University of Sarajevo\\
Zmaja od Bosne 35\\
71 000 Sarajevo\\
Bosnia and Herzegovina\\
lejlas@pmf.unsa.ba
\end{minipage}

\end{document}